\def\VersionDateTime{14/March/2017, 16:35 GMT+9:00. Version $2.1$}
\documentclass[12pt]{amsart}
\usepackage{amssymb,amsthm,amsmath,amscd,rsfs}
\usepackage{latexsym}

\newcommand{\NN}{{\mathbb{N}}}
\newcommand{\ZZ}{{\mathbb{Z}}}
\newcommand{\QQ}{{\mathbb{Q}}}
\newcommand{\RR}{{\mathbb{R}}}
\newcommand{\CC}{{\mathbb{C}}}
\newcommand{\PP}{{\mathbb{P}}}

\newcommand{\KK}{{\mathbb{K}}}
\newcommand{\GG}{{\mathbb{G}}}

\newcommand{\bfK}{\mathbf{K}}

\newcommand{\bfX}{\mathbf{X}}

\newcommand{\OO}{{\mathcal{O}}}

\newcommand{\codim}{\operatorname{codim}}
\newcommand{\ch}{\operatorname{char}}

\newcommand{\Aut}{\operatorname{Aut}}

\newcommand{\Ker}{\operatorname{Ker}}

\newcommand{\Image}{\operatorname{Image}}

\newcommand{\Supp}{\operatorname{Supp}}

\newcommand{\Spec}{\operatorname{Spec}}
\newcommand{\Spf}{\operatorname{Spf}}
\newcommand{\Gal}{\operatorname{Gal}}

\newcommand{\cherncl}{{c}}

\newcommand{\Proof}{{\sl Proof.}\quad}
\newcommand{\QED}{{\unskip\nobreak\hfil\penalty50\quad\null\nobreak\hfil
{$\Box$}\parfillskip0pt\finalhyphendemerits0\par\medskip}}

\newcommand{\ord}{\operatorname{ord}}

\newcommand{\Tr}{\operatorname{Tr}}

\newcommand{\aff}{\mathrm{aff}}

\newcommand{\an}{\mathrm{an}}
\newcommand{\trop}{\mathrm{trop}}
\newcommand{\val}{\operatorname{val}}

\newcommand{\str}{\operatorname{str}}

\newcommand{\relin}{\operatorname{relin}}

\newcommand{\ndr}{\operatorname{nd-rk}}
\newcommand{\pr}{\operatorname{pr}}

\newcommand{\Irr}{\operatorname{\mathrm{Irr}}}

\title
{Survey on the geometric Bogomolov conjecture}
\author
{Kazuhiko Yamaki}
\date{\VersionDateTime}
\subjclass[2000]{Primary~14G40, Secondary~11G50.}
\address
{Institute for Liberal Arts and Sciences,
Kyoto University, Kyoto, 606-8501, Japan}
\email{yamaki.kazuhiko.6r@kyoto-u.ac.jp}

\keywords{Geometric Bogomolov conjecture, Bogomolov conjecture,
canonical heights, small points}

\begin{document}

\theoremstyle{plain}
\newtheorem{Theorem}{Theorem}[section]
\newtheorem{Lemma}[Theorem]{Lemma}
\newtheorem{Proposition}[Theorem]{Proposition}
\newtheorem{Corollary}[Theorem]{Corollary}
\newtheorem{Main-Theorem}[Theorem]{Main Theorem}
\newtheorem{Theorem-Definition}[Theorem]{Theorem-Definition}
\theoremstyle{definition}
\newtheorem{Definition}[Theorem]{Definition}
\newtheorem{Remark}[Theorem]{Remark}
\newtheorem{Conjecture}[Theorem]{Conjecture}
\newtheorem{Claim}{Claim}
\newtheorem{Example}[Theorem]{Example}
\newtheorem{Key Fact}[Theorem]{Key Fact}
\newtheorem{ack}{Acknowledgments}       \renewcommand{\theack}{}
\newtheorem*{n-c}{Notation and convention}      
\newtheorem{citeTheorem}[Theorem]{Theorem}
\newtheorem{citeProposition}[Theorem]{Proposition}

\newtheorem{Step}{Step}

\renewcommand{\theTheorem}{\arabic{section}.\arabic{Theorem}}
\renewcommand{\theClaim}{\arabic{section}.\arabic{Theorem}.\arabic{Claim}}
\renewcommand{\theequation}{\arabic{section}.\arabic{Theorem}.\arabic{Claim}}

\def\Pf{\trivlist\item[\hskip\labelsep\textit{Proof.}]}
\def\endPf{\strut\hfill\framebox(6,6){}\endtrivlist}

\def\Pfo{\trivlist\item[\hskip\labelsep\textit{Proof of Proposition~\ref{ch-of-hyp}.}]}
\def\endPfo{\strut\hfill\framebox(6,6){}\endtrivlist}

\maketitle


\begin{abstract}
This is a survey paper of the developments on the geometric Bogomolov
conjecture.
We explain the recent results by the author
as well as
previous works concerning the conjecture.
This paper also includes an introduction to the 
height theory over function fields
and a quick review on basic notions on non-archimedean analytic geometry.
\end{abstract}

\section{Introduction} \label{sect:intro}

The Bogomolov conjecture is a
problem in Diophantine geometry;
it is a conjecture which should
characterize
the closed subvarieties
with a dense set of points of small height.
There are several versions of the Bogomolov conjecture.
The conjectures for curves
and for abelian varieties
are
widely studied problems among others,
and the former 
is a special case of the latter.

The Bogomolov conjectures are considered both over number fields
and over function fields.
When we consider the conjecture over function fields
with respect to classical heights, we call 
it
the \emph{geometric} Bogomolov conjecture.
When we do that over number fields,
we sometimes call it the arithmetic Bogomolov conjecture.

The arithmetic
Bogomolov conjectures for curves and for abelian varieties
have been already established as
theorems; that for curves is due to Ullmo
\cite{ullmo},
and that for abelian varieties is due to Zhang
\cite{zhang2}.
On the other hand,
the geometric Bogomolov conjecture for abelian varieties is still open,
but there are some significant partial works,
such as \cite{gubler2, yamaki5,yamaki6, yamaki7, yamaki8}.
Recently, using those results concerning
the conjecture for abelian varieties,
we have shown in \cite{yamaki8} that the conjecture for curves holds
in full generality.

The purpose of this survey is to
give an exposition of 
the recent developments of the geometric Bogomolov conjecture for abelian 
varieties
and to
describe 
the idea of the proof of the geometric Bogomolov
conjecture for curves.

\subsection{Notation and conventions}

We put together the notation and conventions that will be
used throughout this article.

\subsubsection{Sets}
A \emph{natural number} means a rational integer greater than $0$.
Let $\NN$ denote the set of natural numbers.

By the convention of this paper,
the notation
$A \subset B$ allows the case of $A=B$.

\subsubsection{Varieties}

Let $\mathfrak{K}$ be any field.
A variety over $\mathfrak{K}$ means
a geometrically integral scheme that is separated and
of finite type
over $\mathfrak{K}$.

For a variety $X$ over $\mathfrak{K}$,
let $\mathfrak{K} (X)$ denote the function field of $X$.

\subsubsection{Abelian varieties} \label{subsubsect:av}
Let $A$ be an abelian variety 
over any field $\mathfrak{K}$.
For each $n \in \ZZ$, let $[n] : A \to A$ denote the
$n$-times endomorphism on $A$, i.e.,
the morphism given by $x \mapsto nx$.

Let $L$ be a line bundle on $A$.
We say that $L$ is \emph{even}
if $[-1]^{\ast} (L) \cong L$.
This property is called symmetry in some literatures.
We remark that for any line bundle $L$ on $A$,
$L \otimes [-1]^{\ast} (L)$ is even.

\subsubsection{Models}

Let $R$ be any ring.
Let $X$ be a scheme over $R$.
Let $S$ be a scheme with a dominant morphism $\Spec (R) \to S$.
A \emph{model over $S$ of $X$ }
is
a morphism of schemes 
$\mathcal{X} \to S$ equipped with an isomorphism $\mathcal{X}
\times_S \Spec ( R ) \cong X$.
A model is said to be \emph{projective} (resp. \emph{proper}, \emph{flat})
if the morphism $\mathcal{X}
\to S$ is projective (resp. proper, flat).
Furthermore, let $L$ be a line bundle on $X$.
A \emph{model over $S$ of $(X , L)$} is a pair of a model $\mathcal{X}$
over $S$
of $X$
and
a line bundle $\mathcal{L}$ on 
$\mathcal{X}$ with an identification $\mathcal{L}|_{X} = L$,
where $\mathcal{L}|_{X}$ is the pull-back of $\mathcal{L}$ by the
morphism $X \to \mathcal{X}$ induced from 
the given identification
$\mathcal{X}
\times_S \Spec ( R ) \cong X$.

Let $R_0$ be a subring of $R$.
We say that \emph{$X$ can be defined over $R_0$}
if there exists a model over $R_0$ of $X$.

\subsubsection{Ground fields} \label{subsubsect:NCgroundfields}

Let $k$ be an algebraically closed field.
Let $\mathfrak{B}$ be a projective normal variety over $k$
with $\dim (\mathfrak{B}) \geq 1$.
When $\dim (\mathfrak{B}) \geq 2$, we fix an ample line
bundle
$\mathcal{H}$ on $\mathfrak{B}$,
which is needed to define height in this case.

Let $K$ denote the function field of $\mathfrak{B}$
or a number field.
We fix an algebraic closure
$\overline{K}$ of $K$.
Any finite extension of $K$ should be taken in $\overline{K}$.

\subsection{Background and history}

The geometric Bogomolov conjecture is not an isolated topic;
there are many related topics and previous results.
Let us begin with its background together with brief historical notes.

\subsubsection{Manin--Mumford conjecture}
According to Lang \cite{lang0},
Manin and Mumford conjectured around 1963 (independently) the following.
Let $\mathfrak{K}$ be an algebraically closed field.
Assume that $\ch (\mathfrak{K}) = 0$.
Let $C$ be a smooth projective curve over $\mathfrak{K}$ of genus $g \geq 2$.
Fix a divisor $D$ on $C$ of degree $1$
and let $\jmath_D : C \to J_C$ be the embedding 
of $C$ into its Jacobian $J_C$
given by $\jmath_D (x) = x-D$.
Then the set
\[
\{ x \in C ( \mathfrak{K} ) \mid \jmath_D (x) \in J_C( \mathfrak{K})_{tor} \}
\]
is finite,
where 
$J_C( \mathfrak{K})_{tor}$ is the set of torsion points of $J_C( \mathfrak{K})$.
We call this conjecture the \emph{Manin--Mumford conjecture}
for curves
over $\mathfrak{K}$.

\subsubsection{Bogomolov conjecture for curves} \label{subsubsect:BC}
Inspired by the conjecture by Manin and Mumford, 
Bogomolov proposed in 1980 an arithmetic analogue of
the conjecture (cf. \cite{bogomolov}),
which is now called the \emph{Bogomolov conjecture for curves}.
Let $K$ be a number field
or a function field.
Fix an algebraic closure $\overline{K}$ of $K$.
Let $C$ be a smooth projective curve
over $\overline{K}$ of genus
$g \geq 2$.
Let $h_{NT}$ be the N\'eron--Tate height on $J_C$.
Remark that $h_{NT}$
is a positive semidefinite quadratic form on the additive group 
$J_C ( \overline{K} )$.
For each real number $\epsilon$,
set
\[
C ( \epsilon ) :=
\{ x \in C ( \mathfrak{K} ) \mid h_{NT}
(\jmath_D (x)) \leq \epsilon \}
\]
When $K$ is a function field,
assume
that
$C$ is non-isotrivial.
Then the Bogomolov conjecture for $C$ asserts that
there should exist an $\epsilon > 0$ such that
$C ( \epsilon )$ is finite.

The Bogomolov conjecture for curves over number fields
generalizes
the Manin--Mumford conjecture for curves
over $\overline{\QQ}$.
Indeed,
since $h_{NT}$ is a quadratic form, we have 
$h_{NT} ( \tau ) = 0$ for any $\tau \in J_C ( \overline{\QQ})_{tor}$.
It follows that
\[
\{ x \in C ( \overline{\QQ} ) \mid \jmath_D (x) \in J_C(
 \overline{\QQ})_{tor} \}
\subset C ( \epsilon )
\]
for any $\epsilon > 0$.
This shows that the Manin--Mumford conjecture over $\overline{\QQ}$
is
deduced from the
Bogomolov conjecture over number fields.

\subsubsection{Raynaud's theorem} \label{subsubsect:raynaud}

In 1983, 
Raynaud proved that the Manin--Mumford conjecture
holds.
In fact, he proved the following.

\begin{Theorem} [cf. \cite{raynaud1}] \label{thm:raynaud-curve}
Let $\mathfrak{K}$ be an algebraically closed field.
Assume that $\ch ( \mathfrak{K}) = 0$.
Let $A$ be an abelian variety over $\mathfrak{K}$.
Let $X$ be a closed subvariety of $A$ of dimension $1$.
Suppose that $X$ is not a torsion subvariety.
Then $X \cap A ( \overline{K} )_{tor}$ is a finite set.
\end{Theorem}

Here, \emph{torsion subvariety} is by definition
the translate of an abelian subvariety by a torsion point.
The above theorem generalizes the Manin--Mumford conjecture
for curves,
because a smooth projective curve of genus at least $2$
embedded in its Jacobian is not a torsion subvariety.

In the same year,
Raynaud generalized Theorem~\ref{thm:raynaud-curve}
to obtain the following theorem,
which is often also called the Manin--Mumford conjecture.

\begin{Theorem} [Raynaud's theorem, cf. \cite{raynaud2}] \label{thm:raynaud-AV}
Let $\mathfrak{K}$ be an algebraically closed field.
Assume that $\ch ( \mathfrak{K}) = 0$.
Let $A$ be an abelian variety over $\mathfrak{K}$.
Let $X$ be a closed subvariety of $A$.
Suppose that $X$ is not a torsion subvariety.
Then $X \cap A ( \overline{K} )_{tor}$ is 
not Zariski dense in $X$.
\end{Theorem}

\subsubsection{Admissible paring and some results} \label{subsubsect:adimsiiblepairing}

While 
the Manin--Mumford conjecture
had been established in a generalized form,
the Bogomolov conjecture for curves
was still open.

Let $K$ be a number field or the function field
of $\mathfrak{B}$ with $\dim ( \mathfrak{B}) = 1$.
In \cite{zhang1},
Zhang defined an \emph{admissible pairing} 
$(\cdot , \cdot )_a$ on a curve
and
introduced the admissible dualizing sheaf $\omega_a$.
Further,
he proved that
the Bogomolov conjecture for a curve $C$ amounts to 
the positivity of the admissible self-pairing 
$( \omega_a , \omega_a )_a$
of $\omega_a$ on $C$.
In that paper, this admissible pairing is described
as 
the self-intersection of the relative dualizing sheaf
of the stable model of the curve
minus a positive number 
that is a combinatorial data arising from
the dual graphs of the singular fibers of the
stable model.

By using Zhang's admissible pairing,
some partial affirmative answers to the Bogomolov conjecture for curves
over function fields
(of transcendence degree $1$) have been obtained, such as
\cite{zhang1, moriwaki0, moriwaki3, moriwaki4, moriwakiRBI, yamaki1, yamaki4}.
However, the conjecture was not proved in full generality.

\subsubsection{Ullmo's theorem and Zhang's theorem}

In 1998,
Ullmo finally proved that the Bogomolov conjecture over number fields
holds.

\begin{Theorem}
Assume that $K$ is a number field.
Let $C$ and $C(\epsilon)$ be as in \S~\ref{subsubsect:BC}.
Then there exists an $\epsilon > 0$
such that $C ( \epsilon )$ is finite.
\end{Theorem}

The proof is \emph{not} given by showing the positivity of the
admissible pairing.
It is obtained by a different approach,
namely, the equidistribution of small points due to Szpiro, Ullmo,
and Zhang \cite{SUZ}.

Just after
Ullmo established the above theorem, 
Zhang proved the following generalized version
of the Bogomolov conjecture over number fields.
To state that, we fix a notation.
Let $A$ be an abelian variety over $\overline{K}$
and
let $L$ be an even ample line bundle on $A$
(cf. \S~\ref{subsubsect:av}).
Let $\widehat{h}_L$ be the canonical height 
on $A$ associated to $L$;
see \S~\ref{sect:canonicalheights}.
It is a semipositive definite quadratic form on $A ( \overline{K} )$.
Let $X$ be a closed subvariety of $A$.
For any real number $\epsilon$,
set $X ( \epsilon ; L)$ to be
$
\{ x \in X ( \overline{K} )
\mid
\widehat{h}_L ( x ) \leq \epsilon
\}
$.

\begin{Theorem} [Corollary~3 of \cite{zhang2}] \label{thm:zhangstheoremintro}
Assume that $K$ is a number field.
Let $A$, $L$, and $X$ be as above.
Suppose that $X ( \epsilon ; L)$ is dense in $X$ for any $\epsilon >0$.
Then $X$ is a torsion subvariety.
\end{Theorem}

Zhang's theorem (Theorem~\ref{thm:zhangstheoremintro})
generalizes Ullmo's theorem.
Indeed, let $C$ be a smooth projective curve over $\overline{K}$.
With the notation in \S~\ref{subsubsect:BC},
consider the case where $A = J_C$ and $X = \jmath_D ( C )$.
Then $h_{NT}$ is the canonical height associated
to a symmetric theta divisor
(cf. Remark~\ref{rem:NTheight}),
and thus
one easily sees that Ullmo's theorem is
a part of Zhang's theorem.
We remark that Zhang's theorem generalizes Raynaud's theorem over
$\overline{\QQ}$.

The proof uses equidistribution theory,
like the proof of Ullmo's theorem.
We will give an outline of the proof of Zhang in the sequel.

\subsubsection{Moriwaki's generalization and the Manin--Mumford conjecture}

In 2000,
Moriwaki generalized Zhang's theorem.
Let $F$ be a finitely generated field
over $\QQ$.
Moriwaki constructed in
\cite{moriwaki5}
arithmetic height functions over $F$.
There, a pair of an
arithmetic variety 
with function field $F$
and an arithmetic line bundle on the arithmetic variety
is called
a polarization of $F$,
and the polarization is said to be big
if the arithmetic line bundle is big.
Once a polarization of $F$ is fixed,
Moriwaki's arithmetic height over $F$ 
is defined.
He proved 
basic properties of this height,
and he constructed the canonical height on
abelian varieties.
Furthermore, he established the following theorem.

\begin{Theorem} [Theorem~B of \cite{moriwaki5}]
\label{thm:moriwaki}
Let $F$ be a finitely generated field over $\QQ$.
Let $A$ be an abelian variety over $\overline{F}$
and let $L$ be an even ample line bundle on $A$.
Fix a big polarization on $F$,
and let $\widehat{h}_L^{arith}$ be the canonical height
on $A$
associated to $L$.
Let $X$ be a closed subvariety of $A$.
Suppose that 
for any $\epsilon > 0$,
the set
\[
\left\{
x \in 
X ( \overline{F} )
\mid
\widehat{h}_L^{arith} (x) \leq \epsilon
\right\}
\]
is dense in $X$.
Then $X$ is a torsion subvariety.
\end{Theorem}

Moriwaki's theorem 
(Theorem~\ref{thm:moriwaki})
generalizes Raynaud's theorem
over any algebraically closed field of characteristic $0$;
we omit the detail here but will see an analogous argument in
\S~\ref{sect:MM}.

We notice that the
arithmetic height 
arising from a big polarization,
which is used in Moriwaki's theorem above, is different from the 
classical ``geometric'' heights over function field.
Therefore, 
the geometric Bogomolov conjecture,
which will be formulated in the sequel,
addresses a problem that is not the same as in Moriwaki's theorem.

\subsubsection{Gubler's theorem}

Let $K$ be 
any function field
(cf. \S~\ref{subsubsect:NCgroundfields}).
In 2007, Gubler established a result over a function field 
which is analogous to Zhang's theorem,
under the assumption that the abelian variety 
is totally degenerate at some
place.
Let $A$ be an abelian variety over $\overline{K}$.
We say that $A$ is \emph{totally degenerate at some place}
if there exists 
a codimension one point $v$ of $\mathfrak{B}$
satisfying the following condition:
there exist a discrete valuation ring $R' \subset \overline{K}$ 
that dominates $\OO_{\mathfrak{B},v}$
and a 
group scheme 
$\mathscr{A}^{\circ} \to \Spec (R')$ with 
geometric generic fiber $A$
such that
the spacial fiber $\widetilde{\mathscr{A}^{\circ}}$
is an algebraic torus;
see also \S~\ref{subsect:degeneration} for details
on degeneracy of abelian varieties.

\begin{Theorem} [Theorem~1.1 of \cite{gubler2}]\label{intthm:gubler'stheorem}
Let $A$ be an abelian variety over $\overline{K}$
and
let $L$ an even ample line bundle on $A$.
Assume that $A$ is totally degenerate at some place of $K$.
Let $X$ be a closed subvariety of $A$.
Suppose that $X ( \epsilon ; L)$ is dense in $X$
for any $\epsilon > 0$.
Then $X$ is a torsion subvariety.
\end{Theorem}

The proof is given by a non-archimedean analogue of
the the proof of Ullmo and Zhang,
which
is a quite important argument.
We will explain it
later in detail.

\subsubsection{Cinkir's theorem} \label{subsubsect:cinkir}

In 2011, 
the geometric Bogomolov conjecture for curves
was proved by Cinkir
under the assumption that $\ch (k) = 0$ and that
$K$ is the function field
of a curve, i.e., $\dim (\mathfrak{B}) = 1$;
see \cite[Theorem~2.13]{cinkir}.
We remark that Cinkir's theorem is effective in the sense that
he gave explicitly a positive number $r$ such that
$C ( r )$ is
finite when $C$ has a semistable model over $\mathfrak{B}$.

To explain Cinkir's proof,
let $K$ be the function field of a curve 
and let $C$ be a smooth projective curve
over $\overline{K}$ of genus $g \geq 2$.
In 2010, 
Zhang proved a new description of the admissible pairing of
the admissible dualizing sheaf
(cf. \cite{zhang3}).
This description uses the height of the Gross--Schoen cycle.
More precisely,
he showed that the admissible pairing of the 
admissible dualizing sheaf equals the sum
of the height of the Gross--Schoen cycle and the
``$\varphi$-invariants'' of the dual graphs of the semistable model of 
the curve.
It is known that if $\ch (k) = 0$,
 the height of the Gross--Schoen cycle are non-negative,
so that 
the Bogomolov conjecture amounts to the positivity
of $\varphi$-invariants.
In 2009, using Zhang's work, Faber proved 
in \cite{faber}
the geometric Bogomolov conjecture
for curves over $\overline{K}$ of small genus.
In 2011, after contributions on the study of graph invariants
by many authors, Cinkir proved 
in \cite{cinkir}
that
the $\varphi$-invariants for non-trivial graphs are positive.
Since the case where every reduction graph is trivial,
that is, the case of everywhere good reduction case had been known,
this proved the geometric Bogomolov conjecture for curves over such 
$\overline{K}$.

However, that is not the finial answer to the conjecture.
Cinkir's theorem needs the assumption that $\dim ( \mathfrak{B} )
= 1$ and $\ch (k) = 0$.
The first assumption is needed because Zhang's description
in \cite{zhang3}
of the admissible pairing is obtained under this assumption.
The second assumption on the characteristic is more crucial.
If $\ch (k) = 0$,
the positivity of the Gross--Schoen cycle follows from the Hodge
index theorem,
but if that is not the case,
the Hodge index theorem,
which is part of the standard conjectures,
is not known.
Therefore the positivity of the $\varphi$-invariants
is not enough for the positivity of the admissible pairing,
and thus 
the Bogomolov conjecture in positive characteristic
cannot be deduced in the same way.

\subsection{Geometric Bogomolov conjecture for abelian varieties} \label{subsect:initroGBC}

It is natural to ask whether or not
the same statement as Theorem~\ref{thm:zhangstheoremintro}
holds for any abelian variety over function fields,
but in fact,
 it does not holds in general.
For example, suppose that $B$ is a constant abelian variety, that is,
$B = \widetilde{B} \otimes_k \overline{K}$ for
some abelian variety $\widetilde{B} $ over $k$.
Let $Y$ be a constant closed subvariety of $B$,
that is, $Y := \widetilde{Y} \otimes_k \overline{K}$ for some
closed subvariety $\widetilde{Y} \subset \widetilde{B}$.
Then
$\widetilde{Y} (k)$, which is naturally a 
subset of $Y ( \overline{K} )$,
is dense in $Y$, and 
for any point $y$ of this set, we have $\widehat{h}_{M} (y) = 0$,
where $M$ is an even ample line bundle on $B$
(cf. Example~\ref{ex:height0points}).
Furthermore, if $\phi : B \to A$ is a homomorphism of abelian varieties,
then we see that $\phi ( \widetilde{Y}(k) )$ is a dense subset
of $\phi (Y)$
and that for any $x \in \phi ( \widetilde{Y}(k) )$
we have $\widehat{h}_{L} (x) = 0$,
where $L$ is  an even ample line bundle
on $A$
(cf. Remark~\ref{rem:imageofconstantpoints}).
This suggests that an abelian variety over a function field
in general has a non-torsion closed subvariety with 
a dense set of height $0$ points.

While we cannot expect the same statements as Zhang's theorem
holds over function fields in general,
it is still natural to ask
how can we characterize the closed subvarieties $X$
such that $X (\epsilon ; L)$ is dense in $X$ for any $\epsilon > 0$.
In 2013, we proposed 
in \cite{yamaki5} a conjecture
that
the ``special subvarieties'' should be the only 
such closed subvarieties.

The special subvarieties are the subvarieties which 
are the sum of a torsion subvariety and a
closed subvariety which is the image of a constant closed subvariety
of a constant abelian variety.
To be precise,
let $A$ be an abelian variety over $\overline{K}$.
Let $X$ be a closed subvariety of $A$.
We say that $X$ is \emph{special}
if there exist a torsion subvariety $T$ of $A$,
a constant abelian variety $B= \widetilde{B} \otimes_k \overline{K}$,
a constant closed subvariety $Y = \widetilde{Y} \otimes_k \overline{K}$
of $B$,
and a homomorphism $\phi : B \to A$
such that $X = \phi (Y) + T$.
Remark that in this definition
of special subvarieties, we take any constant abelian variety $B$,
but it is enough to consider
the universal one among homomorphisms
from constant abelian varieties to $A$, which is called the
$\overline{K}/k$ trace; refer to \S~\ref{subsect:statementGBC} for details.

It is not difficult to see that if $X$ is a
special subvariety,
then for any $\epsilon > 0 $,
$X(\epsilon ; L)$ is dense in $X$.
(cf. Remark~\ref{rem:specialhasdense0}).
The \emph{geometric Bogomolov conjecture for abelian varieties}
asserts that the converse should also hold.

\begin{Conjecture} [cf.
Conjecture~0.3 of \cite{yamaki5}] \label{conj:introGBC}
Let $A$ be an abelian variety over $\overline{K}$,
where $K$ is a function field.
Let $L$ be an even ample line bundle on $A$.
Let $X$ be a closed subvariety of $A$.
Suppose that for any $\epsilon > 0$,
$X ( \epsilon ; L)$ is dense in $X$.
Then $X$ is a special subvariety.
\end{Conjecture}

The geometric Bogomolov conjecture for abelian varieties is still open, 
but the study of this conjecture has been developed,
and there are significant partial answers.
Gubler's theorem is of course an important one.
In \cite{yamaki5, yamaki6, yamaki7,yamaki8,yamaki9},
we have generalized Gubler's theorem,
seeing that the conjecture
holds for a large class of abelian varieties.
This is the main topic of this paper, and
we will explain it in detail
(cf. Theorems~\ref{thm:yamaki7}, \ref{thm:GBCyamaki8},
and \ref{thm:yamaki9}).

By virtue of the development of the study of the geometric Bogomolov
conjecture for 
abelian varieties,
we have very recently proved that the geometric Bogomolov
conjecture for curves holds in full generality.
More generally, we have shown the following.

\begin{Theorem} [cf. Theorem~1.3 of \cite{yamaki8}] \label{thm:GBC_dim1}
Let $A$ be an abelian variety over $\overline{K}$
and let $L$ be an even ample line bundle on $A$.
Let $X$ be a closed subvariety of $A$.
Assume that $\dim (X) = 1$.
Suppose that $X(\epsilon ; L)$ is dense in $X$ for any $\epsilon > 0$.
Then 
$X$ is a special subvariety.
\end{Theorem}

The geometric Bogomolov conjecture for curves 
is deduced
from the above theorem.
Indeed, in the setting of the conjecture for curves,
put $X := \jmath_D ( C )$.
Then $\dim (X) =1$, and since $C$ is non-isotrivial,
$X$ is non-special.
Remark that the N\'eron--Tate height on the Jacobian is the
canonical height associated to the theta divisor.
By Theorem~\ref{thm:GBC_dim1},
$C ( \epsilon )$ is not dense for some $\epsilon > 0$.
Thus we obtain the conjecture for curves.

Our result on the geometric Bogomolov conjecture
for curves is not effective 
in contrast to Cinkir's theorem
(when
$K$ is a function field of transcendence degree $1$
over $k$ of characteristic $0$
and $C$ has semistable reduction over $K$). 
On the other hand,
we have an advantage in  working over any function field
of any characteristic.
One benefit
is that we have an application to the
Manin--Mumford conjecture in positive characteristic
(cf. \S~\ref{sect:MM}).

We give a remark on the proof of Theorem~\ref{thm:GBC_dim1}.
To prove this theorem,
we first show 
that Conjecture~\ref{conj:introGBC} holds under the assumption 
of $\codim ( X , A) = 1$ (cf. Theorem~\ref{thm:GBC_codimX=1})
by using the recent partial results on the Conjecture~\ref{conj:introGBC},
and then we reduce Theorem~\ref{thm:GBC_dim1} to this 
codimension $1$ case.
It should be remarked that our proof of the 
geometric Bogomolov conjecture for curves 
works via 
partial results
of the conjecture for abelian varieties.

\subsection{Organization}
This paper consists of twelve sections
including the introduction, with an appendix.
In \S~\ref{sect:canonicalheights},
we recall the notion of canonical heights on abelian varieties,
where
we mainly focus on the heights over function fields.
In \S~\ref{sect:GBC}, we formulate the geometric
Bogomolov conjecture for abelian varieties,
and we
present some partial results on the conjecture.
In a part of the argument,
we use an non-archimedean analogue of the proofs
of Ullmo's and Zhang's theorem.
Therefore, in \S~\ref{sect:zhang'sproof},
we recall the idea of Zhang's proof,
and in
\S~\ref{sect:nonarchimedeangeom},
we recall
the basic ideas of non-archimedean geometry.
In \S~\ref{sect:proofofgubler},
we explain the proof of Gubler's theorem.
This theorem is the first one where
the non-archimedean analogue of Zhang's proof worked well,
and we will use the idea of the proof of this theorem.
In \S~\ref{sect:structuretheorem},
we describe the structure of the canonical measures.
This structure theorem plays a crucial role
in our argument to reduce the geometric Bogomolov conjecture
to that for nowhere degenerate abelian varieties.
\S~\ref{sect:reductiontoNDcase}
is the first main part of this paper.
There,
we prove that the geometric Bogomolov conjecture
for any abelian variety is reduced to the conjecture for
its maximal nowhere degenerate abelian subvariety.
By this result,
in particular, the conjecture is reduced to that for abelian
varieties that is nowhere degenerate.
In \S~\ref{sect:heightofsubvarietiels},
we recall the notion of canonical heights of closed subvarieties
of an abelian variety.
In \S~\ref{sect:proofofGBCforcurves},
we prove that the geometric Bogomolov conjecture
for nowhere degenerate abelian varieties are reduced to the conjecture
for those with trivial $\overline{K}/k$-trace.
\S~\ref{sect:proofforcurve}
is the second main part of this paper,
where
we give an outline of the proof of Theorem~\ref{thm:GBC_dim1}.
In the last section \S~\ref{sect:MM},
we give a remark of the Manin--Mumford conjecture in positive characteristic.
In the appendix,
we give a summary of some ideas in non-archimedean geometry.

\subsection*{Acknowledgments}
This survey paper is the proceeding 
of my talk at the conference
``Nonarchimedean analytic Geometry: Theory and Practice''
held at
Papeete from 24 to 28, August, 2015.
I thank the organizers for inviting me to the conference
and giving me an opportunity to give a talk.
I thank Professor J\'er\^ome Poineau
for recommending me to write this survey.
Further, I thank the referee for reading the manuscript carefully
and giving me many helpful comments.
This work was partly supported by the Japan Society for the 
Promotion of Science through KAKENHI 26800012.

\section{Canonical heights on abelian varieties} \label{sect:canonicalheights}

In this section, we recall basic properties of canonical heights
on abelian varieties.
For the theory of heights,
\cite[Chapter~1--Chapter~6]{lang2} is a basic reference;
see also \cite{BG}.

\subsection{Height over function fields}
Before describing the notion of canonical heights on abelian varieties,
let us recall the idea of heights.
The theory of heights is developed over
over 
number fields and over function fields.
Because our main topics concern the
\emph{geometric} Bogomolov conjecture,
we mainly focus on
the heights over function fields
in this subsection.
At the end of this subsection,
we will give a comment on the height over number fields.

Let $K$ be the function field of $\mathfrak{B}$
(cf. \S~\ref{subsubsect:NCgroundfields}).
We begin with the notion of height
arising from a model.
Let $X$ be a projective variety over $\overline{K}$
and let $L$ be a line bundle on $X$.
Let $K'$ be a finite extension of $K$
and let $\mathfrak{B}'$ be the normalization of $\mathfrak{B}$
in $K'$.
Let $( \mathcal{X} , \mathcal{L})$ be a
model of $(X , L)$ over $\mathfrak{B}'$
satisfying the following conditions:
the morphism
$\pi : \mathcal{X} \to \mathfrak{B'}$
is proper; 
and
there exists an open subset $\mathfrak{U} \subset \mathfrak{B}'$
with $\codim ( \mathfrak{B}' \setminus \mathfrak{U} , \mathfrak{B}') \geq 2$
over which $\pi : \mathcal{X} \to \mathfrak{B'}$ is flat.
When $\dim (\mathfrak{B}) = 1$,
the last condition is equivalent to saying that $\pi$ is flat.
There always exists such a model.
When we say a model in this subsection,
we assume that it satisfies the above conditions.

Let $\mathcal{H}'$ be the pull-back of $\mathcal{H}$ by
the finite morphism $\mathfrak{B}' \to \mathfrak{B}$.
We take a point $x \in X( \overline{K} )$. 
Let $\Delta_x$ be the closure of $x$ in $\mathscr{X}$
and let $K'(x)$ be the function field of $\Delta_x$.
Then we set
\begin{align*}
h_{(\mathcal{X} , \mathcal{L})} (x)
:=
\frac{\deg_{\mathcal{H}'} \pi_{\ast} \left(
\cherncl_1 ( \mathcal{L}) \cdot [ \Delta_x ]
\right)}{[K'(x) : K]}
\end{align*}
where
$\cherncl_1 ( \mathcal{L}) \cdot [ \Delta_x ]$ is a cycle class on 
$\mathscr{X}$,
$\pi_{\ast} \left(
\cherncl_1 ( \mathcal{L}) \cdot [ \Delta_x ]
\right)$ is the pushout by $\pi$,
and
$\deg_{\mathcal{H}'}$ means the degree with respect to $\mathcal{H}'$.
This defines a function 
\[
h_{(\mathcal{X} , \mathcal{L})} : 
X ( \overline{K} ) \to \RR,
\]
called the \emph{height (function) arising from a
model $(\mathcal{X},\mathcal{L})$}.

The following lemma immediately follows from the definition.

\begin{Lemma} \label{lem:modelheight-functoriality}
Let $X$ and $Y$ be projective varieties over $\overline{K}$.
Let $f : Y \to X$ be a morphism of varieties.
Let $L$ and $L'$ be line bundles on $X$.
Let $\mathfrak{B}'$ be the normalization of $\mathfrak{B}$ 
in a finite extension of $K$
and let $(\mathcal{X} , \mathcal{L})$ and $(\mathcal{X} , \mathcal{L}')$ be
models over $\mathfrak{B}'$ of $(X,L)$ and $(X , L')$, respectively.
Let $\mathcal{Y}$ be a model over $\mathfrak{B}'$ of $Y$
and let $\varphi : \mathcal{Y} \to \mathcal{X}$ be a morphism
over $\mathfrak{B}$
which extends $f$.
Then we have the following:
\begin{enumerate}
\item
$h_{(\mathcal{X} , \mathcal{L} \otimes \mathcal{L}')} =
h_{(\mathcal{X} , \mathcal{L} )}
+
h_{(\mathcal{X} , \mathcal{L}' )}
$;
\item
$f^{\ast} h_{(\mathcal{X} , \mathcal{L})} =
h_{(\mathcal{Y} , \varphi^{\ast} \mathcal{L})}$;
\item
if $\mathcal{L}$ 
is relatively ample, then there exists a constant $C \in \RR$ such that
$h_{(\mathcal{X} , \mathcal{L})} \geq C$.
\end{enumerate}
\end{Lemma}

Indeed, (1) follows from the linearity 
of the intersection product with respect to the line bundle,
and (2) follows from the projection formula.
To see (3),
first note that for any line bundle $\mathcal{M}$ on $\mathfrak{B}'$,
we have
\addtocounter{Claim}{1}
\begin{align} \label{eq:heightlowerbounded}
h_{(\mathcal{X}, \pi^{\ast} (\mathcal{M}))} (x)
=
\frac{\deg_{\mathcal{H}'} \pi_{\ast} \left(
\cherncl_1 ( \pi^{\ast} (\mathcal{M})) \cdot [ \Delta_x ]
\right)}{[K'(x) : K]}
=
\frac{\deg_{\mathcal{H}'} ( \mathcal{M} )}{[K' : K]}
.
\end{align}
Suppose that $\mathcal{L}$ is relatively ample.
Then there exists a line bundle $\mathcal{M}$ on 
$\mathfrak{B}'$ such that $\mathcal{L} \otimes \pi^{\ast} (\mathcal{M})$
is semiample on $\mathcal{X}$,
that is, $\left( \mathcal{L} \otimes \pi^{\ast} (\mathcal{M}) \right)^{\otimes a}$
is basepoint free for some positive integer $a$.
Since a semiample line bundle is nef,
we obtain $h_{(\mathcal{X} , \mathcal{L} \otimes \pi^{\ast} (\mathcal{M}))} \geq 0$.
It follows that
\[
h_{(\mathcal{X} , \mathcal{L})}
=
h_{(\mathcal{X} , \mathcal{L} \otimes \pi^{\ast} (\mathcal{M}))} 
+
h_{(\mathcal{X} , \pi^{\ast} (\mathcal{M}^{\otimes -1}))}
\geq
\frac{\deg_{\mathcal{H}'} ( \mathcal{M}^{\otimes (-1)} )}{[K' : K]}
,
\]
where we use (\ref{eq:heightlowerbounded}).
Thus we have
(3).

For a pair $(X, L)$ of a projective variety $X$ over $K$
and a line bundle $L$ on $X$,
one would hope, naively, to define the height function $h_L$ on $X$ associated
to $L$ by using a model.
In fact,
such a model is not unique,
so that the function arising from a model depends on the choice of
a model.
Indeed, 
in the argument to show Lemma~\ref{lem:modelheight-functoriality}~(3),
both $( \mathcal{X} , \mathcal{L})$
and $( \mathcal{X} , \mathcal{L} \otimes \pi^{\ast} ( \mathcal{M}))$
are models of $(X,L)$,
and $h_{(\mathcal{X} , \mathcal{L})} \neq h_{( \mathcal{X} , \mathcal{L} \otimes \pi^{\ast} ( \mathcal{M}))}$ in general.

However, the following lemma
indicates that the notion of heights with respect to line bundles
can be defined as functions up to bounded functions:

\begin{Lemma} \label{lem:difference_bounded}
Let $X$ be a projective variety over $\overline{K}$
and
let $L$ be line a bundle on $X$.
Let $\mathfrak{B}_1'$ and 
$\mathfrak{B}_2'$ be normal varieties over $k$ 
finite over $\mathfrak{B}$,
and let $(\mathcal{X}_1 , \mathcal{L}_1)$ and 
$( \mathcal{X}_2 , \mathcal{L}_2 )$ be 
proper flat models over $\mathfrak{B}'_1$ and $\mathfrak{B}_2'$,
respectively.
Then $h_{(\mathcal{X}_1 , \mathcal{L}_1)} -
h_{( \mathcal{X}_2 , \mathcal{L}_2 )}$
is a bounded function on $X ( \overline{K} )$.
\end{Lemma}

Taking into account Lemma~\ref{lem:difference_bounded},
we define a height associated to a line bundle $L$
as a representative of an equivalence class of functions on 
$X ( \overline{K} )$ modulo bounded functions.
To be precise,
let $h_L : X ( \overline{K} ) \to \RR$ be a function.
We call $h_L$ a \emph{height (function) associated to $L$}
if there exist a normal variety $\mathfrak{B}'$ over $k$
finite over $\mathfrak{B}$ and a 
proper model $(\mathcal{X} , \mathcal{L})$
of $(X , L)$ over $\mathfrak{B}$
such that
$\mathcal{X} \to \mathfrak{B}'$
is flat
over an open subset
whose complement has codimension 
in 
$\mathfrak{B}'$ at least $2$
and 
such that $h_{L} - h_{(\mathcal{X} , \mathcal{L})}$
is a bounded function over $X(\overline{K})$.
Note that
by Lemma~\ref{lem:difference_bounded},
this condition is equivalent to say that 
$h_{L} - h_{(\mathcal{X} , \mathcal{L})}$ is bounded for any 
normal variety $\mathfrak{B}'$ over $k$
finite over $\mathfrak{B}$ and any 
proper model $(\mathcal{X} , \mathcal{L})$
over $\mathfrak{B}'$ as above.
In particular, the height arising from a model is a height function.

In the following, 
when we write $h_L$,
it means a height function associated to $L$.
By definition, it
is not a uniquely determined function from $L$, but it is unique
up to a bounded function.

For any real-valued functions $h_1$ and $h_2$ defined over the same set,
let $h_1 \sim h_2$ mean that $h_1 - h_2$ is a bounded function.
The following proposition is an immediate consequence
of Lemma~\ref{lem:modelheight-functoriality}.

\begin{Proposition} \label{prop:heightproperty}
Let $X$ and $Y$ be projective varieties over $\overline{K}$,
and let $L$ and $L'$ be line bundles on $X$.
\begin{enumerate}
\item
We have $h_{L \otimes L'} \sim h_{L} + h_{L'}$.
\item
Let $f : Y \to X$ be a morphism of varieties over $\overline{K}$.
Then $h_{f^{\ast} L} \sim f^{\ast} h_{L}$.
\item
Suppose that $L$ is ample.
Then $h_{L}$ is positive up to a bounded function,
that is,
there exists a real constant $C$ such that $h_{L} \geq C$.
\end{enumerate}
\end{Proposition}

\Proof
Assertions (1) and (2) follows from Lemma~\ref{lem:modelheight-functoriality}.
To see (3), suppose that $L$ is ample.
There exist a finite extension $K'$,
a projective variety $X'$ over $K'$,
and a line bundle $L'$ on $X'$ such that $X = X' \otimes_{K'} K$
and $L = L' \otimes_{K'} K$.
Since $L$ is ample, so is $L'$, and hence there exists a positive
integer $N$ such that $(L')^{\otimes N}$ is very ample.
Let $\jmath : X' \hookrightarrow \PP^{N'}_{K'}$ be the closed embedding
associated to global sections of $(L')^{\otimes N}$.
Let $\mathfrak{B}'$ be the normalization of 
$\mathfrak{B}$ in $K'$.
Noting that $\PP^{N'}_{K'}$ is the generic fiber of the canonical projection
$\PP^{N'} \times \mathfrak{B}' \to \mathfrak{B}'$,
we take the closure $\mathcal{X}$ 
of $\jmath (X')$ in $\PP^{N'} \times \mathfrak{B}'$.
Then $\mathcal{X}$ is a proper model of $X$ over $\mathfrak{B}'$,
and it is flat over any point of codimension $1$.
Further, let $\mathcal{L}$ be the restriction of the tautological
line bundle on $\PP^{N'} \times \mathfrak{B}'$ to $\mathcal{X}$.
Then $(\mathcal{X} , \mathcal{L})$ is a model of $(X,L)$,
and $\mathcal{L}$ is relatively ample.
It follows from Lemma~\ref{lem:modelheight-functoriality}~(3)
and Lemma~\ref{lem:difference_bounded}
that $h_L$ is bounded below.
Thus (3) holds.
\QED

We remark that
there is another approach to the height theory,
which begin with the Weil height;
see \cite[Chapter~3]{lang2} for the detail.
The notion of heights we describe here and 
that in \cite[Chapter~3]{lang2} are same
(cf. \cite[Chapter~4]{lang2}).

Also over number fields,
one can define the notion of heights.
In introducing the notion of heights
over number fields,
one can define the heights arising from
models by  using the arithmetic intersection theory;
see \cite{KMY} for example.
Another way, which may be more standard way,
is to use the Weil heights;
see \cite[Chapter~3]{lang2} for details.
We just 
remark that
the heights over number fields 
also satisfy the properties of
Proposition~\ref{prop:heightproperty}.
Further, the notion of canonical height,
which is given in the next subsection,
can be well defined not only over function fields but also
over number fields.

\subsection{Canonical heights on abelian varieties}
For a given line bundle,
a height function 
associated to the line bundle
are only determined up to a bounded function.
However, in some cases, we can make a canonical choice
of height functions among them.
In fact,
we have a notion of canonical height over
abelian varieties,
as we are going to explain.

In this subsection, let $K$ be a function field or a number field.
Let $A$ be an abelian variety over $\overline{K}$
and let $L$ be an even line bundle on $A$.
Fix a positive integer $n$.
Since we have 
$[n]^{\ast} (L) \cong L^{\otimes n^2}$
by the theorem of cube (cf. \cite[\S~6, Corollary~3]{mumford}),
Proposition~\ref{prop:heightproperty} gives us
\addtocounter{Claim}{1}
\begin{align} \label{eq:ate}
[n]^{\ast} h_L \sim n^2 h_{L}
,
\end{align}
where $\sim$ means that they are equal up to bounded functions,
namely,
an equality modulo bounded functions.
The following proposition shows that
there exists a unique height associated to $L$
with which (\ref{eq:ate}) is actually an equality:

\begin{Proposition} \label{prop:ch}
Fix an integer $n > 1$.
Then
there exists a unique height function $h_L$ such that
$[n]^{\ast} h_L = n^2 
h_{L}$.
\end{Proposition}

We can construct such an
$h_L$ as in the above proposition
by limiting process.
Fix an integer $n \geq 2$.
Let $h_1$ be any height function associated to $L$.
We define a sequence $(h_m)_{m \in \NN}$ of function on $A ( \overline{K})$
inductively by
$h_{m+1} := \frac{1}{n^{2}} [n]^{\ast} h_m$.
Then each $h_m$ is a height function on $A$ associated to $L$.
Furthermore, one can show that this sequence
converges to a height function $h_L$ associated to $L$.
By this construction, one also checks $[n]^{\ast} h_L = n^2 
h_{L}$.

The uniqueness is shown as follows.
Suppose that $h_L'$ is a height function associated to $L$
such that $[n]^{\ast} h_L' = n^2 
h_{L}'$.
Set $f:= h_L - h_L'$.
Since $h_L$ and $h_L'$ are height functions associated
to $L$, $f$ is a bounded function.
We prove $f = 0$ by contradiction.
Suppose that there exists an
$x \in X ( \overline{K} )$ such that $f(x) \neq 0$.
Then 
\[
 f (nx) = [n]^{\ast} h_L (x) - [n]^{\ast} 
h_L' (x) = n^{2} h_L (x) - n^{2} h_L' (x)
= n^{2} f (x)
.
\]
It follows inductively that for any positive integer $m$,
we have $f (n^m x) = n^{2m}f (x)$.
Taking $m \to \infty$, we see that this equality indicates that 
$f$ is not bounded,
and that is contradiction.

More strongly than Proposition~\ref{prop:ch},
the following theorem
holds.

\begin{Theorem} \label{thm:canonicalheight}
There exists a unique bilinear form $b_L : A ( \overline{K} )
\times A ( \overline{K} ) \to \RR$ such that the function
$\widehat{h}_L :  A ( \overline{K} ) \to \RR
$ defined by $\widehat{h}_L (x) = \frac{1}{2} b_{L} ( x , x)$ 
is a height function associated to $L$.
\end{Theorem}

The above theorem shows that there exists a height function associated
to $L$ that is a quadratic form.
Since
$[n]^{\ast} \widehat{h}_L = n^{2} \widehat{h}_L$ for any $n \in \ZZ$,
the uniqueness in
Theorem~\ref{thm:canonicalheight} follows from the
uniqueness assertion in Proposition~\ref{prop:ch}.
Remark that
by the uniqueness,
Theorem~\ref{thm:canonicalheight} indicates that
the height function in Proposition~\ref{prop:ch}
is actually a quadratic form
and the heights in Proposition~\ref{prop:ch} and 
Theorem~\ref{thm:canonicalheight} are the same.
For the proof of Theorem~\ref{thm:canonicalheight},
we refer to \cite[Chapter~5]{lang2}.

The height $\widehat{h}_L$ in 
Theorem~\ref{thm:canonicalheight}
is called the
\emph{canonical height
associated to $L$}.
For an even line bundle
$L$ on an abelian variety, let $\widehat{h}_L$
always
denote the canonical height associated to $L$.

\begin{Proposition} \label{prop:propertycanonicalheight}
For the canonical heights 
associated to even line bundles on abelian varieties, the following hold.
\begin{enumerate}
\item
Let $L_1$ and $L_2$ be even line bundles on an abelian variety.
Then
$\widehat{h}_{L_{1} \otimes L_2} = \widehat{h}_{L_1} + \widehat{h}_{L_2}$.
\item
Let $\phi : B \to A$ be a homomorphism of abelian varieties
and let $L$ be an even line bundle on $A$.
Then
$\widehat{h}_{\phi^{\ast} (L)} = \phi^{\ast} \widehat{h}_{L}$.
\item
Let $L$ be an even line bundle on an abelian variety.
Suppose that $L$ is ample.
Then $\widehat{h}_{L} \geq 0$.
\end{enumerate}
\end{Proposition}

\Proof
By Proposition~\ref{prop:heightproperty},
all the assertions hold up to bounded functions.
By the uniqueness assertion of Theorem~\ref{thm:canonicalheight},
we see that (1) and (2) are really equalities,
and (3) also follows from the fact that the canonical height
is a quadratic form
(cf. 
Theorem~\ref{thm:canonicalheight}).
\QED

Next,
we recall basic facts on points of height $0$.

\begin{Remark} \label{rem:torsion_height0}
Let $L$  be an even line bundle on an abelian variety $A$.
Since $\widehat{h}_{L}$ is a quadratic form,
we have $\widehat{h}_{L} (a) = 0$ for any torsion point $a$
of $A ( \overline{K} )$.
\end{Remark}

\begin{Lemma} \label{lem:height0points}
Let $L_1$ and $L_2$ be even ample line bundles on $A$.
Let $a \in A ( \overline{K} )$.
Then $\widehat{h}_{L_1} (a) = 0$ if and only if
$\widehat{h}_{L_2} (a) = 0$.
\end{Lemma}

\Proof
It suffices to show that $\widehat{h}_{L_1} (a) = 0$
implies $\widehat{h}_{L_2} (a) = 0$.
Suppose that $\widehat{h}_{L_1} (a) = 0$.
Since $L_1$ is ample, there exists a positive integer $N$ such that
$L_2^{\otimes -1} \otimes L_1^{\otimes N}$ is ample.
By Proposition~\ref{prop:propertycanonicalheight}~(3), we have
$
\widehat{h}_{L_2^{\otimes -1} \otimes L_1^{\otimes N}} \geq 0
$,
and by Proposition~\ref{prop:propertycanonicalheight}~(1),
we obtain $\widehat{h}_{L_2} \leq N \widehat{h}_{L_1}$.
It follows that $\widehat{h}_{L_2} (a) \leq N \widehat{h}_{L_1} (a) = 0$.
Since $L_2$ is ample,
 $\widehat{h}_{L_2} (a) \geq 0$.
Thus we conclude $\widehat{h}_{L_2} (a) = 0$.
\QED

We say that a point $a \in A ( \overline{K} )$ \emph{has height $0$}
if $\widehat{h}_{L} (a) = 0$, where $L$ is an even ample line bundle
on $A$.
This is well defined by
Lemma~\ref{lem:height0points}.

\begin{Lemma} \label{lem:imageofheith0points}
Let $\phi : B \to A$ be a homomorphism of abelian varieties over $\overline{K}$.
Let $b \in  B ( \overline{K} )$ be a point of height $0$.
Then $\phi (b)$ has height $0$.
\end{Lemma}

\Proof
Let $L$ be an even ample line bundle on $A$
and let $M$ be an even ample line bundle on $B$.
Replacing $M$ with $M^{\otimes a}$ for some positive integer $a$ if necessary,
we may and do assume that
$M \otimes \phi^{\ast} ( L ^{\otimes -1})$ is ample as well as even.
By Proposition~\ref{prop:propertycanonicalheight},
we have
\[
0 \leq
\widehat{h}_L (\phi (b))
= 
\widehat{h}_{\phi^{\ast} (L)} ( b )
\leq
\widehat{h}_{M} ( b ) = 0
,
\]
and thus $\phi (a)$ has height $0$.
\QED

We end with a remark on the N\'eron--Tate height on the Jacobian varieties.

\begin{Remark} \label{rem:NTheight}
Let $J_C^{(g-1)}$ be the Jacobian variety of degree $g-1$ divisor
classes on $C$.
Let $\Theta$ be the theta divisor on $J_C^{(g-1)}$,
that is,
the image of $C^{g-1} \to J_C^{(g-1)}$ given by $(p_1 , \ldots , p_{g-1})
\mapsto p_1 + \cdots + p_{g-1}$.
Let $c_0$ be a divisor on $C$ of degree $1$
such that $(2g - 2) c_0$ is a canonical divisor on $C$.
Let $\lambda : J_C \to J_{C}^{(g-1)}$ be the isomorphism defined
by $a \mapsto a + (g-1) c_0$.
Let $\theta$ be the pullback of $\Theta$ by $\lambda$.
It is known that $\theta$ is ample.
Since $(2g-2)c_0$ is a canonical divisor on $C$,
the line bundle $L := \OO_{J_C} ( \theta )$ is even.
This $\theta$ is called a symmetric theta divisor.
Then
the N\'eron--Tate height $h_{NT}$
is defined to be $\widehat{h}_L$.
We refer to \cite[Chaper~5, \S~5]{lang2} for details.
\end{Remark}

\subsection{Canonical height on the generic fiber of abelian scheme} \label{subsect:gfofabsch}
In this subsection,
we
assume that $K$ is a function field.
Let $A$ be an abelian variety over $\overline{K}$.
Assume that there exist 
a finite covering $\mathfrak{B}' \to \mathfrak{B}$
with $\mathfrak{B}'$ normal
and an abelian scheme $\pi : \mathcal{A} \to \mathfrak{B}'$ 
whose geometric generic fiber equals $A$.
Under this assumption,
we describe the canonical height in terms of intersection 
products on $\mathcal{A}$.
Using that description,
we furthermore see that in general,
a constant abelian variety has a lot of points
which are non-torsion but of canonical height $0$.
The description of the canonical heights
given here will be generalized in \S~\ref{subsect:modelforNDAV}.

Let $L$ be an even ample line bundle on $A$.
Then replacing $\mathfrak{B}'$ by a further finite covering,
we may assume that
there exists a line bundle $\mathcal{L}$ on $\mathcal{A}$
whose restriction to $A$ coincides with $L$.
By tensoring the pull-back by $\pi$ of a line bundle on $\mathfrak{B}'$,
we may take $\mathcal{L}$ such that $0_{\pi}^{\ast} ( \mathcal{L})$ is trivial,
where $0_{\pi}$ is the zero-section of the abelian scheme 
$\pi : \mathcal{A} \to \mathfrak{B}'$.

Then we have $[n]^{\ast} (\mathcal{L}) \cong \mathcal{L}^{\otimes n^2}$.
Indeed, since 
$\mathcal{L}|_{A} = L$
is even, there exists a line bundle on $\mathcal{N}$
on $\mathfrak{B}'$ such that
$[n]^{\ast} ( \mathcal{L}) \cong \mathcal{L}^{\otimes n^2} \otimes \pi^{\ast} 
(\mathcal{N})$.
Note that $0_{\pi}^{\ast} ( \pi^{\ast} 
(\mathcal{N})) = \mathcal{N}$,
and we only have to show
that $0_{\pi}^{\ast} ( \pi^{\ast} 
(\mathcal{N})) \cong \OO_{\mathfrak{B}'}$.
Since $0_{\pi}^{\ast} (\mathcal{L}) 
\cong \OO_{\mathfrak{B}'}$, 
we have $0_{\pi}^{\ast} ( \mathcal{L}^{\otimes n^2}) \cong \OO_{\mathfrak{B}'}$.
On the other hand, we see 
\[
0_{\pi}^{\ast}( [n]^{\ast} (\mathcal{L}))
= ([n] \circ 0_{\pi})^{\ast} (\mathcal{L}) = 0_{\pi}^{\ast} (\mathcal{L}) \cong 
\OO_{\mathfrak{B}'}
.
\]
Since $[n]^{\ast} ( \mathcal{L}) \cong \mathcal{L}^{\otimes n^2} \otimes \pi^{\ast} 
(\mathcal{N})$,
we obtain $0_{\pi}^{\ast} ( \pi^{\ast} 
(\mathcal{N})) \cong \OO_{\mathfrak{B}'}$,
as required.

Let $\mathcal{H}'$ be the pull-back of $\mathcal{H}$ 
by the morphism $\mathfrak{B}' \to \mathfrak{B}$.
We prove that for any $x \in A ( \overline{K} )$, we have
\addtocounter{Claim}{1}
\begin{align} \label{eq:heightnondegenerate}
\widehat{h}_L (x)
=
\frac{\deg_{\mathcal{H}'} 
\left(
\pi_{\ast} 
( \cherncl_1 ( \mathcal{L}) \cdot [\Delta_x])
\right)}{[K'(x) : K]}
\end{align}
where $\Delta_x$ is the closure
of $x$ in $\mathcal{A}$,
$[\Delta_x]$ is the cycle that gives $\Delta_x$,
and
$K'(x)$ is the function field of $\Delta_x$.
Note that
the right-hand side 
in (\ref{eq:heightnondegenerate})
gives a height function associated to $L$.
By Proposition~\ref{prop:ch}.
and Theorem~\ref{thm:canonicalheight},
we only have to show that
\addtocounter{Claim}{1}
\begin{align} \label{eq:heightnondegenerate2}
\frac{\deg_{\mathcal{H}'} 
\left(
\pi_{\ast}  ( \cherncl_1 ( \mathcal{L}) \cdot [\Delta_{n x}])
\right)}{[K'(nx) : K]}
=
n^{2} 
\frac{
\deg_{\mathcal{H}'} 
\left(
\pi_{\ast} 
( \cherncl_1 ( \mathcal{L}) \cdot [\Delta_x])
\right)}{[K'(x) : K]}
.
\end{align}
Since $[n]_{\ast} ([\Delta_{x}])
=
[K'(x) : K'(nx)][\Delta_{nx}]$,
it follows from the projection formula that
as cycle classes, 
\[
[K'(x) : K'(nx)]
 \cherncl_1 ( \mathcal{L}) \cdot [\Delta_{n x}]
= 
 \cherncl_1 ( [n]^{\ast} ( \mathcal{L} )) \cdot [\Delta_{x}]
.
\]
Since
\[
\cherncl_1 ( [n]^{\ast} ( \mathcal{L} )) \cdot [\Delta_{x}]
=
 \cherncl_1 ( \mathcal{L}^{\otimes n^{2}}) \cdot [\Delta_{x}]
=
n^{2}
( \cherncl_1 ( \mathcal{L}) \cdot [\Delta_{x}])
,
\]
it follows that $[K'(x) : K'(nx)]
 \cherncl_1 ( \mathcal{L}) \cdot [\Delta_{n x}]
= n^{2}
( \cherncl_1 ( \mathcal{L}) \cdot [\Delta_{x}])$.
Thus we obtain 
(\ref{eq:heightnondegenerate2}).

\subsection{Height $0$ points}
\label{subsection:height0}

As we mentioned in Remark~\ref{rem:torsion_height0},
torsion points have canonical height $0$.
Over number fields,
it is classically known that the converse also holds;
the torsion points are the only points with canonical height $0$.
Over function fields, however, there are
height $0$ points other than torsion points, in general.
A typical example of an abelian variety which can have
non-torsion height $0$ points is a constant abelian variety.
A \emph{constant abelian variety}
is an abelian variety $B$ over $\overline{K}$
endowed with an identification 
$B = \widetilde{B} \otimes_k \overline{K}$
for some abelian variety $\widetilde{B}$ over $k$. 
Note that for a constant abelian variety 
$B = \widetilde{B} \otimes_k \overline{K}$, 
we regard $\widetilde{B} (k) \subset B (\overline{K})$ naturally.

\begin{Example} \label{ex:height0points}
Let $B:= \widetilde{B} \otimes_k \overline{K}$
be a constant abelian variety.
Note $\widetilde{B} (k) \subset B (\overline{K})$.
Then any point $b \in \widetilde{B} (k)$ has height $0$.
To see that,
we take an  even ample line bundle
$\widetilde{M}$ 
on 
$\widetilde{B}$
and set $M := \widetilde{M} \otimes_k \overline{K}$.
Let  
$\pi : \widetilde{B} \times_{\Spec (k)} \mathfrak{B} \to \mathfrak{B}$ 
and 
$\pr_{\widetilde{B}} : \widetilde{B} \times_{\Spec (k)} \mathfrak{B} \to 
\widetilde{B}$
be the projections.
Then the pair
$( \widetilde{B} \times_{\Spec (k)} \mathfrak{B} , \pr_{\widetilde{B}}^{\ast} 
(\widetilde{M}) )$ is a proper flat model of $(B , M)$
such that $\pi : \widetilde{B} \times_{\Spec (k)} \mathfrak{B} \to \mathfrak{B}$
is an abelian scheme.
(We call such a model a \emph{standard model}.)
For any $b \in \widetilde{B} (k)$,
the closure $\Delta_b$ of $b$
(as a point of $B ( \overline{K})$) 
in 
$\widetilde{B} \times_{\Spec (k)} \mathfrak{B}$ equals
$\pr_{\widetilde{B}}^{-1} (b)$.
Note that $
\cherncl_1 (
\pr_{\widetilde{B}}^{\ast} 
(\widetilde{M}) ) \cdot [\pr_{\widetilde{B}}^{-1} (b)] = 0$
as a cycle class on $\widetilde{B} \times_{\Spec (k)} \mathfrak{B}$.
By (\ref{eq:heightnondegenerate}),
we then obtain
\[
\widehat{h}_{M} (b) 
= \deg_{\mathcal{H}} \pi_{\ast} 
\left( \cherncl_1 ( \pr_{\widetilde{B}}^{\ast} 
(\widetilde{M}) ) \cdot [\pr_{\widetilde{B}}^{-1} (b)] 
\right)
= 0
,
\]
as required.
\end{Example}

\begin{Remark} \label{rem:imageofconstantpoints}
Let $B = \widetilde{B} \otimes_k \overline{K}$ be as
in the above example.
Let $A$ be an abelian variety over $\overline{K}$
and
let $\phi : B \to A$ be a homomorphism.
Then any point of $\phi ( B (k) ) \subset A ( \overline{K} )$
has height $0$ by Lemma~\ref{lem:imageofheith0points}. 
\end{Remark}

To describe the points of height $0$ over function fields, 
we recall the notion of $\overline{K}/k$-trace
of an abelian variety.
It will also be used to give a definition of special subvariety.

A $\overline{K}/k$-trace
is universal
among the homomorphisms from a constant abelian variety to the abelian
variety.
To be precise,
let $A$ be an abelian variety over $\overline{K}$.
A \emph{$\overline{K}/k$-trace} of $A$ is a pair
$\left( \widetilde{A}^{\overline{K}/k},
\Tr_A \right)$ of an abelian variety over $k$
and a homomorphism $\Tr_A : \widetilde{A}^{\overline{K}/k} \otimes_{k}
\overline{K}
\to A$ 
having the following universal property:
for each constant abelian variety $B = \widetilde{B} \otimes_k \overline{K}$ be
and each homomorphism $\phi : B \to A$,
there exists a unique homomorphism $\phi^{t} : 
\widetilde{B} \to \widetilde{A}^{\overline{K}/k}$
such that $\Tr_A \circ ( \phi^{t} \otimes_k \overline{K} ) = \phi$.
It is known that there exists a unique $\overline{K}/k$-trace
of $A$, the uniqueness following from the universal property.
Further, it is known that the $\overline{K}/k$-trace
homomorphism $\Tr_A$ is finite and purely inseparable.
We refer to \cite[Ch.VIII, \S~3]{lang1} for details.

\begin{Proposition} [Theorem~5.4 of Chapter~6 of \cite{lang2}] \label{prop:height0points_ff}
With the above notation, we regard
\[
\widetilde{A}^{\overline{K}/k} (k) \subset
\left(
\widetilde{A}^{\overline{K}/k} \otimes_{k}
\overline{K} 
\right)
( \overline{K} )
.
\]
Then we have
\[
\left\{
a \in A ( \overline{K} )
\left|
\ 
\widehat{h}_L (a) = 0
\right.
\right\}
=
\Tr_A \left(
\widetilde{A}^{\overline{K}/k} (k)
\right)
+ A ( \overline{K} )_{tor}
.
\]
\end{Proposition}

Indeed,
The inclusion $\supset$ is shown as follows.
We take an even ample line bundle $L$ on $A$.
Since $\widehat{h}_L$ is a quadratic form
and since the torsion points have height $0$,
we only have to show that for $a \in \Tr_A \left(
\widetilde{A}^{\overline{K}/k} (k)
\right)$, we have $\widehat{h}_L (a) = 0$;
in fact it is Remark~\ref{rem:imageofconstantpoints}.
For the proof of the other inclusion, we refer to \cite[Chapter~6]{lang2}.

\subsection{Density of small points} \label{subsect:densityofsmallpoints}
In this section, let $K$ be a function field or a number field.
Let $A$ be an abelian variety over $\overline{K}$.
We introduce the notion of
density of small points
and give basic facts on this notion.

The notion of density of small points can be defined
due to the following lemma.

\begin{Lemma} [cf. the proof Lemma~2.1 of \cite{yamaki5}] \label{lem.DOSM}
Let $A$ be an abelian variety over $\overline{K}$
and let $X$ be a closed subvariety.
Let $L_1$ and $L_2$ be even ample line bundles on $A$.
Then the following are equivalent:
\begin{enumerate}
\item
$X ( \epsilon_1 ; L_1 )$ is dense in $X$ for any 
$\epsilon_1 > 0$;
\item
$X ( \epsilon_2 ; L_2 )$ is dense in $X$ for any 
$\epsilon_2 > 0$.
\end{enumerate}
\end{Lemma}

\Proof
The proof is given by an argument similar to the proof of 
Lemma~\ref{lem:height0points}, so that we omit the detail.
\QED

We say that \emph{$X$ has dense small points} if for any $\epsilon > 0$
\[
X ( \epsilon ; L) :=
\left\{ 
x \in X ( \overline{K} )
\left|
\ 
 \widehat{h}_L (x) \leq \epsilon
\right.
\right\}
\]
is dense in $X$.
This notion does not depend on the choice of $L$ 
by Lemma~\ref{lem.DOSM}.
If $X$ has a dense subset of points of height $0$, then it has
dense small points.

Using the terminology of density of small points,
Zhang's theorem is restated as follows.
Let $K$ be a number field.
Let $A$ be an abelian variety over $\overline{K}$
and
let $X$ be a closed subvariety of $A$.
Suppose that $X$ has dense small points.
Then $X$ is a torsion subvariety.

We put together some lemmas concerning the density of small points.
The following lemma indicates that the image of 
a closed subvariety with dense small points by a homomorphism
has the same property.

\begin{Lemma} [cf. Lemma~2.1 of \cite{yamaki5}] \label{lem:density-quotient}
Let $\phi : B \to A$ be a homomorphism of abelian varieties
over $\overline{K}$.
Let $Y$ be a closed subvariety of $B$.
If $Y$ has dense small points,
then $\phi (Y)$ has dense small points.
\end{Lemma}

\Proof
Let $L$ be an even ample line bundle on $A$ and let $M$ be an even ample line
bundle on $B$.
Since $M$ is ample,
$M^{\otimes m} \otimes \phi^{\ast} (L)^{\otimes -1}$
is ample for some $m \in \NN$.
By Proposition~\ref{prop:propertycanonicalheight}~(3),
we have
$\widehat{h}_{M^{\otimes m} \otimes \phi^{\ast} (L)^{\otimes -1}} \geq 0$.
Further, by Proposition~\ref{prop:propertycanonicalheight}~(1)
and (2),
we obtain
$\widehat{h}_{M^{\otimes m}}  \geq \phi^{\ast} \widehat{h}_{L}$.
Put $X := \phi (Y)$.
The last inequality gives us
$\phi (Y ( \epsilon ; M^{\otimes m})) \subset 
X (\epsilon ; L )$.
This shows that if $Y$ has dense small points,
then so does $X$.
\QED

\begin{Lemma} \label{lem:product-density}
Let $A_1$ and $A_2$ be abelian varieties over $\overline{K}$
and let $X_1$ and $X_2$ be closed subvarieties of $A_1$ and $A_2$,
respectively.
Suppose that $X_1$ and $X_2$ have dense small points.
Then $X_1 \times X_2$ has dense small points.
\end{Lemma}

\Proof
This lemma follows from Proposition~\ref{prop:propertycanonicalheight}~(1).
See also \cite[Lemma~2.4]{yamaki5}.
\QED

In the argument later,
we will use the equidistribution theorem of small points.
There, we will use a small generic net of points.
Let $X$ be a closed subvariety of $A$.
An element of $(x_i)_{i \in I} \in X( \overline{K})^{I}$,
where  $I$ is a directed set,
is called a
\emph{net} on $X( \overline{K})$.
We say that $(x_i)_{i \in I}$ is \emph{generic}
if for any proper closed subset $Y$ of $X$,
there exists an $i_0 \in I$ such that $x_{i} \notin Y$
for any $i \geq i_0$.
We say that $(x_i)_{i \in I}$ is \emph{small}
if $\lim_i \widehat{h}_L (x_i) = 0$,
where $L$ is an even ample line bundle on $A$.
The notion of small does not depend on the choice of $L$.

We will later use the following lemma, which
asserts that
there exists a small generic net on a closed subvariety
if the closed subvariety has dense small points.

\begin{Lemma} \label{lem:genericsmallnets}
Let $A$ be an abelian variety over $\overline{K}$
and let $X$ be a closed subvariety of $A$.
Suppose that $X$ has dense small points.
Then there exists a small generic net $( x_i)_{i \in I}$ on 
$X ( \overline{K} )$.
\end{Lemma}

\Proof
If $\dim (X) = 0$, then the assertion is trivial.
Therefore we may assume that $\dim (X) > 0$.
Let $\mathfrak{S}$ be the
set of closed subsets $Z$ of $X$
such that each irreducible component of $Z$ has codimension $1$ in $X$.
Let $I$ be an index set of $\mathfrak{S}$,
which means there exists a bijective map from $I$ to $\mathfrak{S}$.
For each $i \in I$, let $Z_i$ denote the closed subset corresponding to $i$.
We put a partial order on $I$ in such a way that
for $i_1 , i_2 \in I$,
$i_1 \leq i_2$ if and only if $Z_{i_1} \subset Z_{i_2}$.
Then $I$ is a directed set with respect to this order.

We construct a net which will be checked to be generic and small.
We take an even ample line bundle $L$ on $A$.
For any $i \in I$, 
let $n_i$ denote the number of irreducible components of $Z_i$.
Since $X$ has dense small points,
$X ( 1/n_i ; L)$ is dense in $X$.
It follows that for each $i \in I$ there exists a point 
$x_i \in X ( 1/ n_i ; L) \setminus Z_i$.
By the axiom of choice,
this constructs a net $(x_i)_{i \in I}$ on $X ( \overline{K} )$.

We prove that $(x_i)_{i \in I}$ is generic and small.
It is generic by construction. 
Indeed, if we take any closed subset $Y$ of $X$
and an $i_0 \in I$ with $Y \subset Z_{i_0}$,
then for any $i \geq i_0$,
$x_i \notin Z_i \supset Z_{i_0} \supset Y$,
which shows that $(x_i)_{i \in I}$ is generic.
To see that it is small,
we take any $\epsilon > 0$.
There exists $i_0 \in I$ such that $n_{i_0} > \epsilon^{-1}$.
For any $i \geq i_0$, we have $Z_i \supset Z_{i_0}$.
Since any irreducible component of $Z_i$ and $Z_{i_0}$ has codimension
$1$ in $X$, it follows that $n_i \geq n_{i_0}$.
Thus for any $i \geq i_0$,
we have 
\[
0
\leq 
\widehat{h}_{L} ( x_i ) \leq 1/n_i \leq 1/n_{i_0} < \epsilon
,
\]
which shows that $(x_i)_{i \in I}$ is small.
\QED

\section{Geometric Bogomolov conjecture} \label{sect:GBC}

Throughout this section, let $K$ be a function field,
that is, $K$ is the function field
of a normal projective variety $\mathfrak{B}$
over a fixed algebraically closed field $k$
(cf. \S~\ref{subsubsect:NCgroundfields}).

\subsection{Place of $\overline{K}$}

We define the set $M_{\overline{K}}$
of $\overline{K}$.
For a finite extension $K'$ of $K$,
let $\mathfrak{B}'$ denote the normalization of $\mathfrak{B}$ in $K$.
Let $M_{K'}$ be the set of points of $\mathfrak{B}'$ of codimension $1$.
We call each $v \in M_{K'}$ a \emph{place} of $K'$.
(Remark that the notion of place of $K'$ depends 
not only on $K'$ but also on $\mathfrak{B}$ unless $\dim ( \mathfrak{B} )
= 1$.)
If $K''$ is a finite extension of $K'$, then
it is well known that there exists a natural surjective map $M_{K''} \to M_{K'}$,
and thus we have an inverse system $( M_{K'})_{K'}$,
where $K'$ runs through the finite extensions of $K'$ in $\overline{K}$.
We define $M_{\overline{K}}$ to be
$\varprojlim M_{K'}$.
We call each element of $M_{\overline{K}}$ a \emph{place of $\overline{K}$}.

Each place $v \in M_{\overline{K}}$ gives a unique 
equivalence class of non-archimedean absolute
values on $\overline{K}$ such that
if $v_{K'}$ denote the natural projection of $v$ to $M_{K'}$,
then restriction of the value to $K'$ is equivalent to
the value on $K'$ given by the discrete valuation ring
$\OO_{\mathfrak{B}' , v_{K'}}$.
Let $\overline{K}_v$ denote the completion of $\overline{K}$
with respect to this value.
Let $\overline{K}_v^{\circ}$ denote the ring of integer of $\overline{K}_v$.
This has residue field $k$.

For an algebraic variety $X$ over $\overline{K}$,
we write $X_v := X \otimes_{\overline{K}} \overline{K}_v$.

\subsection{Degeneration of abelian varieties} \label{subsect:degeneration} 

We recall the notion of degeneracy of abelian varieties.
Let $A$ be an abelian variety over $\overline{K}$.
We take a $v \in M_{\overline{K}}$.
Then by the semistable reduction theorem,
there exists a unique semiabelian scheme $\mathscr{A}^{\circ} \to 
\Spec ( \overline{K}_v^{\circ} )$
whose generic fiber equals $A_v$.
(Remark that $A$ can be defined over the quotient field
of some discrete valuation ring.)
Let $\widetilde{\mathscr{A}^{\circ}}$ be the special fiber.
Since it is a semiabelian variety,
there exist a nonnegative integer $r$
and
an exact sequence
\addtocounter{Claim}{1}
\begin{align}
\label{align:raynaud:special}
1 \to (\GG_{m}^r)_k \to \widetilde{\mathscr{A}^{\circ}} \to \widetilde{B} \to 0
\end{align}
of algebraic group over $k$,
where $(\GG_{m}^r)_k$ is the algebraic torus over $k$ of dimension $r$
and $\widetilde{B}$ is an abelian variety over $k$.
We call $\widetilde{B}$ the abelian part of 
the reduction of $A_v$.
Put $b (A_v) := \dim ( \widetilde{B} )$,
which
is well defined from $A_v$.

We say that $A$ is \emph{degenerate at $v$} if 
$b (A_v) < \dim (A)$.
We say $A$ is \emph{non-degenerate at $v$} if it
is not degenerate at $v$,
i.e., $b (A_v) = \dim (A)$.
Further, we say that $A$ is 
\emph{totally degenerate at $v$} if 
$b (A_v) = 0$.
Notice that in this terminology,
if $\dim (A) = 0$, then $A$ is non-degenerate 
and totally degenerate at any place.

We restate Gubler's theorem here,
as we have just defined the notion of total degeneracy.

\begin{Theorem} [Gubler's theorem, Theorem~1.1 of \cite{gubler2}]
\label{thm:gubler:text}
Assume that $A$ is totally degenerate at some place,
that is, there exists $v \in M_{\overline{K}}$
such that $A_v$ is totally degenerate.
Let $X$ be a closed subvariety of $A$.
Suppose that $X$ has dense small points.
Then $X$ is a torsion subvariety.
\end{Theorem}

\subsection{Statement of the conjecture} \label{subsect:statementGBC}

For a general abelian variety $A$ over $\overline{K}$,
the statement of Gubler's theorem does not hold:
let $( \widetilde{A}^{\overline{K}/k} ,\Tr_A )$
be the $\overline{K}/k$-trace of $A$;
then the image of a constant subvariety of $\widetilde{A}^{\overline{K}/k}$
by $\Tr_A $ has dense small points,
and in general it is not a torsion subvariety.
This fact indicates that if we wish to characterize the closed subvarieties
with dense small points,
we have to define a suitable counterpart
of torsion subvariety.
The candidate is the class of
special subvarieties, which we are going to define.

Let $X$ be a closed subvariety of $A$.
We say that $X$ is \emph{special}
if there exist an abelian subvariety $G$ of $A$,
a torsion point $\tau \in A ( \overline{K} )_{tor}$,
and a closed subvariety $\widetilde{Y} \subset 
\widetilde{A}^{\overline{K}/k}$ such that
\[
X = \Tr_A \left(
\widetilde{Y} \otimes_{k} \overline{K}
\right)
+ G  + \tau
.
\]
This definition of special subvariety coincides with the definition
in \S~\ref{subsect:initroGBC} by the universal property of the
$\overline{K}/k$-trace.

One shows that 
if $A$ is totally degenerate at some place,
then it has trivial $\overline{K}/k$-trace.
In the setting of Gubler's theorem, therefore, the 
notion of special subvarieties
is the same as that of torsion subvarieties.

A point $x \in A ( \overline{K} )$ is said to be 
\emph{special} if $\{ x \}$ is a special subvariety.
By definition, $x$ is a special point if and only if
\[
x \in 
\Tr_A \left(
\widetilde{A}^{\overline{K}/k} (k)
\right)
+ A ( \overline{K} )_{tor}
,
\]
where we regard $\widetilde{A}^{\overline{K}/k} (k) \subset
\widetilde{A}^{\overline{K}/k} ( \overline{K})$
naturally.

\begin{Remark} \label{rem:height0point=specialpoint}
By Proposition~\ref{prop:height0points_ff},
it follows that a point is special if and only if it has height $0$.
\end{Remark}

\begin{Remark} \label{rem:specialhasdense0}
Any special subvariety has dense small points.
Indeed, with the above expression of $X$,
\[
\Tr_A \left(
\widetilde{Y} (k)
\right)
+ G ( \overline{K} )_{tor} + \tau
\]
is a dense subset of $X$,
and 
by Proposition~\ref{prop:height0points_ff},
each point of this subset has height $0$.
\end{Remark}

Now, we state the geometric Bogomolov conjecture for abelian varieties,
which has been proposed as \cite[Conjecture~2.9]{yamaki5}
and
asserts that the converse of Remark~\ref{rem:specialhasdense0}
should hold.

\begin{Conjecture} [Geometric Bogomolov conjecture for abelian varieties]
\label{conj:GBCintext}
Let $A$ be an abelian variety over $\overline{K}$,
where $K$ is a function field.
Let $X$ be a closed subvariety of $A$.
Suppose that $X$ has dense small points.
Then $X$ is a special subvariety.
\end{Conjecture}

In view of Conjecture~\ref{conj:GBCintext},
Gubler's theorem is a partial answer to the conjecture
because a spacial subvariety is a torsion subvariety in 
the totally degenerate setting.

\begin{Remark} \label{rem:trivialcaseforGBC}
We keep the setting of 
Conjecture~\ref{conj:GBCintext}.
\begin{enumerate}
\item
Assume that $\dim (X) = 0$.
Suppose that $X$ has dense small points.
Then $X$ is special.
Indeed, since $\dim (X) = 0$ 
and $X$ has dense small points, we can write $X = \{ x \}$
with a point $x$ of height $0$.
Then $X$ is special by Remark~\ref{rem:height0point=specialpoint}.
\item
Assume that $\dim (A) = 1$.
Then the geometric Bogomolov conjecture holds
for $A$,
because any proper closed subvariety 
of $A$
has dimension $0$. 
\end{enumerate}
\end{Remark}

\subsection{Partial answers to the conjecture}

Although the geometric Bogomolov conjecture is still open,
there are some important partial answers.
In
\cite{yamaki5}, where Conjecture~\ref{conj:GBCintext}
is proposed, 
we prove the following theorem.

\begin{Theorem} [Theorem~0.5 of \cite{yamaki5}] \label{thm:yamaki5}
Let $A$ be an abelian variety over $\overline{K}$.
Assume that $b (A_v) \leq 1$ for some $v \in M_{\overline{K}}$.
Suppose that $X$ has dense small points.
Then $X$ is a special subvariety.
\end{Theorem}

The above theorem generalizes Gubler's theorem
because
in his setting, we have 
$b (A_v) = 0$ for some $v \in M_{\overline{K}}$.

In \cite{yamaki6},
we generalize Theorem~\ref{thm:yamaki5},
where we use the maximal nowhere degenerate abelian subvariety.
An abelian variety over $\overline{K}$ is said to be 
\emph{nowhere degenerate} 
if it is non-degenerate at any place of $\overline{K}$.
For any abelian variety $A$ over $\overline{K}$,
there exists a unique maximal nowhere degenerate abelian subvariety
$\mathfrak{m}$ of $A$;
$\mathfrak{m}$ is an abelian subvariety that is characterized by
the conditions that
$\mathfrak{m}$ is nowhere degenerate, and
that
if $\mathfrak{m}'$ is an abelian subvariety of $A$ that is nowhere degenerate,
then $\mathfrak{m}' \subset \mathfrak{m}$.
By considering the dimension,
the existence is obvious;
the uniqueness is actually proved in 
\cite[\S~7.3]{yamaki6}.

\begin{Theorem} [Theorem~E of \cite{yamaki6}] \label{thm:yamaki6}
Let $A$ be an abelian variety over $\overline{K}$
and let $\mathfrak{m}$ be the maximal nowhere degenerate
abelian subvariety of $A$. Then the following are equivalent.
\begin{enumerate}
\renewcommand{\labelenumi}{(\alph{enumi})}
\item
The geometric Bogomolov conjecture holds for $A$.
\item
The geometric Bogomolov conjecture holds for $\mathfrak{m}$.
\end{enumerate}
\end{Theorem}

The essential part in Theorem~\ref{thm:yamaki6}
is that (b) implies (a).
In fact, the other implication holds
because
if the geometric Bogomolov conjecture holds for $A$,
then it holds for any abelian subvariety of $A$.

As in Remark~\ref{rem:trivialcaseforGBC}~(2),
the conjecture holds for abelian varieties of dimension 
at most
$1$.
Thus
we have the following corollary.

\begin{Corollary} [Theorem~F of \cite{yamaki6}] \label{cor:ndr=1}
With the notation as in Theorem~\ref{thm:yamaki6},
suppose that $\dim ( \mathfrak{m} ) \leq 1$.
Then the geometric Bogomolov conjecture holds for $A$.
\end{Corollary}

It is not difficult to see that 
Corollary~\ref{cor:ndr=1} generalizes Theorem~\ref{thm:yamaki5}
if you use the Lemma~\ref{lem:ndrandb} below;
see \cite[Proposition~3.3]{yamaki6} for a proof.

\begin{Lemma} \label{lem:ndrandb}
Let $A$ be an abelian variety over $\overline{K}$
and let $A'$ be an abelian subvariety of $A$.
Then for any $v \in M_{\overline{K}}$,
we have $b (A'_v ) \leq b (A_v)$.
\end{Lemma}

Indeed,
by this lemma, we have $\dim ( \mathfrak{m} )
= b (\mathfrak{m}_v) \leq b (A_v)$
for any $v \in M_{\overline{K}}$.
Therefore,
$b (A_v) \leq 1$ implies $\dim ( \mathfrak{m}) \leq 1$,
and thus Theorem~\ref{thm:yamaki6} generalizes Theorem~\ref{thm:yamaki5}.

By Theorem~\ref{thm:yamaki6}, the geometric Bogomolov
conjecture
is reduced to the conjecture
for nowhere degenerate abelian varieties.
In \cite{yamaki7}, furthermore, we show that
the conjecture is reduced to the case 
where the abelian variety has trivial 
$\overline{K}/k$-trace.
To be precise,
let $A$ be an abelian variety over $\overline{K}$
and let $\mathfrak{m}$ be the maximal nowhere degenerate
abelian variety of $A$.
Let $\left( \widetilde{A}^{\overline{K}/k } , \Tr_A \right)$ be the
$\overline{K}/k$-trace.
We set $\mathfrak{t} := \Image ( \Tr_A )$,
the image of the trace homomorphism $\Tr_A : 
\widetilde{A}^{\overline{K}/k } \otimes_{k} \overline{K} \to A$.
One can show that
the image of the maximal nowhere degenerate abelian subvariety
by a homomorphism is contained the maximal nowhere degenerate abelian 
subvariety (cf. \cite[Lemma~7.8~(2)]{yamaki6}).
Since the maximal nowhere degenerate abelian subvariety
of $\widetilde{A}^{\overline{K}/k } \otimes_k \overline{K}$ equals itself,
we have $\mathfrak{t} \subset \mathfrak{m}$,
and hence we can take the quotient $\mathfrak{m} / \mathfrak{t}$.

\begin{Theorem} \label{thm:yamaki7}
Let $A$ be an abelian variety over $\overline{K}$
and let $\mathfrak{m}$ be the maximal nowhere degenerate
abelian subvariety of $A$. Let $\mathfrak{t}$ be the image of
the $\overline{K}/k$-trace homomorphism.
Then the following are equivalent.
\begin{enumerate}
\renewcommand{\labelenumi}{(\alph{enumi})}
\item
The geometric Bogomolov conjecture holds for $A$.
\item
The geometric Bogomolov conjecture holds for $\mathfrak{m}$.
\item
The geometric Bogomolov conjecture holds for $\mathfrak{m}/\mathfrak{t}$.
\end{enumerate}
\end{Theorem}

The equivalence between the first two statement is nothing but
Theorem~\ref{thm:yamaki6}.
By \cite[Lemma~7.7]{yamaki6},
we see that (b) implies (c), which is not difficult.
Thus
the main contribution the theorem
the assertion that (c) implies (b)

Since the quotient of nowhere degenerate abelian variety
by an abelian subvariety is again nowhere degenerate,
$\mathfrak{m}/\mathfrak{t}$ is nowhere degenerate
(cf. \cite[Lemma~7.8~(2)]{yamaki6}).
Furthermore, $\mathfrak{m}/\mathfrak{t}$ has trivial $\overline{K}/k$-trace;
see \cite[Remark~5.4]{yamaki7}.
Thus the geometric Bogomolov conjecture for abelian varieties
is reduced the conjecture
for nowhere degenerate abelian varieties with trivial $\overline{K}/k$-trace.

In a large part of the sequel, we will 
explain the idea of the proof of Theorems~\ref{thm:yamaki6}
and \ref{thm:yamaki7}.
As for Theorem~\ref{thm:yamaki6},
we give an idea of the proof in \S~\ref{sect:reductiontoNDcase};
we will not prove this theorem but give an outline the proof
of a little weaker result.
We recall 
in \S~\ref{sect:proofofgubler}
the proof of Gubler's theorem
because
the idea of our proof is inspired by the proof of Gubler's theorem.
Since the Gubler's proof is a non-archimedean  analogue of the proof
of Zhang's theorem,
we also recall the idea of Zhang in \S~\ref{sect:zhang'sproof}
and fundamental facts on 
non-archimedean geometry in 
\S~\ref{sect:nonarchimedeangeom}.
We recall the structure of the canonical measures in 
\S~\ref{sect:structuretheorem},
as these measures 
are the key ingredients of the proof of Theorem~\ref{thm:yamaki6}.
As for Theorem~\ref{thm:yamaki7},
we give an outline of the proof in \S~\ref{sect:proofofGBCforcurves}
after recalling the notion of canonical heights of closed subvarieties
in \S~\ref{sect:heightofsubvarietiels}.

\section{Proof of 
Zhang's theorem} \label{sect:zhang'sproof}


We
recall Zhang's theorem:

\begin{Theorem} [Zhang's theorem, restated]
Assume that $K$ is a number field.
Let $A$ be an abelian variety over $\overline{K}$.
Let $X$ be a closed subvariety of $A$.
Suppose that $X$ has dense small points.
Then $X$ is a torsion subvariety.
\end{Theorem}

The proof is based on the equidistribution of small points
argued in \cite{SUZ}.
In this section, we give an outline of the proof of Zhang's theorem
with an emphasis on how to use the equidistribution theorem.
The proof of Zhang will be a prototype of the proof of many
results concerning
the geometric Bogomolov conjecture.
The basic reference is Zhang's original paper \cite{zhang2}.
The survey paper \cite{KMY} will be an accessible exposition
for Zhang's theorem.

\subsection{Archimedean canonical measures}
In this subsection, we work over $\CC$.
We begin by recalling the notion of canonical metrics on an
abelian variety,
which is studied in \cite{MB}.
We refer to \cite[\S~4]{KMY} for details of this subsection.
Let $A$ be the complex analytic space
associated to
an abelian variety over $\CC$
and let $L$ be an even line bundle on $A$.
For an integer $n > 1$,
fix an isomorphism $\phi : [n]^{\ast}( L) \to L^{\otimes n^{2}}$.
Let $|| \cdot ||$ be a metric on $L$.
We call $|| \cdot ||$ a 
\emph{canonical metric}
on $L$ if $[n]^{\ast} || \cdot || = || \cdot ||^{\otimes n^{2}}$ 
holds via the isomorphism $\phi$.
Later we consider a nonarchimedean analogue of this metric,
so when we emphasize that we are working over $\CC$,
we call it an archimedean canonical metric.

\begin{Theorem}
For a fixed isomorphism $\phi : [n]^{\ast}( L) \to L^{\otimes n^{2}}$,
there exists a unique canonical metric $|| \cdot ||$ on $L$.
Furthermore, it is a smooth metric.
\end{Theorem}

We denote by $|| \cdot ||_{can}$ the canonical metric.
Further, 
we write $\overline{L} = ( L , || \cdot ||_{can} )$
in this subsection.
Since the canonical metric is a smooth metric,
we consider the curvature form
$\cherncl_1 ( \overline{L} )$.

\begin{Remark} \label{rem:differencecanonicalmetric}
The canonical metric depends on the choice of the isomorphism 
$\phi : [n]^{\ast} L \to L^{\otimes n^{2}}$,
but $\cherncl_1 ( \overline{L} )$ does not.
Indeed, one shows that
a different choice of $\phi$ makes the canonical metric change
only by positive constant multiple,
so that it does not change $\cherncl_1 ( \overline{L} )$.
\end{Remark}

Next,
let us give an explicit description 
of the curvature forms of canonical metrics
in terms of the universal
covering of $A$.
Let $p : \CC^{n} \to  A$ be the universal covering
which is a homomorphism
with respect to the additive structure on $\CC^n$ and the group structure
on $A$.
Let $z_1 , \ldots , z_n$ be the standard coordinates of $\CC^{n}$.
Then there exists a unique
hermitian matrix $(c_{ij}) \in M_{n} ( \CC )$ such that
\addtocounter{Claim}{1}
\begin{align} \label{eq:archimedeancanonicalmeasure-uc}
p^{\ast}
(\cherncl_1 ( \overline{L} ))
= \sum_{1 \leq i,j \leq n}
\sqrt{-1} c_{ij}  dz_i \wedge d \overline{z_j}
,
\end{align}
where $\overline{z_j}$ denotes the complex conjugate of $z_j$;
see \cite[\S~4]{KMY} for more details.

Furthermore,
suppose that
$L$ is ample.
Then
the matrix $(c_{ij})$ is positive definite.
Note that this shows that $\cherncl_1 ( \overline{L} )$
is a positive $(1,1)$-from.

\begin{Remark} \label{rem:positiveform}
Suppose that $L$ is ample.
Let $U$ be a submanifold of $A$ of dimension $d$.
Then $\cherncl_1 ( \overline{L} )^{\wedge d} |_{U}$ is positive
 in the following
sense.
Take any $p \in U$.
Let $u_1 , \ldots , u_d$ be a system of local
holomorphic coordinates
of $U$ around $p$.
We write
\[
\cherncl_1 ( \overline{L} )^{\wedge d} |_{U} 
= (\sqrt{-1})^{d} \varphi du_1 \wedge d\overline{u_1} 
\wedge
\cdots \wedge du_d \wedge d\overline{u_d} 
,
\] 
where $\varphi$ is a smooth function on $U$ around $p$.
Then
$\varphi$ is a positive real valued function.
(This follows from the positivity of 
$\cherncl_1 ( \overline{L} )$.)
\end{Remark}

Finally, we define the canonical measures.
Let $L$ be an even line bundle on $A$.
Let $\overline{L}$ denote the line bundle $L$
equipped with a canonical metric on $L$.
Assume that $L$ is ample.
Let $X$ be a closed subvariety of $A$.
Then
\[
\mu_{X , L} := \frac{\cherncl_1( \overline{L})^{\wedge d} |_{X}}{\deg_L (X)}
\]
is naturally a positive 
regular Borel measure on $X$ with total volume $1$.
We call this probability measure the 
\emph{canonical measure associated} to $L$.
By Remark~\ref{rem:differencecanonicalmetric},
the canonical measure 
does not depend on the choice of a canonical metric on $L$,
and it is well-defined for $L$.

\subsection{Archimedean equidistribution theorem}

Let $K$ be number field.
Let $A$ be an abelian variety over $\overline{K}$
and let $X$ be a closed subvariety of $A$.
Fix a finite extension $K'$ of $K$ such that
$A$ and $X$ can be defined over $K'$.
Set $G_{K'} := \Gal ( \overline{K} / K')$.
Then $G_{K'}$ acts on $X ( \overline{K} )$.
For any $x \in X ( \overline{K} )$,
let $O(x)$ denote the $G_{K'}$-orbit of $x$.
Let $\sigma$ be an archimedean place of $\overline{K}$,
that is, an embedding 
$\overline{K} \hookrightarrow \CC$ of fields.
Let $X_{\sigma}^{\an}$ be the complex analytic space associated to 
$X \otimes_{\overline{K}} \CC$,
where $\CC$ is regarded as an $\overline{K}$-algebra via $\sigma$.
We
regard $X ( \overline{K} ) \subset X_{\sigma}^{\an}$,
and hence $O(x) \subset X_{\sigma}^{\an}$
for any $x \in X ( \overline{K} )$.

\begin{Theorem} [Theorem~2.1 of \cite{zhang2}] \label{thm:archimedeanequidistribution}
With the notation above,
let $(x_{i})_{i \in I}$ be a small generic net on $X ( \overline{K} )$.
Then we have a weak convergence
\[
\lim_{i}
\frac{1}{\# O (x_i)}
\sum_{z \in O (x_i)} \delta_{z}
=
\mu_{X_{\sigma}^{\an},L}
,
\]
where $\delta_z$ denotes the Dirac measure of $z$.
\end{Theorem}

We do not mention the proof of the equidistribution theorem.
We refer to the original paper \cite{zhang2} or
\cite[\S~6]{KMY}.

\subsection{Difference morphism}

Let $A$ be an abelian variety over any algebraically closed field
$\mathfrak{K}$.
Let $X$ be a closed subvariety of $A$.
Let $N$ be a natural number with $N > 1$.
We define  $\alpha_N : X^{N} \to A^{N-1}$ by
\[
\alpha_N (x_1 , \ldots , x_N) = (x_1 - x_2 , \ldots , x_{N-1} - x_N)
.
\]
We call $\alpha_N$ the \emph{difference morphism}.
In this subsection, we give remarks on 
this morphism.
Lemma~\ref{lem:genericallyfiniteN} below will be significantly used not only
in the proof of Zhang but also in the proof
of some results on the geometric Bogomolov conjecture.

The \emph{stabilizer} of $X$,
denoted by $G_X$, is a reduced closed subgroup scheme of $A$
characterized by
\[
G_X (\mathfrak{K}) = \{ a \in A (\mathfrak{K})
\mid X + a \subset X \}
.
\]

\begin{Lemma} [cf. Lemma~4.1 of \cite{abbes}] \label{lem:genericallyfiniteN}
Suppose that $G_X = 0$.
Then there exists an $N \in \NN$ such that
$\alpha_{N} : X^{N} \to A^{N-1}$
is generically finite.
\end{Lemma}

\Proof
For a subset $S$ of $X ( \mathfrak{K} )$,
we set
\[
G_{X , S} := 
\{
a \in A ( \mathfrak{K} )
\mid
S + a \subset X
\}
.
\]
Further,
for any $x_1 , \ldots , x_m \in X ( \mathfrak{K} )$,
we write $G_{X , x_1 , \ldots , x_m} := G_{X, \{ x_1 , \ldots , x_m \}}$.
We prove that
for any $y_1 , \ldots , y_N \in X ( \mathfrak{K} )$,
we have
\[
\alpha_{N}^{-1}
( \alpha_N ( y_1 , \ldots , y_N ) 
(\mathfrak{K} ))=
\{
(y_1 + a , \ldots , y_N + a)
\in A ( \mathfrak{K} )
\mid
a \in G_{X , y_1 , \ldots , y_N}
\}
.
\]
The inclusion ``$\supset$''
is obvious.
To show the other inclusion,
we
take any $(y_1' , \ldots , y_N')
\in \alpha_{N}^{-1}
( \alpha_N ( y_1 , \ldots , y_N ) 
(\mathfrak{K} ))$.
Then for any $i = 1 , \ldots , N-1$,
we have
$y_i' - y_{i+1}' = y_{i} - y_{i+1}$,
and hence $y_1' - y_1 = \cdots = y_N' - y_N$.
Let $a$ denote this element.
Then we have $a \in G_{X , y_1 , \ldots , y_N}$
and $(y_1' , \ldots , y_N' ) = (y_1 + a , \ldots , y_N + a)$,
which shows ``$\subset$''.

For any $S \subset X ( \mathfrak{K} )$,
we note that
$G_X ( \mathfrak{K} ) \subset G_{X , S}$.
Further, we
see that
$G_X ( \mathfrak{K} ) = \bigcap_{S} G_{X , S}$
where $S$ runs through all the finite subset of $X ( \overline{K} )$.
Since $G_X ( \mathfrak{K} )$ and $G_{X , S}$
are closed subsets of $A ( \mathfrak{K} )$,
it follows that there exist $x_1 , \ldots , x_N \in X ( \mathfrak{K} )$
such that $G_X ( \mathfrak{K} ) = G_{X , x_1 , \ldots , x_N}$.
Now,
since $G_X = 0$,
there exist an $N \in \NN$
and $x_1 , \ldots , x_N$ such that
$G_{X, x_1, \ldots , x_N } = G_X = 0$.
It follows 
that
\[
\dim
\left(
 \alpha_{N}^{-1}
(\alpha_N ( x_1 , \ldots , x_N ) )
\right)
=
0
.
\]
This shows that $\alpha_N$ is generically finite.
(In fact, the above argument
shows that $\alpha_N$ is generically injective.)
\QED

\subsection{Proof of Zhang} \label{subsect:proofofZhang}
We start the proof of Zhang.
It is argued by contradiction;
suppose that we have a counterexample to 
Zhang's theorem,
that is,
suppose that
there exist an abelian variety $A$ over $\overline{K}$ 
and
a closed subvariety $X$ 
such that
$X$ is not a torsion subvariety but has dense small points.
Let $G_X$ be the stabilizer of $X$
and
consider the quotient $\phi : A \to A/G_X$.
Then $X / G_X$ 
is a closed subvariety and has trivial stabilizer.
Further, 
$X / G_X$ has dense small points
by Lemma~\ref{lem:density-quotient}.

We prove that $\dim (X/G_X) > 0$
by contradiction;
suppose that $\dim (X/G_X) = 0$.
Then $X/G_X = \{ \phi (x) \}$ for some point $x$.
Since $X/G_X$ has dense small points,
$\phi (x)$ has height $0$,
and hence it is a torsion point.
Since $\phi$ induces a surjective homomorphism between the 
subgroups of torsion points
(cf. \cite[Proof of Lemma~2.10]{yamaki5}), we 
may take $x$ to be a torsion point of $A$.
Then we see that $X = G_X + x$,
which means that $X$ is a torsion subvariety.
That is a contradiction.

The above argument suggests that
replacing $A$ and $X$ by $A/G_X$ and $X/G_X$ respectively if necessary,
we have an abelian variety $A$ and a closed subvariety
$X$ with $\dim (X) > 0$
such that
$X$ has dense small points and has trivial stabilizer.
Put $d := \dim (X)$.
Since $X$ has trivial stabilizer,
there exists,
by Lemma~\ref{lem:genericallyfiniteN},
an integer $N \geq 2$
such that
\[
\alpha_N : X^{N} \to A^{N-1};
\quad 
(x_1 , \ldots , x_N ) \mapsto
(x_1 - x_2 , \ldots ,x_{N-1} - x_N)
\]
is generically finite.
Set $Z:= X^{N}$
and $Y:= \alpha_{N} (Z)$.
Let $\alpha : Z \to Y$ denote the restriction of $\alpha_N$,
which is a generically finite surjective morphism.
We remark that $\dim (Z) = \dim (Y) = dN$.
Since $X$ has dense small points,
so does $Z$ by Lemma~\ref{lem:product-density}.
By Lemma~\ref{lem:genericsmallnets},
there exists a small generic net $(z_i)_{i \in I}$
on $Z(\overline{K})$.
The image $( \alpha (z_{i}))_{i \in N}$
is a generic net on $Y(\overline{K})$.
Further,
if follows from
Proposition~\ref{prop:propertycanonicalheight}~(2)
that this net is small as well.
Take an archimedean place $\sigma$ of $\overline{K}$,
that is, an embedding $\sigma : \overline{K} \hookrightarrow \CC$.
Fix even ample line bundles $M$ and $L$ on $A^{N}$ and $A^{N-1}$, respectively.
Let $\mu_{Z_{\sigma}^{\an} , M}$
and $\mu_{Y_{\sigma}^{\an} , L}$
be the canonical measures on $Z_{\sigma}^{\an}$
and $Y_{\sigma}^{\an}$,
respectively.
By
the (archimedean) equidistribution theorem 
(Theorem~\ref{thm:archimedeanequidistribution}),
we have
\[
\lim_{i}
\frac{1}{\# O (z_i)}
\sum_{u \in O (z_i)} \delta_{u}
=
\mu_{Z_{\sigma}^{\an},M}
\]
and
\[
\lim_{i}
\frac{1}{\# O (\alpha (z_i))}
\sum_{v \in O ( \alpha(z_i) )} \delta_{v}
=
\mu_{Y_{\sigma}^{\an},L}
.
\]
Let $\alpha^{\an} : Z_{\sigma}^{\an} \to Y_{\sigma}^{\an}$ denote the
morphism 
of analytic spaces associated to $\alpha$.
Since
\[
\alpha^{\an}_{\ast}
\left(
\frac{1}{\# O (z_i)}
\sum_{u \in O (z_i)} \delta_{u}
\right)
=
\frac{1}{\# O (\alpha (z_i))}
\sum_{v \in O ( \alpha(z_i) )} \delta_{v}
,
\]
we obtain
$\alpha^{\an}_{\ast} ( \mu_{Z_{\sigma}^{\an} ,M} )= \mu_{Y_{\sigma}^{\an},L}$.

Let $V$ be the nonsingular locus of $X_{\sigma}^{\an}$.
Then $V^N$ is a non-empty nonsingular open subset of
$Z_{\sigma}^{\an}$.
Since $\alpha : Z \to Y$ is generically finite
and
 $\dim (Z) = \dim (Y) = dN$,
$\alpha^{\an}_{\ast} ( \mu_{Z_{\sigma}^{\an},M} ) = \mu_{Y_{\sigma}^{\an},L}$
gives us an equality
\addtocounter{Claim}{1}
\begin{align} \label{eq:equalityofforms}
\frac{\cherncl_1 ( \overline{M})^{\wedge dN}|_{V^{N}}}{\deg_{M} (Z)}
=
\frac{(\alpha^{\an}|_{V^{N}})^{\ast}  \cherncl_1 ( \overline{L})^{\wedge dN}}{\deg_{L} (Y)}
\end{align}
of smooth $(dN , dN)$-forms on $V^N$.
Let $p \in V^{N}$ be a point on the diagonal.
Then the left-hand side on (\ref{eq:equalityofforms}) is positive at $p$
(Remark~\ref{rem:positiveform}).
On the other hand,
since $\alpha^{\an}|_{V^{N}}$ is ramified at $p$,
the right-hand side is not positive at $p$.
That is a contradiction.
Thus the proof
of Zhang's theorem
is complete.

\section{Nonarchimedean geometry} \label{sect:nonarchimedeangeom}

In the proof of Zhang, the following are crucially used:
\begin{enumerate}
\item
analytic spaces over an archimedean place;
\item
canonical measures over analytic spaces;
\item
equidistribution theorem of small points over analytic spaces.
\end{enumerate}
Over function fields, we do not have (1) above:
there does not exist
an archimedean
place over function fields.

Therefore,
if we wish to follow the idea
of Zhang (and Ullmo), we need counterparts.
The counterparts of (1) that we will use
are the
Berkovich analytic spaces over some nonarchimedean place.
Further,
over Berkovich analytic spaces, we can consider ``canonical measures''.

In this section, we briefly review 
the notions of Berkovich analytic spaces,
skeleta, and measures on Berkovich spaces.
The reader familiar with those notions
can skip this section.

Basic references are Berkovich's original papers
\cite{berkovich1,berkovich2,berkovich3,berkovich4}.
For the exposition in this section,
we refer to
Nicaise's exposition \cite{nicaise}.
This
is very accessible to non-experts.

\subsection{Notation and convention}

Throughout this section,
let $\KK$ be an algebraically 
closed field complete with respect to 
a non-trivial non-archimedean value $| \cdot | =| \cdot |_{\mathbb{K}}$.  
Let $\KK^\circ$ denote the ring of integers of $\KK$
and
let $k$ denote the residue field.

Let $\mathbf{K}$ be a subfield of $\mathbb{K}$.
The restriction $| \cdot |_{\bfK}$ of $| \cdot |_{\KK}$
is an absolute value on $\bfK$.
We always assume that
$\mathbf{K}$ is complete with respect to $| \cdot |_{\bfK}$.
Let $\bfK^{\circ}$ denote the ring of integers of $\mathbf{K}$.

We say that $\mathbf{K}$
is a \emph{complete discretely valued subfield} 
if
$| \cdot |_{\bfK}$
is a complete discrete absolute value on $\bfK$.
We abbreviate the name to a \emph{CDV subfield}.
A subring that is the ring of integers of some
CDV subfield is called the \emph{complete discrete valuation
subring}, abbreviated to a \emph{CDV subring}.

In this paper, we mainly use non-archimedean geometry
over 
$\mathbf{K} = \mathbb{K}$.
However, we will sometimes consider
non-archimedean geometry over a CDV subfield
because it will help us to make more accessible 
description of  basic notions
on non-archimedean geometry.

Let $X$ be a scheme over $\KK$.
We say that $X$ \emph{can be defined over a discrete valued
field}
if it has a model over some CDV subfield.

\subsection{Berkovich analytic spaces and skeleta} \label{subsect:berkovich}
Let $\mathbf{K}$ be 
$\KK$ or
a CDV subfield of $\mathbb{K}$.
Let $X$ be an algebraic variety over $\mathbf{K}$.
For a point $p \in X$, let 
$\kappa(p)$ denote the residue field at $p$. 
We mean by the \emph{(Berkovich) analytic space}
associated to 
$X$,
which we denote by $X^{\an}$,
a topological space 
given as follows.
As a set, 
\[
  X^{\an} := 
  \left\{
  (p, |\cdot|) \mid \text{$p \in X$ and $|\cdot|$ is an absolute value of 
$\kappa(p)$ extending $|\cdot|_{\mathbf{K}}$}
  \right\}.
\]
We endow $X^{\an}$ with the weakest topology such that 
the map $\iota: X^{\an} \to X, (p, |\cdot|) \mapsto p$ is continuous and 
such that 
for any Zariski open subset $U$ of $X$ and any regular function 
$g \in \OO_X(U)$, 
the map $\iota^{-1}(U) \to \RR, (p, |\cdot|) \mapsto |g(p)|$ is continuous. 
It is known that $X^{\an}$ is a Hausdorff, locally compact,
and locally path-connected space.

Let $x \in X ( \mathbf{K} )$ be a $\mathbf{K}$-valued point.
Then this gives a point of $X^{\an}$.
Indeed,
the natural homomorphism $\mathbf{K} \hookrightarrow
\kappa (x)$ is isomorphism in this case,
so that $\kappa (x)$ has a unique absolute value 
extending $|\cdot|_{\mathbf{K}}$
via this isomorphism.
Thus we have $X ( \mathbf{K} ) \subset X^{\an}$ naturally.
We call a point in $X ( \mathbf{K} )$ a \emph{classical point}.

Let $f : X \to Y$ be a morphism of algebraic varieties.
Then $f$ induces a continuous map $f^{\an} : X^{\an} \to Y^{\an}$ naturally.
Indeed, 
since we have a inclusion $\kappa (f (p)) \hookrightarrow \kappa (p)$
for any $p \in X$,
we assign $(p , | \cdot | ) \in X^{\an}$ to 
$( f (p) ,  | \cdot |_{\kappa (f (p))})$,
where $| \cdot |_{\kappa (f (p))}$ is the restriction of $| \cdot |$
to $\kappa (f (p))$.

We make a remark on the base-change.
Set $X_{\mathbb{K}} := X \otimes_{\mathbf{K}} \mathbb{K}$.
Then we have a natural map $\rho : X_{\mathbb{K}} \to X$.
Furthermore, for each $p \in X_{\mathbb{K}}$,
we have $\kappa ( \rho (p)) \subset \kappa ( p)$,
and
taking the restriction of the absolute value of 
$\kappa ( p)$ to $\kappa ( \rho (p))$,
we obtain the natural map $X_{\mathbb{K}}^{\an} \to X^{\an}$.
One sees that this is continuous and surjective.

\subsubsection{Divisorial points} \label{subsect:DP}
A Berkovich space in general has many other points than classical points.
The divisorial points associated to the generic points of 
irreducible components of the special fiber of a model 
are very important class of them.
To make the explanation simple,
we assume that
$\KK = \overline{K}_v$,
where $K$ is a function field over $k$
and $v \in M_{\overline{K}}$ is a fixed place.
Remark that $k \subset \KK$
and the residue field of $\KK$ equals $k$.

First, we define divisorial points in analytic spaces
over a CDV subfield.
Let $\mathbf{K}$ be a CDV subfield.
Let $X$ be an algebraic variety over $\mathbf{K}$.
Let
$\mathscr{X} \to \Spec (\bfK^{\circ})$ be a 
flat model of
$X$ over $\bfK^{\circ}$,
where $\bfK^{\circ}$ is a discrete valuation ring
by the assumption.
Let $\varpi$ be a uniformizer of $\bfK^{\circ}$.
Let $\Irr ( \widetilde{\mathscr{X}} )$
be the set of irreducible components of $\widetilde{\mathscr{X}}$.
For each $V \in \Irr ( \widetilde{\mathscr{X}} )$,
let $\xi_V$ denote the generic point of $V$ in $\mathscr{X}$.

Assume that 
for any $V \in \Irr ( \widetilde{\mathscr{X}} )$,
$\mathscr{X}$ 
is normal
at $\xi_V$.
Take any $V \in \Irr ( \widetilde{\mathscr{X}} )$.
Then
the local ring $\mathcal{O}_{\mathscr{X} , \xi_V}$ is a discrete valuation ring,
whose fraction field equals the function field $\mathbf{K} ( X )$
of $X$.
Remark that $\mathbf{K} ( X )$ is the residue field at the generic point of $X$.
Let $\ord_{\xi_{V}} : \mathbf{K} ( X ) \to \ZZ \cup \{ + \infty \}$
be the order function ($\ord_{\xi_V} (0) = + \infty$ by convention).
Let $m_V$ denote the multiplicity of $V$ in $\widetilde{\mathscr{X}}$.
Then we define the value 
$| \cdot |_{\xi_V} : \mathbf{K} ( X ) \to \RR$
by $| f |_{\xi_V} :=  | \varpi |^{\ord_{\xi_V} (f)/m_V}$.
It is easy to see that
the restriction of
this absolute value equals $| \cdot |_{\bfK}$,
and thus
it gives a point of $X^{\an}$.
We denote this point by $\xi_V^{\an}$
and call it the \emph{divisorial point associated to $\xi_V$ (or $V$)}.
Set
$\mathrm{DP}(\mathscr{X}) := 
\left\{ \xi_V^{\an} \in X^{\an} 
\left|
\ 
 V \in \Irr ( \widetilde{\mathscr{X}}) 
\right.
\right\}$,
whose points are called
divisorial points with respect to $\mathscr{X}$.

Next, 
we introduce the divisorial points over $\KK$.
Although we can define a similar kind of points in general setting,
we only consider the models that can be
defined over CDV subring,
because that is easy to describe and will be enough for our later use.
Let $X$ be an algebraic variety over $\KK$.
Assume that 
there exist a CDV subfield with $k \subset \mathbf{K}$
and a flat model
$\mathscr{X} \to 
\Spec ( \bfK^{\circ} )$
of $X$
such that the special fiber is reduced.
For any $V \in \Irr ( \widetilde{\mathscr{X}} )$,
$\mathscr{X}$ is then normal at the generic point $\xi_V$ of $V$.
The base-change $\mathscr{X}_{\KK^{\circ}} \to \Spec ( \KK^{\circ} )$
is a flat model of $X$ over $\KK^{\circ}$.
Let $\rho :
\mathscr{X}_{\KK^{\circ}} \to \mathscr{X}$
be the natural map.
Let $\eta$ be the generic point of an irreducible component of 
the special fiber of $\mathscr{X}_{\KK^{\circ}}$.
Then $\rho ( \eta )$ is the generic point of some irreducible component
of $\widetilde{\mathscr{X}}$.
Thus we have a point $\rho ( \eta )^{\an} \in ( \mathscr{X}_{\mathbf{K}})^{\an}$,
where $\mathscr{X}_{\mathbf{K}}$
is the generic fiber of $\mathscr{X} \to \Spec ( \bfK^{\circ})$.
Since the natural map $X^{\an} \to ( \mathscr{X}_{\mathbf{K}})^{\an}$
is surjective, there exists a point in $X^{\an}$
that lies over $\rho ( \eta )$.
By our assumption,
any irreducible component of $\widetilde{\mathscr{X}}$
is geometrically irreducible.
From that,
one can show that
a point of
$X^{\an}$ over $\rho ( \eta )$
is unique.
This unique point
in $X^{\an}$
is called the \emph{divisorial point associated to 
$\eta$}, denoted by $\eta^{\an}$.
Let $\mathrm{DP}(\mathscr{X}_{\KK^{\circ}})$
be 
the set of those divisorial points.
One can show that $\mathrm{DP}(\mathscr{X}_{\KK^{\circ}})$
does not depend on the choice of $\mathbf{K}$,
depending only on $\mathscr{X}_{\KK^{\circ}}$.

\subsection{Skeleta}
In this subsection,
we explain skeleta.
A skeleton is a compact subset of an analytic space
associated to a strictly semistable model.
It has a canonical structure of simplicial set.

\subsubsection{Strictly semistable scheme}

We begin by recalling the notion of strictly semistable schemes
over a discrete valuation ring.
Let $\mathbf{K}$ be a CDV subfield.
Let $f : \mathscr{X} \to \Spec ( \mathbf{K}^{\circ} )$
be a 
flat morphism of finite type.
We say that $f$ is a \emph{strictly semistable scheme}
if it satisfies the following conditions:
\begin{enumerate}
\renewcommand{\labelenumi}{(\alph{enumi})}
\item
$\mathscr{X}$ is regular;
\item
$f$ is generically smooth;
\item
the special fiber $\widetilde{\mathscr{X}}$
is a normal crossing divisor
whose
irreducible components 
are geometrically integral and
smooth.
\end{enumerate}

Next, we define a strictly semistable scheme over $\KK^{\circ}$.
Let $\mathscr{X} \to \Spec ( \KK^{\circ} )$
a morphism of scheme.
It is called a strictly semistable scheme 
if there exist a CDV subfield $\mathbf{K}$
and strictly semistable scheme 
$\mathscr{X}' \to  \Spec ( \mathbf{K} ^{\circ})$
such that $\mathscr{X} = \mathscr{X}'_{\KK^{\circ}}$.

The above definition is an expedient definition;
we have a more general 
definition which makes sense for those
not necessarily defined over a discrete valuation ring.
However, 
we only consider those which can be defined over a discrete valuation
subring, because that is easy to describe and
will be enough for our use.

\subsubsection{Monomial points}

We recall the notion of strata of a reduced algebraic scheme
$Z$
over $k$.
We put $Z^{(0)} := Z$.
For $r \in \ZZ_{\geq 0}$,
we inductively define $Z^{(r+1)} \subset Z^{(r)}$ 
to be the set of non-normal points of $Z^{(r)}$.
Then we obtain a descending sequence of closed subsets
\[
Z = Z^{(0)} \supsetneq
Z^{(1)} \supsetneq
\cdots 
\supsetneq
Z^{(s)}
\supsetneq
Z^{(s + 1)} = \emptyset
.
\]
The irreducible components of $Z^{(r)} \setminus Z^{(r+1)}$ are
called the \emph{strata of $Z$ of codimension $r$}.
Let $\str ( Z )$ be the set
of strata of $Z$.

Let
$X$ be a variety over $\bfK$
with a strictly semistable model
$\mathscr{X}$.
The monomial points are the points
of $X^{\an}$
which are
determined by 
$S \in \str (\mathscr{X})$ 
and positive real numbers
$(u_{0} , \ldots , u_r )$ with $u_0 + \cdots + u_r= 1$
where $r$ is the codimension of $S$.
To give a more precise description,
we first assume that $\bfK$ is a CDV subfield,
and then we explain them for $\KK$.

Assume that $\bfK$ is a CDV subfield.
Let $\varpi$ be a uniformizer of $\bfK^{\circ}$.
Let $X$ be an algebraic variety over $\bfK$.
Assume that $X$ has a strictly semistable model 
$\mathscr{X}$
over $\bfK^{\circ}$,
that is,
there exists a strictly semistable scheme
$\mathscr{X} \to \Spec ( \bfK^{\circ} )$ with generic fiber $X$.
Let $S$ be a stratum
of $\widetilde{\mathscr{X}} $ of codimension $r$ in 
$\widetilde{\mathscr{X}}$
and let
$\xi$ be the generic point of $S$.
Locally at $\xi$ in $\widetilde{\mathscr{X}}$,
$S$ is given by the intersection
$\bigcap_{i = 0}^{r} V_i$
of
$r+1$ irreducible components $V_0 , \ldots , V_r 
\in \Irr ( \widetilde{\mathscr{X}} )$.
(In other words, the closure of $S$ is an irreducible component
of $\bigcap_{i = 0}^{r} V_i$. )
Let $\mathcal{O}_{\mathscr{X} , \xi}$ be the local
ring of $\mathscr{X}$ at $\xi$.
Let $t_0 , \ldots , t_r$ be elements
of $\mathcal{O}_{\mathscr{X} , \xi}$ defining $V_0 , \ldots , V_r$
around $\xi$,
respectively.
Then since $\mathcal{O}_{\mathscr{X} , \xi}$ is regular
and $V_0 , \ldots , V_r$ are normally crossing,
$t_0 , \ldots , t_r$ 
is a regular system of parameters of $\mathcal{O}_{\mathscr{X} , \xi}$.
Note that there exists a unit $\lambda$ in $\mathcal{O}_{\mathscr{X} , \xi}$
such that
$\varpi = \lambda t_0 \cdots t_r$.

Let $\widehat{\OO_{\mathscr{X} , \xi}}$ be the completion of $\OO_{\mathscr{X} , \xi}$.
Let $k ( \xi )$ be the residue field of $\OO_{\mathscr{X} , \xi}$.
By Cohen's structure theorem of complete regular local rings,
there exist a section of 
the quotient homomorphism $\widehat{\OO_{\mathscr{X} , \xi}}
\to k ( \xi )$
and 
an isomorphism
\[
k ( \xi ) [[ 
t_0 , \ldots , t_r ]]
\cong
\widehat{\OO_{\mathscr{X} , \xi}}
.
\]
We call this isomorphism a Cohen isomorphism
for $\widehat{\OO_{\mathscr{X} , \xi}}$.

Put $\Delta^r := 
\{ \mathbf{u} = (u_0 , \ldots , u_r )
\in \RR_{\geq 0}^{r+1}
\mid
u_{0} + \cdots + u_r = 1
\}$.
Then under the setting above,
each $\mathbf{u} \in \Delta^r$ determines a point of $\bfX^{\an}$
due to the following proposition.

\begin{Proposition} [cf. \S~2.3 of \cite{nicaise}] \label{prop:monomialpoints}
Let $\bfK ( X )$ be the function field of $X$.
For any $\mathbf{u} = (u_0 , \ldots , u_r ) \in \Delta^{r}$,
there exists a unique absolute value $| \cdot |_{\mathbf{u}} : \bfK ( X ) \to \RR_{\geq 0}$
with the following property.
Fix a Cohen isomorphism for $\widehat{\mathcal{O}_{\mathscr{X} , \xi}}$.
Take any $f \in \OO_{\mathscr{X} , \xi}$.
Regarding it as a formal power series
via the Cohen isomorphism,
we write
\[
f = \sum_{\mathbf{m} \in \ZZ_{\geq 0}^{r+1}} c_{\mathbf{m}} t^{\mathbf{m}}
.
\]
Then
we have
\[
\log | f |_{\mathbf{u}} = ( \log | \varpi | )
\min
\{
\mathbf{u} \cdot \mathbf{m}
\mid
\mathbf{m} \in \ZZ_{\geq 0}^{r+1}
,
c_{\mathbf{m}} \neq 0
\}
.
\]
\end{Proposition}

We remark that in the above proposition,
the absolute value $|\cdot|_{\mathbf{u}}$
is independent of the choice of the Cohen isomorphism.

Since $\bfK (X)$ is the 
residue field of the generic point of $X$,
the absolute value $| \cdot |_{\mathbf{u}}$ 
in the above proposition gives us a point
of
$X^{\an}$.
We call it the \emph{monomial point associated to
$( \mathscr{X}, \xi , \mathbf{u})$}.
The monomial point associated to 
$( \mathscr{X}, \xi , \mathbf{u})$
for some stratum $S$ 
and $\mathbf{u} \in \Delta^{r}$ 
is called a monomial point with respect to $\mathscr{X}$.

\subsubsection{Skeleta (over CDV subfields)}

By Proposition~\ref{prop:monomialpoints},
we obtain an injective map
$\rho_S : \Delta^{r} \to \bfX^{\an}$
which maps $\mathbf{u} = (u_0 , \ldots , u_r )$ to the 
monomial point
given by
$| \cdot |_{\mathbf{u}}$.
One sees that this map is continuous.
Let $\Delta_S$ denote the image of $\Delta^{r}$ by this map.
Then
it is a compact subset of $X^{\an}$.
We set $S ( \mathscr{X} ) :=
\bigcup_{S \in \str ( \widetilde{\mathscr{X}} )} \Delta_S$.
This is the set of monomial points
with respect to $\mathscr{X}$.
It is called the \emph{skeleton} of $X^{\an}$
associated to $\mathscr{X}$.

The skeleton is a compact subset.
Furthermore,
it has a canonical structure of simplicial complex
which reflects the incidence relation of $\str ( \widetilde{\mathscr{X}} )$.
First,
we put a structure of a simplex on
 $\Delta_S$ 
via the homeomorphism $\rho : \Delta^r \to \Delta_S$.
Next, to see the incidence relation,
let $S$ and $S'$ be strata of $\widetilde{\mathscr{X}}$
and suppose $S' \subset \overline{S}$,
where $\overline{S}$ is the Zariski closure of $S$ in 
$\widetilde{\mathscr{X}}$.
Then we have $\Delta_{S'} \supset \Delta_S$
and $\Delta_S$ is a face of $\Delta_{S'}$.
Indeed,
let $\xi$ and $\xi'$ be the generic points of $S$ and $S'$,
respectively.
Let
$V_0 , \ldots , V_r$ be the irreducible components of 
$\widetilde{\mathscr{X}}$ such that
$S = V_0 \cap \cdots \cap V_r$ around $\xi$.
Since $S' \subset \overline{S}$,
there are irreducible components  $V_{r+1} , \ldots , V_{r'}$ of 
$\widetilde{\mathscr{X}}$
such that $S' = V_0 \cap \cdots \cap V_{r'}$ around $\xi'$.
Let $t_0 , \ldots , t_{r'} \in \OO_{\mathscr{X} , \xi'}$ be 
local equations of $V_0 , \ldots , V_{r'}$.
Then $t_0 , \ldots , t_{r'}$ form a 
regular system of parameters of $\OO_{\mathscr{X} , \xi'}$,
and
$t_0 , \ldots , t_{r}$ form that of $\OO_{\mathscr{X} , \xi}$.
Let $( u_0 , \ldots , u_r) \in \Delta^r$.
Note $(u_0 , \ldots , u_r , 0 , \ldots , 0) \in \Delta^{r'}$.
Let $\rho_S : \Delta^{r} \to \Delta_S$
and $\rho_{S'} : \Delta^{r'} \to \Delta_{S'}$ be the homeomorphism
defined above.
Then by the uniqueness assertion of Proposition~\ref{prop:monomialpoints},
one can show that
$\rho_S ( u_0 , \ldots , u_r) = \rho ( u_0 , \ldots , u_r , 0, \ldots , 0)$.
This means that $\Delta_{S'}$ contains $\Delta_S$ as a face.
To the contrary,
suppose 
that $\Delta_{S'} \supset \Delta_S$.
Then one can show that
$\Delta_{S}$ is a face of $\Delta_{S'}$
and 
$S' \subset \overline{S}$.
Thus
$S ( \mathscr{X} ) :=
\bigcup_{S \in \str ( \widetilde{\mathscr{X}} )} \Delta_S$
is a simplicial complex which reflects
the incidence relation of
$\str ( \widetilde{\mathscr{X}} )$.

This will not be really used in the sequel,
but let us give a remark when
$\dim (X) = 1$.
In this case, we have a notion of dual graph by configuration of the 
special fiber of $\mathscr{X}$.
The skeleton $S ( \mathscr{X} )$ is the realization of the 
dual graph inside the analytic space $X^{\an}$.

\subsubsection{Skeleta (over $\KK$)}

Now we consider skeleta over $\KK$.
Let $X$ be a smooth variety over $\KK$.
Assume that there exist
a CDV subfield $\bfK$
and
 a strictly semistable
model $ \mathscr{X} \to \Spec ( \bfK^{\circ} )$
of $X$.
The generic fiber $\mathscr{X}_{\bfK}$ is a smooth variety over $\mathbf{K}$.
Let $\rho : X^{\an} \to (\mathscr{X}_{\bfK})^{\an}$ be the natural 
surjective map.
Then one can show that for any $u \in S ( \mathscr{X} )
\subset (\mathscr{X}_{\bfK})^{\an}$,
i.e., a monomial point of $(\mathscr{X}_{\bfK})^{\an}$
with respect to $\mathscr{X}$,
$\rho^{-1} ( u )$ consists of one point.
It follows that $\rho^{-1} ( S ( \mathscr{X} ) )$
is a compact subset of $X^{\an}$ and has a structure of simplicial complex.
Furthermore,
one can show that this does not depend on the choice of $\bfK$
but depends only on the model $\mathscr{X}_{\KK^{\circ}}$
over $\KK^{\circ}$.
We denote this set by $S ( \mathscr{X}_{\KK} )$
and call it the \emph{skeleton}
of $X^{\an}$ associated to $\mathscr{X}$.
This skeleton reflects the incidence relation of
the special fiber of $\mathscr{X}_{\KK^{\circ}}$.

We can actually define the skeleton of the analytic space
associated to any strictly semistable model over $\KK$,
but we do not do that,
because
in our use,
the strictly semistable model can be defined over a CDV subring.

\subsection{Metrics on line bundles}
In this subsection, we recall basic notions of metrics on line bundles.
We refer to \cite{zhang1-2}.

Let $V$ be a $1$-dimensional vector space over $\KK$.
A \emph{metric $|| \cdot ||_{V}$ on $V$} 
is a nontrivial function $V \to \RR_{\geq 0}$
such that $|| \lambda v ||_V = |\lambda| || v ||_V$
for any $\lambda \in \KK$,
where $| \lambda |$ is the absolute value of $\lambda$.

Let $X$ be an algebraic variety over $\KK$
and let $L$ be a line bundle on $X$.
A metric on $L$ is a collection $\{ || \cdot ||_{L (x)} \}_{x \in X( \KK )}$
of metrics on the $1$-dimensional
$\KK$-vector spaces $L(x)$ for all $x \in X( \KK )$,
where $L(x) = L \otimes_{\OO_{X}} \kappa (x ) $
denotes the fiber of $L$
over the $\KK$-valued point $x$.

Among the metrics on line bundles,
there are two very important classes of metrics:
the algebraic metrics and the semipositive metrics.
The algebraic metrics are determined from
models over $\KK^{\circ}$,
and the semipositive metrics are 
the ``uniform limits'' of algebraic metrics 
arising from vertically nef models.

First, we recall the notion of algebraic metric.
Let $X$ be a projective variety over $\KK$.
Let $( \mathscr{X} \to \Spec ( \KK^{\circ}) , \mathscr{L} )$
be a proper flat model of $( X, L)$.
Take any point $x \in X( \KK )$.
By the valuative criterion of properness,
$x$ extends to a section $\sigma_{x} : \Spec ( \KK^{\circ} )
\to \mathscr{X}$.
There exists an open neighborhood 
$\mathscr{U}$
of the image of $\sigma_{x}$
with a trivialization
$\varphi : \mathscr{L}|_{\mathscr{U}} \to 
\OO_{\mathscr{U}}$.
Restricting this trivialization to the generic fibers,
we obtain an isomorphism
$\varphi_{U} : \mathscr{L}|_{U} \to 
\OO_{U}$,
where $U$ is the generic fiber of $\mathscr{U}$.
For any $s(x) \in L(x)$, $\varphi_U (s(x))$ is an element of $\OO_{U} (x) =
\KK$, and
we set $|| s(x) ||_{\mathscr{L}(x)} := | \varphi_U (s(x)) |_{\KK}$.
It is not difficult to
see that $|| s(x) ||_{\mathscr{L}(x)}$
does not depend on any choices other than the model $\mathscr{L}$.
Thus we obtain a metric $|| \cdot ||_{\mathscr{L}}
=
\{ || \cdot ||_{\mathscr{L}(x)}\}_{x \in X(\mathbb{K})}$
on $L$.
This metric is called the \emph{algebraic
metric arising from $(\mathscr{X} , \mathscr{L})$}.

We sometimes consider a model $( \mathscr{X} , \mathscr{L})$
of $(X , L)$ over a CDV subring and consider the
algebraic metric arising from
$( \mathscr{X}_{\KK^{\circ}} , \mathscr{L}_{\KK^{\circ}})$,
the base-change of $( \mathscr{X} , \mathscr{L})$ to $\KK$
over the CDV subring.
To ease notation,
this algebraic metric is denoted by $|| \cdot  ||_{\mathscr{L}}$
instead of $|| \cdot  ||_{\mathscr{L}_{\KK^{\circ}}}$.
We say that such an
algebraic metric \emph{can be
defined over a CDV subring}.

A line bundle $\mathscr{L}$ on $\mathscr{X}$
is said to be \emph{vertically nef} if
$\widetilde{\mathscr{L}} := 
\mathscr{L}|_{\widetilde{\mathscr{X}}}$
is nef.
A metric that is of the form $|| \cdot ||_{\mathscr{L}}$
for some vertically nef $\mathscr{L}$
is said to be \emph{semipositive}.

The notion of algebraic metric
is generalized to the notion of $\QQ$-algebraic metric.
Let $L$ be a line bundle on a projective variety over $\KK$.
Let $|| \cdot ||$ be a metric on $L$.
We call $|| \cdot ||$ a \emph{$\QQ$-algebraic metric}
if there exists a positive integer $N$
such that $||  \cdot ||^{\otimes N}$,
which is a metric on $L^{\otimes N}$,
is an algebraic metric.
Furthermore,
if $N$ can be taken in such a way that
$||  \cdot ||^{\otimes N}$ can be defined over a CDV subring,
we say that 
$|| \cdot ||$ 
can be defined over a CDV subring.
A $\QQ$-algebraic metric $|| \cdot ||$
is said to be \emph{semipositive}
if $|| \cdot ||^{\otimes N}$ is an algebraic 
and semipositive metric
for some $N \geq 1$.

We would like to define the notion of uniform limit of a sequence of metrics.
Before that, we define a function arising
from two metrics on a line bundle.
Let $|| \cdot ||$ and $|| \cdot ||'$
be metrics on a line bundle $L$.
We define a function
$
- \log ( || \cdot ||/ || \cdot ||')
$
on $X (\KK)$
by 
\[
- \log \frac{|| \cdot ||}{|| \cdot ||'} (x) :=
- \log \frac{|| s(x) ||}{|| s (x) ||'}
,
\]
where $s$ is a local section of $L$ that is non-zero at $x$.
This does not depend on the choice of $s$ and gives a well-defined function
on $X (\KK)$.
Now, let $|| \cdot ||$ be a metric on $L$
and let $( || \cdot ||_n )_{n \in \NN}$ be a sequence of metrics on $L$.
We say that 
$( || \cdot ||_n )_{n \in \NN}$ \emph{uniformly converges to
$|| \cdot ||$}
if 
$- \log ( || \cdot ||/ || \cdot ||_n )$
is a bounded function on $X ( \overline{K} )$
and if
\[
\lim_{n \to \infty}
\sup_{x \in X ( \KK )}
\left\{  - \log \frac{|| \cdot ||}{|| \cdot ||_n} (x) 
\right\}
= 0
.
\]
If this is the case,
we call $|| \cdot ||$
is the \emph{uniform limit}
of $( || \cdot ||_n )_{n \in \NN}$.

Let $(L , || \cdot ||)$ be a metrized line bundle.
We say that $|| \cdot ||$
is \emph{semipositive},
if it is a uniform limit of a sequence of semipositive $\QQ$-algebraic
metrics.

\begin{Lemma} \label{lem:functorialsemipositive}
Let $f : X \to Y$ be a morphism of projective variety over $\KK$.
Let $\overline{L}$ and $\overline{M}$ be line bundles on $X$ and $Y$
with semipositive metrics, respectively.
Then the metric of
$\overline{L} \otimes f^{\ast} ( \overline{M} )$
is semipositive.
\end{Lemma}

Let us give some comments.
In the above, we consider metrics over
$X ( \KK )$,
but in fact,
we can do that over the analytic space $X^{\an}$.
Indeed, we can define algebraic metrics as a metric
on line bundles on $X^{\an}$.
The algebraic metrics are continuous with respect to
the topology on $X^{\an}$.
It follows that
the $\QQ$-algebraic metrics are also
continuous metrics.
Furthermore,
a
semipositive metric is actually defined to be the uniform limit
of a sequence
of $\QQ$-algebraic metrics over $X^{\an}$
in the usual sense.
Therefore a semipositive metric is continuous on $X^{\an}$
in the usual sense.
We refer to \cite[\S~7]{gubler-1} for details.

\subsection{Chambert-Loir measures} \label{subsect:chamberloirmeasure}

Let $X$ be a projective variety over $\KK$
of dimension $d$.
Let $\overline{L} = ( L , || \cdot ||)$ be a line bundle on $X$
equipped
with a semipositive metric.
To these data, we can associate a regular Borel measure 
$\cherncl_1 ( \overline{L})^{\wedge d}$,
which we are going to explain.

In this subsection, we assume that $X$ can be defined over a
CDV subfield.
Further,
we only considers semipositive metrics
which are the limit of $\QQ$-algebraic metrics
that can be defined over a CDV subring.

First, we consider the case where the metric $|| \cdot ||$ is algebraic.
By definition, there exist
a discrete valued subfield $\bfK$
and a proper flat model
$( \mathscr{X} , \mathscr{L})$ 
of $(X,L)$
over $\bfK^{\circ}$
such that
$|| \cdot || = || \cdot ||_{\mathscr{L}}$.
Further, replacing $\bfK^{\circ}$
by a finite extension,
we can take
a proper flat morphism 
$\varphi : \mathscr{X}' \to \Spec (\bfK^{\circ})$
with reduced special fiber
and a generically finite surjective morphism 
$\nu : \mathscr{X}' \to \mathscr{X}$ over $\bfK^{\circ}$.
(In fact, de Jong's semistable alteration theorem will
give us such $\varphi : \mathscr{X}' \to \Spec (\bfK^{\circ})$;
see \S~\ref{subsect:SA} in the sequel.)
Set $\mathscr{L}' := \nu^{\ast} ( \mathscr{L} )$.
Put $X' := \mathscr{X}' \otimes_{\bfK^{\circ}}
\bfK$ and $X'_{\KK} := 
X' \otimes_{\bfK} \KK$.
Then the metric of the pull-back
$\overline{L'} := \nu|_{X'_{\KK}}^{\ast} ( \overline{L} )$
is 
the algebraic metric
associated to the model $( \mathscr{X}' , \mathscr{L}' )$.
Let $\Irr (\widetilde{\mathscr{X}'})$ be the set of
irreducible component of the special fiber of
$\mathscr{X}'$.
For each $V \in \Irr (\widetilde{\mathscr{X}'})$,
let $\xi_V$ denote the generic point of $V$.
Since $\varphi$ is flat and $\widetilde{\mathscr{X}'}$ is reduced at
$\xi_V$,
$\mathscr{X}'$ is normal at $\xi_V$.
Therefore we have a
corresponding point
$\xi_V^{\an} \in (X')^{\an}$.
As we noted in \S~\ref{subsect:DP},
there exists a unique point $\eta_V$ 
which maps to $\xi_V^{\an}$ by the canonical map $(X'_{\KK})^{\an}
\to (X')^{\an}$.
With those notation,
we define a measure $\cherncl_1 ( \overline{L'} )^{\wedge d}$
on $(X'_{\KK})^{\an}$
to be
\[
\sum_{V \in \Irr (\widetilde{\mathscr{X}'})}
\deg_{\mathscr{L}'} (V)
\delta_{\eta_V}
,
\]
where 
$\delta_{\eta_V}$
is the Dirac measure on $\eta_V$.
It is a regular Borel measure
with total mass
$\deg_{L'} (X')$.
To define a regular Borel measure
$\cherncl_1 ( \overline{L} )^{\wedge d}$
on $X^\an$,
let $\nu^{\an} : (X'_{\KK})^{\an} \to X_{\an}$ be the morphism
of analytic spaces associated to $\nu$.
Consider the pushout
$\nu^{\an}_{\ast} ( \cherncl_1 ( \overline{L'} )^{\wedge d} )$.
Then one shows that
this does not depend on the choice of $\mathscr{X}'$ or $\nu$,
and hence it is well defined for $(X , \overline{L})$.
We define $\nu^{\an}_{\ast} ( \cherncl_1 ( \overline{L'} )^{\wedge d} )$
to be $\cherncl_1 ( \overline{L} )^{\wedge d}$.
This measure also
has total mass $\deg_{L} (X)$.
It is
called the \emph{Chambert-Loir measure}
of $\overline{L}$.
It is originally defined by Chambert-Loir \cite{chambert-loir}
and is also developed by Gubler \cite[Section~3]{gubler1}

\begin{Remark}
If one works over the framework of admissible formal schemes,
one can define the Chambert-Loir measure 
without using the generically
finite morphism $\nu$ as above.
In fact, one can take a proper
flat model of $(X , \overline{L})$ 
in the category of admissible formal schemes with reduced fiber,
and using those models will lead to the definition of the 
Chambert-Loir measure.
\end{Remark}

Next assume that $\overline{L}$ is a $\QQ$-algebraically
metrized line bundle on $X$
and let $N$ be a positive integer
such that
$\overline{L}^{\otimes N}$ 
is a line bundle with an algebraic metric associated to a model
that can be defined over a CDV ring.
In this case,
we define
\[
\cherncl_1 ( \overline{L} )^{\wedge d} :=
\frac{1}{N^d} \cherncl_1 ( \overline{L} ^{\otimes N})^{\wedge d}
.
\]
This is a regular Borel measure of total mass $\deg_L (X)$.

Finally, 
let $|| \cdot ||$ be a semipositive metric.
We take a sequence $( || \cdot ||_{n})_{n \in \NN}$
of semipositive $\QQ$-algebraic metrics that can be define over a CDV subring
such that $|| \cdot ||$
is the uniform limit of $( || \cdot ||_{n})_{n \in \NN}$.

\begin{Proposition} [cf. \cite{chambert-loir}] \label{prop:functorialprop}
With the notation above,
the sequence of regular Borel measures
$( \cherncl_1 ( L , || \cdot ||_n)^{\wedge d} )_{n \in \NN}$
weakly converges to a regular Borel measure.
Furthermore, the weak limit does not depend on the
choice of the sequence $( || \cdot ||_{n})_{n \in \NN}$
and depends only on $|| \cdot ||$.
\end{Proposition}

By the above proposition,
we define the \emph{Chambert-Loir measure} 
$\cherncl_1 ( \overline{L} )^{\wedge d}$
of $\overline{L}
= ( L , || \cdot ||)$
to be the weak limit of $( \cherncl_1 ( L , || \cdot ||_n)^{\wedge d} )_{n \in \NN}$.
This has total mass $\deg_L (X)$.

\subsection{Nonarchimedean canonical metrics and canonical measures}

Let $A$ be an abelian variety $\KK$.
Let $L$ be an even line bundle.
A \emph{rigidification of $L$}
is an isomorphism $L ( 0 ) \cong \KK$,
where $0$ is the zero element of $A$.

For a natural number $m > 1$,
there exists an isomorphism
$\phi : [m]^{\ast} (L) \to
L^{\otimes m^2}$,
where $[m]: A \to A$
is the $m$ times endomorphism.
Once we fix a rigidification $L (0) = \KK$
of $L$,
the isomorphism which respects the rigidification is unique.
Indeed, then $[m]^{\ast} (L)$ has the rigidification 
$[m]^{\ast} (L) (0) \cong \KK$
by pullback,
and $L^{\otimes m^2}$ has the rigidification 
$L^{\otimes m^2} (0) \cong \KK$
by tensor product,
so that we have an isomorphism
\[
\begin{CD}
\KK @>{\cong}>> [m]^{\ast} (L) (0) 
@>{\phi}>> L^{\otimes m^2} @>{\cong}>> \KK
.
\end{CD}
\]
We have a unique choice of $\phi$
in such a way that the above isomorphism is an identity.
That is an isomorphism $\phi$ which respects the rigidification.

\begin{Proposition} [cf. Theorem~9.5.4 of 
\cite{BG}] \label{prop:nonarchimedeancanoncalmeasure}
Let $A$ and $L$ be as above.
Fix a rigidification of $L$.
Then
there exists a unique metric $|| \cdot ||$
on $L$
with the following property.
For any $m$,
let $\phi : [m]^{\ast} L \to
L^{\otimes m^2}$
be the isomorphism which respects the rigidification.
Then $\phi$ induces an isometry
$[m]^{\ast} ( L , || \cdot ||) \cong 
( L^{\otimes m^2 } , || \cdot ||^{\otimes m^2 })$
of metrized line bundles.
Furthermore,
suppose that $L$ is ample.
Then
this metric is semipositive.
\end{Proposition}

We call the metric in this proposition
the \emph{canonical metric of $L$ (with respect to the rigidification)}.
The canonical metric depends on the choice of the rigidification,
but the difference is only up to positive constant multiple.

We do not give a precise proof of the above theorem,
but when $L$ is ample, we briefly describe how the 
semipositive canonical metric
is constructed on $L$
under the following assumption:
the abelian variety $A$ and the even ample line bundle $L$ 
can be defined over a
CDV subfield $\bfK$;
further,
the rigidification 
we consider can be defined over $\bfK$.
(Suppose $A = A' \otimes_{\bfK} \KK$
and $L = L' \otimes_{\bfK} \KK$.
Then we say that the rigidification $L ( 0 ) \cong \KK$
can be defined over $\bfK$ if it is the base-change
of an isomorphism $L' (0) \cong \bfK$.)
We construct the metric by making a sequence
of semipositive 
$\QQ$-algebraic metrics
by induction.
We start with any proper flat model $( \mathscr{A}_1 , \mathscr{L}_1)$
over $\bfK^{\circ}$
such that $\mathscr{L}_1$ is vertically nef.
We can take such a model because $L$ is ample.
Fix a natural number $m > 1$.
The $m$-times homomorphism $[m] : A \to A$
does not extend to $\mathscr{A}_1 \to \mathscr{A}_1$ in general,
but if we take a suitable proper model $\mathscr{A}_2$ of $A$,
then $[m] : A \to A$
extends to a morphism $f_1 : \mathscr{A}_2 \to \mathscr{A}_1$.
Set $\mathscr{L}_2 := f_1^{\ast} ( \mathscr{L}_1) $
The pair
$( \mathscr{A}_2, \mathscr{L}_2)$
is a model of $(A , [m]^{\ast} ( L) )$.
Since $\mathscr{L}_1$ is vertically nef,
so is $\mathscr{L}_2$.
Since we have the isomorphism $[m]^{\ast} ( L ) \cong L^{\otimes m^{2}}$
which respects the rigidification,
that model is naturally regarded as a model of 
$( A , L^{\otimes m^{2}})$.
This gives us an algebraic metric on $L^{\otimes m^2}$
and hence defines a
$\QQ$-algebraic metric 
$|| \cdot ||_2$ on $L$.
Since $f^{\ast} ( \mathscr{L}_1 )$ is vertically nef,
$|| \cdot ||_2$ is semipositive.

Next, we take a model $\mathscr{A}_3$ with a
morphism $f_2 : \mathscr{A}_3 \to \mathscr{A}_2$
extending $[m] : A \to A$.
Set $\mathscr{L}_3 := f_2^{\ast} ( \mathscr{L}_2)$,
which is vertically nef, 
and consider $( \mathscr{A}_3 , \mathscr{L_3})$.
Then via the isomorphism $[m]^{\ast} ( L ) \cong L^{\otimes m^{2}}$,
this is a model of $(A , L^{\otimes m^{4}})$.
Thus this defines a semipositive $\QQ$-algebraic metric $|| \cdot ||_3$
on $L$.
Repeating this process,
we obtain a sequence 
\addtocounter{Claim}{1}
\begin{align} \label{eq:sequenceofmodelsforcanonial}
\left(
\pi_n : \mathscr{A}_n 
\to \Spec ( \mathbf{K}^{\circ})
, \mathscr{L}_n
\right)_{n \in \NN}
\end{align}
whose $n$th term is a model of $(A , L^{\otimes m^{2(n -1)}})$
and which gives a
sequence of
semipositive $\QQ$-algebraic metrics $(|| \cdot ||_n)_{n \in \NN}$
on $L$.
One can show that this sequence converges to a metric
$|| \cdot ||$ on $L$.
It is semipositive,
and furthermore it has the required property.

\begin{Remark} \label{rem:canonicalmetricnondeg}
Suppose that the abelian variety $A$ over $\KK$ is nondegenerate.
Then the canonical metric is an algebraic metric.
Indeed, fix a rigidification $L (0) = \KK$.
Since $A$ is nondegenerate, we take a model
$(\mathscr{A} , \mathscr{L} )$ of $(A,L)$ such that $\pi : \mathscr{A}
\to \Spec ( \KK^{\circ} )$ is an abelian scheme.
Let $0_\pi$ be the zero-section of $\pi$.
Then the identity $L (0) = \KK$
extends to an isomorphism $0_{\pi} (\mathscr{L}) \otimes  \mathscr{N} 
\cong \KK^{\circ}$ for some line bundle $\mathscr{N}$ on 
$\Spec ( \KK^{\circ} )$.
Replacing $\mathscr{L}$ by $\mathscr{L} \otimes \pi^{\ast} ( \mathscr{N} )$,
we assume that the rigidification extends to
$0_{\pi}^{\ast} ( \mathscr{L})
\cong \KK^{\circ}$.
Then for any $n \in \NN$,
we see that
$[n]^{\ast} ( \mathscr{L} ) \cong \mathscr{L}^{\otimes n^2}$.
It follows that the metric $|| \cdot ||_{\mathscr{L}}$
satisfies $[n]^{\ast} ( || \cdot ||_{\mathscr{L}} ) = 
|| \cdot ||_{\mathscr{L}}^{\otimes n^2}$ on $L$.
Thus the the canonical metric is an algebraic metric.
\end{Remark}

Finally,
we define the canonical measure.
Let $A$ and $L$ be as above.
Fix a rigidification of $L$,
and let $\overline{L}$ be the line bundle $L$ with the canonical metric.
Let $X$ be a closed subvariety of $A$.
Since the canonical metric is semipositive,
the metric of the restriction $\overline{L}|_{X}$
is also semipositive by Proposition~\ref{prop:functorialprop}.
It follows that
we have the Chambert-Loir measure
$
\cherncl_1 ( \overline{L}|_{X} )^{\wedge d}
$
on $X^{\an}$.
By the semipositivity of the canonical metric,
this is a positive measure.
We can show 
that this measure does not depend on the choice of the rigidification,
and hence it is well defined for $L$.
One often calls it the canonical measure,
but
in this paper, we call the probability measure
\[
\mu_{X_v^{\an} , L} := \frac{1}{\deg_L (X)}
\cherncl_1 ( \overline{L}|_{X} )^{\wedge d}
,
\]
given by normalization,
the \emph{canonical measure
on $X^{\an}$ 
associated to $L$}.

\begin{Remark} \label{rem:canonicalmeasurenondeg}
Suppose that the abelian variety $A$ over $\KK$ is nondegenerate.
Then since the canonical metric is an algebraic metric
(cf. Remark~\ref{rem:canonicalmetricnondeg}),
the canonical measure is a linear combination of Dirac measures
of points.
\end{Remark}

\section{Proof of Gubler's theorem} \label{sect:proofofgubler}

Recall that in the proof of Zhang over number fields,
analytic spaces over an archimedean place,
canonical measures on analytic spaces,
and
equidistribution theorem
were key ingredients.
Over function fields,
we have 
(candidates for) the counterparts of the first two items;
they are Berkovich analytic spaces over a non-archimedean place
and the canonical measures on them.

However,
we have not yet had that of the third.
Indeed,
the equidistribution theorem on Berkovich spaces
had not been established yet
at the time when Gubler proved his theorem.
In such a situation,
He considered the tropicalization of closed subvarieties of 
totally degenerate abelian varieties and canonical measures,
and he established the equidistribution theorem
over tropical varieties with respect to the tropical canonical measure,
instead.

In the sequel, 
we follow \cite[\S~A]{gubler1}
for terminology of basic notions on convex geometry.

\subsection{Uniformization and tropicalization of totally degenerate abelian varieties}
In this subsection,
let $\mathbb{K}$ be an algebraically closed field that is complete
with respect to a non-trivial non-archimedean value $| \cdot |$.
(The most important example in this paper is $\overline{K}_v$.)
The subset $\{ - \log |a| \mid a \in \mathbb{K}^{\times} \}$
is called the value group of $| \cdot |$.
For simplicity, we assume that the value group equals $\mathbb{Q}$.
For an algebraic variety $X$ over $\mathbb{K}$,
let $X^{\an}$ denote the associated Berkovich analytic space.

Let $A$ be an abelian variety over $\mathbb{K}$
of dimension $n$.
Assume that $A$
is totally degenerate.
Then
there exists a homomorphism of group analytic spaces
$p : ( \mathbb{G}_{m }^{n} )_{\mathbb{K}}^{\an} \to A^{\an}$ 
such that $\Ker ( p )$ is a
free abelian group of rank $n$.
(That is an alternative definition of totally degenerate abelian 
varieties.)
This $p$ is the universal covering of $( A^{\an} , 0 )$.
Let $x_1 , \ldots , x_n$ be the standard coordinates of $\mathbb{G}_m^{n} $.
For any $x_i$ ($i= 1, \ldots , n$)
and
$P = ( \iota (P) , | \cdot | ) \in ( \mathbb{G}_m^{n} )_{\mathbb{K}}^{\an}$
(see \S~\ref{subsect:berkovich} for the notation),
write $|x_i ( P )| := | x_i( \iota (P)) |$
where $x_i( \iota (P))$ the image of $x_i$ in the residue field
$\kappa ( \iota (P) )$ of $\mathcal{O}_{( \mathbb{G}_m^{n} )_{\mathbb{K}}}$ at 
$\iota ( P )$.
We consider the map
$\val : ( \mathbb{G}_m^{n} )_{\mathbb{K}}^{\an} \to \RR^{n}$ defined
by
\addtocounter{Claim}{1}
\begin{align}
\label{align:def:val}
\val ( P ) = ( - \log |x_1 (P)| , \ldots , - \log | x_n (P) | )
.
\end{align}
Then one shows that $\Lambda := \val (\Ker ( p ))$ is a complete lattice
of $\RR^{n}$.
Thus $\val$ descends to a homomorphism
$\overline{\val} : A^{\an} \to \RR^{n} / \Lambda$.
This homomorphism called the \emph{tropicalization map} of $A^{\an}$.
To summarize,
we obtain the following commutative diagram
in which all the maps are continuous homomorphisms:
\begin{align*}
\begin{CD}
( \mathbb{G}_{m }^{n} )_{\mathbb{K}}^{\an} @>>> A^{\an}
\\
@VVV @VVV
\\
\RR^{n} @>>> \RR^{n} / \Lambda
.
\end{CD}
\end{align*}

The tropicalization maps are compatible with homomorphisms of abelian 
varieties and direct products.
To be precise,
let $A_1$ and $A_2$ be an abelian varieties over $\mathbb{K}$,
and
suppose that $A_1$ and $A_2$ are totally degenerate.
Let $\phi : A_1 \to A_2$ be a homomorphism of abelian varieties.
Let $p_1 : ( \mathbb{G}_{m }^{n_1} )_{\mathbb{K}}^{\an} \to A_{1}^{\an}$ 
and $p_2 : ( \mathbb{G}_{m }^{n_2} )_{\mathbb{K}}^{\an} \to A_{2}^{\an}$ 
be the uniformizations.
Then $\phi$ lifts to a homomorphism 
$( \mathbb{G}_{m }^{n_1} )_{\mathbb{K}}^{\an} \to
( \mathbb{G}_{m }^{n_2} )_{\mathbb{K}}^{\an}$.
Furthermore, this induces a linear map
$\RR^{n_1} \to \RR^{n_2}$,
which descends to a homomorphism $\phi_{\aff} : \RR^{n_1} / \Lambda_1 \to 
\RR^{n_2} / \Lambda_2$ of real tori.
Further, consider the direct product $A_1 \times A_2$
and let $\overline{\val} : (A_1 \times A_2)^{\an} \to \RR^{n} / \Lambda$
be the tropicalization.
Then
we have $\RR^{n} / \Lambda = \RR^{n_1} / \Lambda_1
\times \RR^{n_2} / \Lambda_2$

The following lemma is used in Gubler's proof.

\begin{Lemma} \label{lem:techttd}
Let $A_1$ and $A_2$ be an abelian varieties over $\mathbb{K}$.
Let $\phi : A_1 \to A_2$ be a homomorphism of abelian varieties.
Suppose that $A_1$ is totally degenerate and that
$\phi$ is surjective.
Then $A_2$ is totally degenerate. 
\end{Lemma}

\Proof
See \cite[Lemma~6.1]{gubler2}.
\QED

The following theorem is a fundamental theorem
of tropical analytic geometry.

\begin{Theorem} [cf. Theorem~1.1 of \cite{gubler1}] \label{thm:tropicalization1}
Let $A$ be a totally degenerate abelian variety 
over $\KK$
and
let $X$ be a closed subvariety of $A$ of dimension $d$.
Then the image
$\overline{\val} (X^{\an})$ is a rational 
polytopal subset of $\RR^{n} / \Lambda$
of pure dimension $d$.
\end{Theorem}

We call $\overline{\val} (X^{\an})$ the \emph{tropicalization of $X^{\an}$}
or the \emph{tropical variety associated to $X$.}

\subsection{Tropical canonical measure}

From here on to the end of this section, let $K$ be a function field
unless otherwise specified.
Let $A$ be an abelian variety over $\overline{K}$.
Let $v$ be a place of $\overline{K}$.
Assume that $A_{v} := A \otimes_{\overline{K}} \overline{K}_v$
is totally degenerate.
Let $\overline{\val} : A_v^{\an} \to \RR^{n} / \Lambda$ be the tropicalization
map of $A_v$, where $n := \dim (A)$.

Let $X$ be a closed subvariety of $A$ of dimension $d$.
We consider the canonical measure $\mu_{X_v^{\an}}$ on $X_v^{\an}$
of an even line bundle on $A$.
The pushout $\overline{\val}_{\ast} ( \mu_{X_v^{\an}} )$ is a measure
on $\overline{\val} ( X_v^{\an})$,
which is called a
\emph{tropical canonical measure}.
We write $\mu_{X_v^{\an}}^\trop :=
\overline{\val}_{\ast} ( \mu_{X_v^{\an}} )$.

The following theorem describes the tropical canonical measures.

\begin{Theorem} \label{thm:tropicalcanonicalmeasure}
Let $A$, $v$, $\overline{\val} : A_v^{\an} \to \RR^{n} / \Lambda$,
and
$X$
be
as above.
Let $L$ be an even line bundle on $A$
and
let $\mu_{X_v^{\an}}$ be the canonical measure associated to
$L$.
Then
there exist rational simplices $\Delta_1 , \ldots , \Delta_m$
of $\RR^{n} / \Lambda$ of dimension $d$ 
and positive real numbers $r_1 , \ldots , r_m$
such that
$\relin ( \Delta_i) \cap \relin ( \Delta_j ) = \emptyset$ for 
$i \neq j$,
$\overline{\val} (X_v^{\an}) =
\bigcup_{i = 1}^{m} \Delta_i$,
and such that
\[
\mu_{X_v^{\an}}^\trop
= 
\sum_{i = 1}^{m} r_i \delta_{\Delta_i}
\]
where $\relin ( \Delta_i)$ is the relative interior of $\Delta_i$
and
$\delta_{\Delta_i}$ is the relative Lebesgue measure on $\Delta_i$.
Furthermore, if $L$ is ample,
then $r_i > 0$ for all $i = 1 , \ldots , m$.
\end{Theorem}

\begin{Remark} \label{rem:positivityofthesiplex}
Let $\Delta$ be a simplex with $\Delta \subset 
\overline{\val} (X_v^{\an})$.
By Theorem~\ref{thm:tropicalcanonicalmeasure},
if $\dim ( \Delta ) = d$, then
$\mu_{X_v^{\an}}^\trop ( \Delta ) > 0$;
if $\dim (\Delta) < d$,
then $\mu_{X_v^{\an}}^\trop ( \Delta ) = 0$.
\end{Remark}

\subsection{Tropical equidistribution theorem} \label{subsect:ted}

Let $A$ be an abelian variety over $\overline{K}$
and let $X$ be a closed subvariety of $A$.
Let $K'$ be a finite extension of $K$ 
over which $A$ and $X$ can be defined.
Let $\Aut_{K'} ( \overline{K} )$
be the group of automorphisms of $\overline{K}$
over $K'$.
For a point $x \in X ( \overline{K} )$,
let $O (x)$ denote the $\Aut_{K'} ( \overline{K} )$-orbit of $x$ in 
$X ( \overline{K} )$.
Let $v \in  M_{\overline{K}}$ be a place
and regard $X ( \overline{K} ) \subset X_v^{\an}$.

\begin{Theorem} [cf. Theorem~5.5 of \cite{gubler3}] \label{thm:tropicalequidistribution}
Let $A$,
$X$, and $O(x)$ for $x \in X ( \overline{K} )$ be as above.
Let $( x_i )_{i \in I}$ be a generic small net on $X ( \overline{K} )$.
Let $\mu_{X_v^{\an}}$ be the canonical measure 
on $X_v^{\an}$
of an even ample line bundle on $A$.
Assume that $A_v$ is totally degenerate.
Let $\overline{\val} : A_v^{\an} \to \RR^{n} / \Lambda$
be the tropicalization map.
Then we have a weak convergence
\[
\lim_{i}
\frac{1}{\# O (x_i)}
\sum_{z \in O (x_i)} \delta_{\overline{\val} (z)}
=
\mu_{X_v^{\an}}^\trop
\]
of regular Borel measures on $\overline{\val} (X_v^{\an} )$,
where $\delta_z$ is the Dirac measure on $z$.
\end{Theorem}

\subsection{Proof of Gubler}

Let us start the proof of Gubler.
To argue by contradiction, suppose that
there exist an abelian variety $A$ over $\overline{K}$
that is totally degenerate at some place $v$
and a non-torsion subvariety $X$ of $A$
with dense small points.
Let $G_X$ be the stabilizer of $X$.
Then $A/G_X$ is totally degenerate at $v$ 
by Lemma~\ref{lem:techttd},
and $X/G_X$ has dense small points
by Lemma~\ref{lem:density-quotient}.
Replacing $A$ and $X$ with $A/G_X$ and $X/G_X$
respectively,
we may assume that $X$ is a non-torsion subvariety of $A$,
has dense small points, and has trivial stabilizer.
Since the theorem holds for $0$-dimensional subvarieties
(cf. Remark~\ref{rem:trivialcaseforGBC}~(1)),
we have $d := \dim (X) > 0$.

For a positive integer $N$,
let
$\alpha_N : X^{N} \to A^{N-1}$ 
be the difference morphism,
given by
$(x_1 , \ldots , x_N ) \mapsto (x_1 - x_2 , \ldots , x_{N-1} - x_{N} )$.
Set $Z:= X^{N} \subset A^{N}$ and 
$Y := \alpha_N (Z)$.
Since $X$ has trivial stabilizer,
there exists an $N$ such that
the restriction $\alpha : Z \to Y$
is generically finite
by Lemma~\ref{lem:genericallyfiniteN}.

Note that
$(A^{N})_v$ and $ (A^{N-1})_v$ 
are totally degenerate.
Let $\overline{\val}_1 : (A^{N})_v^{\an} \to \RR^{Nn} / \Lambda_1$
and $\overline{\val}_2 : (A^{N-1})_v^{\an} \to \RR^{(N-1)n} / \Lambda_2$
be the tropicalization maps,
where $n := \dim (A)$.
The induced affine homomorphism $\RR^{Nn} / \Lambda_1
\to \RR^{(N-1)n} / \Lambda_2$ 
restricts to a surjective piecewise affine map 
$\overline{\alpha}_{\aff}
: \overline{\val}_1 ( Z_v^{\an} ) \to \overline{\val}_2 ( Y_v^{\an} )$
between $Nd$-dimensional polytopal sets.
Let $\mu_{Z_v^\an}$ and $\mu_{Y_v^\an}$ be the canonical measures
on $Z_v^\an$ and $Y_v^\an$
of even ample line bundles,
respectively.
We put $\mu_{Z_v^\an}^{\trop} := (\overline{\val}_1)_{\ast} ( \mu_{Z_v^\an} )$
and $\mu_{Y_v^\an}^{\trop} := ( \overline{\val}_2)_{\ast} ( \mu_{Y_v^\an} )$.
Since $X$ has dense small points, so has $Z$.
By Lemma~\ref{lem:genericsmallnets},
there exists a generic small net $(z_i)_{i \in I}$ of $Z ( \overline{K} )$.
Since $\alpha : Z \to Y$ is surjective,
$(\alpha (z_i))_{i \in I}$ is a generic net on $Y ( \overline{K})$,
and it follows from
Proposition~\ref{prop:propertycanonicalheight}~(2)
that $(\alpha (z_i))_{i \in I}$ is small.
Since we have the equidistribution theorem
Theorem~\ref{thm:tropicalequidistribution},
the same argument as Zhang's proof
using the equidistribution theorem
gives us
$(\overline{\alpha}_{\aff})_{\ast} (  \mu_{Z_v^\an}^{\trop} ) = 
 \mu_{Y_v^\an}^{\trop} $.

Let $\overline{\val}_0 : A_v^{\an} \to \RR^{n} / \Lambda_0$
be the tropicalization map.
Since the tropicalization
is compatible with direct products,
we have 
$\RR^{Nn} / \Lambda_1 = (\RR^n / \Lambda_0)^N$
and
$\overline{\val}_1 ( Z_v^{\an})
= (\overline{\val}_0 ( X_v^{\an}))^{N}$.
Since $d = \dim (X) > 0$,
the diagonal of $\overline{\val}_1 ( Z_v^{\an})$ has positive dimension
(cf. Theorem~\ref{thm:tropicalization1}).
Since $\overline{\alpha}_{\aff}$ is linear and
contracts the diagonal to a point,
there exists an $Nd$-dimensional simplex $\Delta \subset
\overline{\val}_1 ( X_v^{\an})$ such that
$\dim (\overline{\alpha}_{\aff} ( \Delta ) ) < Nd$.
Since $
\mu_{Y_v^\an}^{\trop}$ is a piecewise relative Lebesgue
measure of dimension $Nd$,
it follows that $\mu_{Y_v^\an}^{\trop} 
( \overline{\alpha}_{\aff} ( \Delta )  )= 0$
(cf. Remark~\ref{rem:positivityofthesiplex}).
On the other hand,
since 
$(\overline{\alpha}_{\aff})_{\ast} (  \mu_{Z_v^\an}^{\trop} ) = 
 \mu_{Y_v^\an}^{\trop} $,
we have
\[
\mu_{Y_v^\an}^{\trop}
( \overline{\alpha}_{\aff} ( \Delta )  )
\geq
\mu_{Z_v^\an}^{\trop} ( \Delta)
.
\]
Since $\dim (\Delta) = \dim (Z) =Nd$,
the right-hand sides of this inequality
should be positive (cf. Remark~\ref{rem:positivityofthesiplex}).
However,
that is a contradiction.
Thus the proof is complete.

\subsection{Non-archimedean equidistribution theorem}

Gubler had to establish the tropical equidistribution theorem
(Theorem~\ref{thm:tropicalequidistribution}) because
no equidistribution theorem that could be applied to his setting was
known when he proved the theorem.
However,
inspired by Yuan's work
\cite{yuan},
he established in
\cite{gubler3}
a non-archimedean equidistribution theorem,
afterwards.

\begin{Proposition} [cf. Theorem~1.1 of \cite{gubler3}] \label{prop:ned}
Let $A$
be any abelian variety over $\overline{K}$
and let $X$ be a closed subvariety of $A$.
Let $O(x)$ for $x \in X ( \overline{K} )$ be as
in \S~\ref{subsect:ted}.
Let $v \in  M_{\overline{K}}$ be a place
and regard $X ( \overline{K} ) \subset X_v^{\an}$.
Let $\mu_{X_v^{\an}}$ be the canonical measure 
on $X_v^{\an}$
of an even ample line bundle on $A$.
Let $( x_i )_{i \in I}$ be a generic small net on $X ( \overline{K})$.
Then
we have a weak convergence
\[
\lim_{i} \frac{1}{\# O (x_i)}\sum_{z \in O(x_i)} \delta_z
= \mu_{X_v^{\an}}
,
\]
where $\delta_z$ is the Dirac measure on the point $z \in X_{v}^{\an}$.
\end{Proposition}

\subsection{Naive imitation of Zhang's proof} \label{subsect:imitateZ}

In the last subsection,
we have finally obtained candidates for counterparts of the three ingredients
of Zhang's proof.
Now,
let us 
try to apply Zhang's proof to the setting of
the geometric Bogomolov conjecture.
To show the conjecture by contradiction,
suppose that there exists a counterexample to the 
geometric Bogomolov conjecture for abelian varieties.
Then by the same argument as in the proofs of Zhang's theorem
and Gubler's theorem,
there exist an abelian variety $A$ over $\overline{K}$
and a non-special closed subvariety $X$ of $A$ 
such that $X$ has dense small points, $d := \dim (X) > 0$, and $X$ has trivial
stabilizer.
For a positive integer $N$,
we consider the difference homomorphism 
$\alpha_N : X^{N} \to A^{N-1}$,
and let $\alpha : Z \to  Y$
be the restriction
from $Z:= X^{N}$ to $Y := \alpha_N (Z)$.
For large $N$, $\alpha$ is generically finite 
by Lemma~\ref{lem:genericallyfiniteN},
as well as surjective.
Since $X$ has dense small points, so has $Z$.
Fix a $v \in M_{\overline{K}}$.
Let $\alpha^{\an} : Z_v^{\an} \to Y_v^{\an}$ be the associated
map between analytic spaces.
Since we have Proposition~\ref{prop:ned},
the same argument 
gives us 
\addtocounter{Claim}{1}
\begin{align} \label{eq:equalityofmeasures-nonarch}
\alpha^{\an}_{\ast} ( \mu_{Z_v^{\an}} ) = \mu_{Y_{v}^\an}
,
\end{align}
where $\mu_{Z_v^{\an}}$ and $\mu_{Y_{v}^\an}$ are the canonical
measures on $Z_v^{\an}$ and $Y_v^{\an}$ associated to 
even ample line bundles, respectively.

Then the problem is how to deduce a contradiction from the equality
(\ref{eq:equalityofmeasures-nonarch}).
In Gubler's totally degenerate setting, we can deduce
a contradiction from (\ref{eq:equalityofmeasures-nonarch}).
Indeed, assume that $A_v$ is totally degenerate.
Then the support of $\mu_{Z_v^{\an}}$
and that of $\mu_{Y_{v}^\an}$
have 
structures of polytopal set of pure dimension
dimension $Nd = \dim (Z) = \dim (Y)$.
Further,
for a suitable simplicial decompositions,
$\mu_{Z_v^{\an}}$ and $\mu_{Y_{v}^\an}$
are linear combinations of relative Lebesgue measures
of simplices of dimension $Nd$.
Roughly speaking, this says that
we have the same situation as Gubler's proof without
passing to the tropicalization.
Therefore,
we can get a contradiction by the same way.

However,
the equality
(\ref{eq:equalityofmeasures-nonarch})
does not necessarily lead to a contradiction,
in general.
For example, assume that $A_v$ is nondegenerate.
Then $A^{N}_v$ and $A^{N-1}_v$ are also nondegenerate.
In the nondegenerate case,
the canonical measures are
linear combinations of Dirac measure of points,
as is noted in Remark~\ref{rem:canonicalmeasurenondeg}.
Thus there is nothing strange with the equality
(\ref{eq:equalityofmeasures-nonarch}).

That observation suggests that
if $A$ is nowhere degenerate, then the equidistribution
method is not efficient;
it
is not almighty.
If we wish to apply this method to the geometric Bogomolov
conjecture, we have to restrict the setting where the method
works.
In what follows, we will apply the equidistribution argument
to the
proof of a weaker version of 
Theorem~\ref{thm:yamaki6}
(cf. Proposition~\ref{prop:main111}).
In that setting, a detailed analysis of the canonical measures gives us enough
information to lead to a contradiction in
(\ref{eq:equalityofmeasures-nonarch}).

\section{Structure of canonical measures} \label{sect:structuretheorem}

In this section,
we recall the structure of canonical measure due to
Gubler.
We refer to \cite{gubler4} for the detail.

Let $A$ be an abelian variety over $\overline{K}$
and let $X$ be a closed subvariety of $A$.
Further, fix a $v \in M_{\overline{K}}$ and put $\KK := \overline{K}_v$.
We write $X_v := X \otimes_{\overline{K}} \KK$.
Remark that $X_v$ can be defined over a CDV subfield (of $\KK$).

\subsection{Semistable alteration} \label{subsect:SA}

The canonical measures on $X^{\an}$ 
are described by Gubler
by using a semistable alteration.
Let $A_v$ and $X_v$ be as above.
By
a \emph{semistable alteration} for $X$,
we mean 
a pair $ ( \mathscr{X}' , f : X' \to X_v )$
consisting of a
strictly semistable scheme 
$\mathscr{X}'$
over 
$\KK$ 
that can be defined over a 
CDV subring
and
a 
generically finite surjective morphism $f : X' \to X$
where $X'$ is the generic fiber of $\mathscr{X}'$.
By de Jong's semistable alteration theorem
in \cite{dejong},
there always exists a semistable alteration for $X_v$.

To describe the canonical measures on $X^{\an}$,
actually,
one needs to use a semistable alteration 
$ ( \mathscr{X}' , f : X' \to X_v )$
with an additional
property.
Indeed,
Gubler gives a description of the canonical measure on $X^{\an}$
in terms of the skeleton of $\mathscr{X}'$,
which we are going to explain in this section,
under the condition that 
the semistable alteration
$ ( \mathscr{X}' , f : X' \to X_v )$
for $X$ is
``compatible with
some Mumford model of $A_v$''.

A Mumford model
of $A_v$
is proper admissible formal scheme
$\mathscr{A} \to \Spf ( \KK^{\circ} )$
with ``Raynaud generic fiber'' $A_v^{\an}$
that is constructed 
from a (rational) polytopal decomposition of a real torus via the
valuation map.
In the appendix, we give a quick review of admissible formal schemes,
their Raynaud generic fibers,
and the valuation maps for abelian varieties.
However,
we do not give the definition of Mumford models in this paper.
In fact,
the readers do not need to know the definition in the following arguments
if they assume Lemmas~\ref{lem:exist:MM} and \ref{lem:semistable_alteration} as black boxes.

\begin{Lemma}
\label{lem:exist:MM}
Let $A_1$ and $A_2$ be abelian varieties
over $\overline{K}$, and let
$\phi : (A_1)_v \to (A_2)_v$ be a homomorphism of abelian varieties.
Then there exist Mumford models $\mathscr{A}_1$ of $(A_1)_v$
and $\mathscr{A}_2$ of $(A_2)_v$ and a morphism
$\varphi : \mathscr{A}_1 \to \mathscr{A}_2$
whose restriction $(A_1)_v^{\an} \to (A_2)_v^{\an} $
to the Raynaud generic fibers coincides with the analytification
$\phi^{\an} : 
(A_1)_v^{\an} \to (A_2)_v^{\an}$ of $\phi$.
\end{Lemma}

A semistable alteration
$ ( \mathscr{X}' , f : X' \to X_v )$
for $X$ is said to be
\emph{compatible with
some Mumford model of $A_v$}
if there exists a Mumford model $\mathscr{A}$ of $A_v^{\an}$
such that the morphism 
$f : X' \to X_v \hookrightarrow 
A_{v}$ extends to a
morphism $\widehat{\mathscr{X}'} \to \mathscr{A}$
of admissible formal schemes.
The following lemma says that
for a given Mumford model,
there exists a semistable alteration that is compatible with
the Mumford model.

\begin{Lemma} \label{lem:semistable_alteration}
Let $\mathscr{A} \to \Spf ( \KK^{\circ} )$ be a Mumford model
of $A_v$.
Then there exists 
semistable alteration $( \mathscr{X}' , f : X' \to X_v )$
for $X_v$
such that the morphism $f$
extends to a morphism $\widehat{\mathscr{X}'} \to \mathscr{A}$
of admissible formal schemes.
\end{Lemma}

\begin{Remark} \label{rem:compatiblewithMM}
In the sequel,
we only consider
semistable alterations
that are compatible with
some Mumford model.
\end{Remark}

\subsection{Gubler's description of canonical measures} \label{subsect:structureofCM}

Let $L$ be an even ample line bundle on $A$.
Let $L_v$ denote the pullback of $L$ to $A_v$.
Put a canonical metric $|| \cdot ||$ 
(with respect to a rigidification of $L$)
on $L_v$,
and we have a canonically metrized line bundle $\overline{L_v} =
(L_v , || \cdot ||)$.
The canonical metric is a semipositive metric 
(cf. Proposition~\ref{prop:nonarchimedeancanoncalmeasure}),
and hence the restriction $\overline{L_v}|_{X_v}$
is also a semipositive metric
by Lemma~\ref{lem:functorialsemipositive}.

Let 
$( \mathscr{X}' , f : X' \to X_v )$
be a
semistable alteration 
(cf. Remark~\ref{rem:compatiblewithMM}).
By Lemma~\ref{lem:functorialsemipositive},
$f^{\ast} ( \overline{L_v} )$ is a semipositively metrized line bundle on 
$(X')^{\an}$,
so that
we consider the (normalized) Chambert-Loir measure
$\frac{1}{\deg_{ f^{\ast} ( L )} (X')}
\cherncl_1 ( f^{\ast} ( \overline{L_v} ))^{\wedge d}$
on $(X')^{\an}$, where $d := \dim ( X ) = \dim (X')$.
For a stratum $S$ of $\widetilde{\mathscr{X}'}$,
let $\Delta_S$ denote the canonical simplex corresponding to $S$.
Recall that the skeleton of $(X')^{\an}$ associated to $\mathscr{X}'$
is $S ( \mathscr{X}' ) = \bigcup_{S \in \str ( \widetilde{\mathscr{X}'})}
\Delta_S$.

\begin{Proposition} [Corollary~6.9 of \cite{gubler4}] \label{prop:measureonSA} 
With the above notation,
we can express 
\[
\frac{1}{\deg_{f^{\ast} (L)} (X')}
\cherncl_1 ( f^{\ast} ( \overline{L_v} ))^{\wedge d}
= 
\sum_{S \in \str ( \widetilde{\mathscr{X}'})} r_S \delta_{\Delta_S}
,
\]
where $r_S$ is a nonnegative real number,
$\Delta_S$ is the canonical simplex corresponding to 
$S$,
and $\delta_{\Delta_S}$ is 
the relative Lebesgue measure
on $\Delta_S$.
Furthermore,
whether
$r_S$ is positive or not
does not depend on the choice of the even ample line bundle $L$ on $A$.
\end{Proposition}

\begin{Remark} \label{rem:nondege}
\begin{enumerate}
\item
In \cite{gubler4},
Gubler defines the notion that $\Delta_S$ is \emph{non-degenerate
with respect to $f$}.
We  omit the Gubler's definition in this paper
but remark that
this notion is determined by $( \mathscr{X}' , f : X' \to X_v)$
and does not depend on the choice of $L$.
He proves that $\Delta_S$ is nondegenerate if and only if
$r_S > 0$ in the description of Proposition~\ref{prop:measureonSA}
for any even ample line bundle $L$.
Therefore, if we say
that $\Delta_S$ is non-degenerate with respect to $f$,
this means that $r_S > 0$ in Proposition~\ref{prop:measureonSA}
for one and hence any even ample line bundle $L$.
\item
Let $f^{\an} : (X')^{\an} \to X_v^{\an}$ be the 
morphism of analytic spaces associated to $f$.
If $\Delta_S$ is non-degenerate with respect to $f$,
then $f^{\an}$ restricts to an homeomorphism from $\Delta_S$ to 
$f^{\an} ( \Delta_S)$.
In fact, this property is a part of the definition of non-degeneracy due
to
Gubler.
\end{enumerate}
\end{Remark}

Let $\str_{f\mathrm{-nd}} ( \widetilde{\mathscr{X}'} )$
be the set of strata of $\widetilde{\mathscr{X}'} $
that are non-degenerate with respect to $f$.
It is know that the projection formula holds for Chambert-Loir measures
(cf. \cite[Proposition~3.8]{gubler4}),
so that
\addtocounter{Claim}{1}
\begin{align} \label{eq:projectionformula}
\mu_{X_v^{\an} , L} = 
f_{\ast} 
\left(
\frac{1}{\deg_{f^{\ast} (L)} (X')}
\cherncl_1 ( f^{\ast} ( \overline{L_v} ))^{\wedge d}
\right)
.
\end{align}
By 
Proposition~\ref{prop:measureonSA},
we obtain
\addtocounter{Claim}{1}
\begin{align} \label{eq:canonicalsubset}
\Supp ( \mu_{X_v^{\an} , L} ) = \bigcup_{S \in \str_{f\mathrm{-nd}} ( \widetilde{\mathscr{X}'} )}
f^{\an} ( \Delta_S )
.
\end{align}
Noting Remark~\ref{rem:nondege}~(1),
we see that $\Supp ( \mu_{X_v^{\an} , L} )$ does not depend on the choice
of $L$.
We set $S_{X_v^{\an}} := \Supp ( \mu_{X_v^{\an} , L} )$ and call it 
the \emph{canonical subset of $X_{v}^{\an}$}.

Gubler proved that $S_{X_v^{\an}}$ has a canonical
piecewise rational affine structure such that for any
semistable alteration 
$( \mathscr{X}' , f : X' \to X_v)$
(cf. Remark~\ref{rem:compatiblewithMM})
and for any stratum $S$ 
non-degenerate with respect to $f$,
the homeomorphism
$\Delta_S \to f^{\an} ( \Delta_S )$
(cf. Remark~\ref{rem:nondege}~(2))
given by
$f^{\an}$ is a rational piecewise linear map
(cf. \cite[Theorem~6.12]{gubler4}).
By  Proposition~\ref{prop:measureonSA} 
and
(\ref{eq:projectionformula}),
it follows that $\mu_{X_v^{\an} , L}$ is a linear combination
of relative Lebesgue measures on polytopes on $S_{X_v^\an}$
with this piecewise linear structure of $S_{X_v^\an}$.

\begin{Remark} \label{rem:imageofnondegenerate-poly}
If $\Delta_S$ is non-degenerate with respect to $f$,
then
$f^{\an} ( \Delta_S )$ is regarded as a polytopal set of pure dimension 
$\dim ( \Delta_S)$.
\end{Remark}

\section{Reduction to the nowhere degenerate case} \label{sect:reductiontoNDcase}

Under the preparation so far, we discuss
Theorem~\ref{thm:yamaki6}.
In this theorem,
the nontrivial part is to show that (b) implies (a),
that is:
let $A$ be an abelian variety over $\overline{K}$
and
let $\mathfrak{m}$ be the maximal nowhere degenerate abelian
subvariety of $A$;
if the geometric Bogomolov conjecture holds for $\mathfrak{m}$,
then it holds for $A$.
This assertion is obtained by showing 
the following theorem.

\begin{Theorem} \label{thm:yamaki6_7.21}
Let $A$ be an abelian variety over $\overline{K}$
and
let $\phi : A \to \mathfrak{m}$ be a surjective homomorphism,
where
$\mathfrak{m}$ is the maximal nowhere degenerate abelian subvariety of $A$.
Let $X$ be a closed subvariety of $A$.
Suppose that $X$ has dense small points and that $\phi (X)$ is a special 
subvariety.
Then $X$ is a special subvariety.
\end{Theorem}

Let us
explain 
how Theorem~\ref{thm:yamaki6} follows from Theorem~\ref{thm:yamaki6_7.21}.
As is noted before,
we only have to show that if the geometric Bogomolov conjecture holds
for $\mathfrak{m}$, then it holds for $A$.
First we note that
there exists a surjective homomorphism
$\phi : A \to \mathfrak{m}$.
Indeed,
by the Poincar\'e complete reducibility theorem,
there exists an abelian subvariety $G$ of $A$ such that $\mathfrak{m} + G = A$ 
and $\mathfrak{m} \cap G$ is finite.
Since
the natural homomorphism $\mathfrak{m} \times G \to A$
is an isogeny,
there exists an isogeny $A \to \mathfrak{m} \times G$
(cf. \cite[p. 157]{mumford}).
This homomorphism composed with the natural projection to $\mathfrak{m}$
gives a desired $\phi$.
Now, suppose that $X$ has dense small points.
Then $\phi (X)$ also has dense small points
by Lemma~\ref{lem:density-quotient}.
Assume that the geometric Bogomolov conjecture holds for $\mathfrak{m}$.
It follows that $\phi (X)$ is a special subvariety.
By Theorem~\ref{thm:yamaki6_7.21}, we conclude that $X$ is a special subvariety.

We do not give a complete proof of Theorem~\ref{thm:yamaki6_7.21}
in this paper.
Instead,
we give 
the idea of the proof of the following weaker version of the theorem.
Here we call $\dim ( \mathfrak{m} )$ the nowhere degeneracy rank
of $A$ and denote it by $\ndr (A)$.

\begin{Proposition} \label{prop:main111}
Let $A$ be an abelian variety over $\overline{K}$.
Assume that $\ndr (A) = 0$.
Let $X$ be a closed subvariety of $A$.
Suppose that $X$ has dense small points.
Then $X$ is a torsion subvariety.
\end{Proposition}

The above proposition generalizes
Gubler's theorem (cf. Theorem~\ref{thm:gubler:text}).
Let $A$ be an abelian over $\overline{K}$,
and assume that it is totally degenerate at some place
$v$.
Let $\mathfrak{m}$ be the maximal nowhere degenerate 
abelian subvariety of $A$.
Since $\mathfrak{m}_v$ is a non-degenerate abelian subvariety
of $A_v$ and since this is totally degenerate,
$\mathfrak{m}_v$ is trivial.
It follows
that $\mathfrak{m} = 0$, namely, $\ndr (A) = 0$.
Thus one can apply the proposition to $A$ to obtain the
conclusion.

We remark that we may have $\ndr (A) = 0$ 
even if $A$ is not
\emph{totally} degenerate at any place.
Indeed, there exists a simple abelian variety $A$ 
that is not totally degenerate at any place but is
degenerate
at some place;
since $A$ is simple and degenerate,
we have $\ndr (A) = 0$.

\subsection{Strict supports} \label{subsect:strictsupport}

In our proof of 
Proposition~\ref{prop:main111},
we apply the method of Zhang.
We hope to get a contradiction from the equality
(\ref{eq:equalityofmeasures-nonarch}).
To do that, we define the notion of strict supports of measures
and investigate the structure of the strict supports.

As we noted in \S~\ref{subsect:structureofCM},
the canonical subset $S_{X_v^{\an}}$ has a piecewise affine structure,
so that
we have a notion of polytopal decomposition of $S_{X_v^{\an}}$.
We say that a polytopal decomposition
$\Sigma$ 
of $S_{X_v^{\an}}$
is \emph{$f$-subdivisional}
if 
\[
f^{\an} ( \Delta_S ) \cap \Sigma
:=
\{
\sigma \in \Sigma
\mid f^{\an} ( \Delta_S ) \cap \sigma \neq \emptyset
\}
\]
is a polytopal decomposition of $f^{\an} ( \Delta_S )$.
Note that for any rational polytopal decomposition 
$\Sigma_0$ of $S_{X_v^{\an}}$,
there exists a rational subdivision $\Sigma$ of $\Sigma_0$
that is $f$-subdivisional.

\begin{Proposition} \label{prop:strofcm}
Let $L$ be an even ample line bundle on $A$.
Let $\Sigma$ be an
$f$-subdivisional
rational polytopal decomposition of 
$S_{X_v^{\an}}$.
Then
\[
\mu_{X_v^{\an} , L} =
\sum_{\sigma \in \Sigma} r'_{\sigma} \delta_{\sigma}
\]
for some non-negative real numbers $r'_{\sigma}$,
where $\delta_{\sigma}$ is the relative Lebesgue measure on $\sigma$.
\end{Proposition}

\Proof
This follows from Proposition~\ref{prop:measureonSA},
(\ref{eq:projectionformula}),
and the fact that $f^{\an} |_{\Delta_S} : \Delta_S \to f^\an ( \Delta_S)$
is a piecewise linear map for $\Delta_S$ non-degenerate with respect
to $f$.
\QED

\begin{Remark} \label{rem:forstrictsupport}
It follows from equality
(\ref{eq:projectionformula})
and Remark~\ref{rem:nondege}~(1)
that
 $r'_{\sigma} > 0$ if and only if there exists
a non-degenerate stratum $S$ of $\widetilde{\mathscr{X}'}$ such that
$\sigma \subset \Delta_S$ and
$\dim ( \sigma ) = \dim ( \Delta_S )$.
\end{Remark}

Let $\sigma$ be a polytope in the canonical subset $S_{X_v^{\an}}$
of $X_v^{\an}$.
Let $\mu_{X_v^{\an}}$ be the canonical measure on $X_v^{\an}$ associated
to an even ample line bundle on $A$.
We say that $\sigma$ is a 
\emph{strict support of $\mu_{X_v^{\an}}$}
if there exists an $\epsilon > 0$
such that
$\mu_{X_v^{\an}} - \epsilon \delta_{\sigma}$
is a positive measure,
where $\delta_{\sigma}$ is a relative Lebesgue measure on $\sigma$.
By the last assertion of Proposition~\ref{prop:measureonSA},
this notion 
does not depend on the
choice of the even ample line bundle on $A$
and well defined for $X$.

Let $(\mathscr{X}', f : X' \to X_v)$ be 
a semistable alteration 
(cf. Remark~\ref{rem:compatiblewithMM}).
Let $\Sigma$ be an 
$f$-subdivisional
polytopal decomposition of $S_{X_v^{\an}}$.

By Proposition~\ref{prop:strofcm},
we write
\[
\mu_{X_v^{\an}} =
\sum_{\sigma \in \Sigma}
r'_{\sigma} \delta_{\sigma}
\]
with $r'_\sigma \geq 0$.
Then $\sigma \in \Sigma$
is a strict support of $\mu_{X_v^{\an}}$
if and only if $r'_{\sigma} > 0$.

We would like to see how $f^{\an}|_{\Delta_S} : \Delta_S \to X^\an$
behaves.
If $\Delta_S$ is nondegenerate with respect to $f$,
then
the map $f^{\an} |_{\Delta_S} : \Delta_S \to S_{X_v^{\an}}$
is an isomorphism onto its image,
which
follows from the definition of the piecewise affine
structure for $S_{X_v^{\an}}$.
On the other hand,
if $\Delta_S$ is degenerate
(i.e. not non-degenerate),
then the behavior of $f^{\an} |_{\Delta_S}$
may be complicated.
Nevertheless, 
when $f^{\an} ( \Delta_S ) \subset S_{X_v^{\an}}$,
which does not necessarily means that $\Delta_S$ is nondegenerate
with respect to $f$,
we can describe 
$f^{\an}|_{\Delta_S} : \Delta_S \to X^\an$
in the following sense.

\begin{Lemma} \label{lem:behaviorofdegenerate}
Take an 
$S \in \str ( \widetilde{\mathscr{X}'})$.
Suppose that $f^{\an} ( \Delta_S ) \subset S_{X_v^{\an}}$.
Then we have the following.
\begin{enumerate}
\item
The map $f^{\an} |_{\Delta_S} : \Delta_S \to S_{X_v^{\an}}$
is a piecewise linear map.
Furthermore, there exists a polytope $P \subset \Delta_S$
with $\dim ( P )  = \dim ( \Delta_S  )$
such that $f^{\an} (P  )$ is a polytope and
$f^{\an} |_{P } : P \to f^{\an} (P )$
is an affine map for some rational polytopal decomposition
of $S_{X_v^{\an}}$.
\item
There exists a non-negative integer $r$ 
with the following property.
Let $P$ be a polytope in $\Delta_S$ such that
$\dim (P) = \dim ( \Delta_S )$
and such that
$f^{\an} |_{P} : P \to S_{X_v^{\an}}$
is an affine map for some rational polytopal decomposition
of $S_{X_v^{\an}}$.
Then $\dim ( f^{\an} ( P ) )= r$.
\end{enumerate}
\end{Lemma}

\Proof
It is known that there exists a continuous
map $\overline{\val} : A_v^{\an} \to \RR^{n} / \Lambda$
where $\Lambda$ is a complete lattice of $\RR^{n}$ 
with $\Lambda \subset \QQ^{n}$;
see \S~\ref{subsection:ap2}.
Remark that $\RR^{n} / \Lambda$ has a natural
piecewise affine structure induced from $\RR^{n}$.
The map $\overline{\val}$ has the following properties:
\begin{enumerate}
\renewcommand{\labelenumi}{(\alph{enumi})}
\item
The restriction $\overline{\val} |_{S_{X_v^{\an}}} : S_{X_v^{\an}} \to
\RR^n$ is a finite piecewise affine map;
\item
For any $S \in \str ( \widetilde{\mathscr{X}'})$,
the map $\overline{\val} \circ f^{\an} |_{\Delta_S} : \Delta_S 
\to \RR^{n} / \Lambda$
is an affine map.
\end{enumerate}
We refer to \cite[\S~4]{gubler4}.

By (a) above,
$f^{\an} ( \Delta_S )$ is a finite union of polytopes of dimension $r$
if and only if so is $\overline{\val} ( f^{\an} ( \Delta_S ))$,
and in fact, $f^{\an} ( \Delta_S )$ is a polytope by (b).
Thus we have the first assertion of (1).
The second assertion is obvious.

To show (2), put $r:= \dim ( ( \overline{\val} \circ f^{\an}) ( \Delta_S ) )$
by (b).
Then by (a), 
$\dim ( f^{\an} ( P ) )= r$.
This shows the assertion.
\QED

Suppose that
$f^{\an} ( \Delta_S ) \subset S_{X_v^{\an}}$.
Let $r$ be the integer as in the above lemma.
Then
$f^{\an} ( \Delta_S )$ is a finite union of
$r$-dimensional polytopes.
Thus $\dim ( f^{\an} ( \Delta_S ) )$ makes sense,
i.e., $\dim ( f^{\an} ( \Delta_S ) ) = r$.

The following is an immediate consequence of
Gubler's description of canonical measures.

\begin{Lemma} \label{lem:existenceofnondegstr}
Let $(\mathscr{X}', f : X' \to X_v)$ be 
a semistable alteration 
(cf. Remark~\ref{rem:compatiblewithMM}).
Let $\Sigma$ be an 
$f$-subdivisional
polytopal decomposition of $S_{X_v^{\an}}$.
Let $\sigma \in \Sigma$ be a strict support of 
$\mu_{X_v^{\an}}$.
Then there exists a stratum $S \in \str ( \widetilde{\mathscr{X}'})$
with the following properties:
\begin{enumerate}
\renewcommand{\labelenumi}{(\alph{enumi})}
\item
$\sigma \subset f^{\an} ( \Delta_S)$;
\item
$f^{\an} ( \Delta_S ) \subset S_{X_v^{\an}}$;
\item
$
\dim ( \sigma) = \dim ( \Delta_S)$.
\end{enumerate}
\end{Lemma}

\Proof
By Remark~\ref{rem:forstrictsupport},
there exists a non-degenerate stratum $S \in \widetilde{\mathscr{X}'}$
such that 
$\sigma \subset f^{\an} ( \Delta_S)$ and $\dim ( \sigma) = \dim ( \Delta_S)$.
By 
(\ref{eq:canonicalsubset}),
we also have $f^{\an} ( \Delta_S ) \subset S_{X_v^{\an}}$.
\QED

Remark that if $\sigma$ and $\Delta_S$ satisfy the three conditions
in Lemma~\ref{lem:existenceofnondegstr},
then $\dim (\sigma) = \dim ( f^{\an} ( \Delta_S))$.
Take $\sigma \in \Sigma$.
We say tentatively that \emph{$\Delta_S$ is over $\sigma$}
if $\Delta_S$ satisfies conditions (a) and (b) in 
Lemma~\ref{lem:existenceofnondegstr} as well as 
$\dim (\sigma) = \dim ( f^{\an} ( \Delta_S))$.
With this word, 
Lemma~\ref{lem:existenceofnondegstr}
means that if $\sigma$ is a strict support,
then there exists $\Delta_S$ over $\sigma$ such that
$
\dim ( \sigma) = \dim ( \Delta_S)$, which is condition (c).

Conversely,
the following proposition shows that
if $\sigma $ is a strict support, then
any canonical simplex $\Delta_S$ over $\sigma$
should satisfy $
\dim ( \sigma) = \dim ( \Delta_S)$.

\begin{Proposition} \label{prop:anyhasgooddimension}
Let $(\mathscr{X}', f : X' \to X_v)$ be 
a semistable alteration (cf. Remark~\ref{rem:compatiblewithMM}).
Let 
$\Sigma$ be a 
$f$-subdivisional
polytopal decomposition of $S_{X_v^{\an}}$.
Let $\sigma \in \Sigma$ be a strict support of 
the canonical measure $\mu_{X_v^{\an}}$
of an even ample line bundle on $A$.
Let $S$ be a stratum of $\widetilde{\mathscr{X}'}$
with the following properties:
\begin{enumerate}
\renewcommand{\labelenumi}{(\alph{enumi})}
\item
$\sigma \subset f^{\an} ( \Delta_S)$;
\item
$f^{\an} ( \Delta_S ) \subset S_{X_v^{\an}}$;
\item
$\dim ( \sigma) =
\dim f^{\an} ( \Delta_S ) $.
\end{enumerate}
Then $
\dim ( \sigma)
= \dim ( \Delta_S)$.
\end{Proposition}

The above proposition 
is essentially given in \cite[Proposition~5.12, Lemma~5.13]{yamaki6}.
In the proof,
we need detailed analysis on the non-degenerate strata
by using Mumford models, subdivision of polytopes,
toric method, etc.
We omit the proof and refer to \cite{yamaki6}.

\subsection{Proof of Proposition~\ref{prop:main111}}

In this subsection, we give an outline of the
proof of  Proposition~\ref{prop:main111}.
Let $A$ be an abelian variety over $\overline{K}$,
where $K$ is a function field.
Let $\mathfrak{m}$ be the maximal nowhere degenerate
abelian subvariety of $A$.
Set $\ndr (A) := \dim ( \mathfrak{m})$,
called
the nowhere degeneracy rank of $A$.

We recall basic properties of the nowhere degeneracy rank.

\begin{Lemma} [cf. Corollary~7.12 in \cite{yamaki6}] \label{lem:ndr}
Let $\phi : A \to A'$
be a surjective homomorphism of abelian varieties.
Then $\ndr (A) \geq \ndr (A')$.
\end{Lemma}

In Proposition~\ref{prop:main111},
we assume that $\ndr (A) = 0$.
The following lemma shows that under this assumption,
the support of the canonical measure has positive dimension.
It is essentially \cite[Proposition~7.16]{yamaki6}.

\begin{Lemma} \label{lem:canoicalsubsethaspositivedimension}
Let $A$ be an abelian variety over $\overline{K}$.
Assume that $\ndr (A) = 0$.
Let $X$ be a closed subvariety of $A$.
Suppose that $\dim (X) > 0$
and that $X$ has trivial stabilizer.
Then 
there exists a place $v$ of $\overline{K}$
such that $A$ is degenerate at $v$ and
such that
$\dim ( S_{X_v^{\an}}) > 0$,
namely, $S_{X_v^{\an}}$ contains a positive dimensional polytope.
\end{Lemma}

\Proof
Recall that for any $v \in M_{\overline{K}}$,
the valuation map $\overline{\val} : A_v \to \RR^{n} / \Lambda$
exists (cf. \S~\ref{subsection:ap2}).
Recall also the fact
that this restricts to a finite surjective map
$\overline{\val}|_{S_{X_v^{\an}}} : S_{X_v^{\an}} \to \overline{\val (X_v)}$,
which is a part of \cite[Theorem~1.1]{gubler4}.
By the implication from (b) to (a) in 
\cite[Proposition~7.16]{yamaki6},
we see that if $\dim (X) > 0$,
then $\dim ( \overline{\val (X_v)} ) > 0$ 
for some $v \in M_{\overline{K}}$.
Thus
$\dim ( S_{X_v^{\an}}) > 0$ for some $v \in M_{\overline{K}}$.
\QED

Let us prove Proposition~\ref{prop:main111}.
We argue by contradiction:
suppose that there exist an abelian variety over $\overline{K}$
of non-degeneracy rank $0$
and a closed subvariety  that is not a torsion subvariety
but has dense small points.
We take the quotient of the abelian variety by the stabilizer
of the closed subvariety.
Taking into account Lemma~\ref{lem:ndr},
we then construct an abelian $A$ 
with $\ndr (A) = 0$
and a 
closed subvariety $X$
that is non-torsion, has trivial stabilizer, has dense small points,
and has positive dimension.
For some large natural number $N$,
the difference morphism
$\alpha_N : X^{N} \to A^{N-1}$
given $\alpha_N ( x_1 , \ldots , x_N ) =
(x_1 - x_2 , \ldots , x_{N-1} - x_N )$
is generically finite (cf. Lemma~\ref{lem:genericallyfiniteN}).
Put $Z := X^N$ 
and $Y:= \alpha_N ( Z )$,
and
let $\alpha : Z \to Y$ be the morphism induced by $\alpha_N$.
Note that $\alpha$
is a generically finite surjective morphism.

Let $L$ be an even ample line bundle on $A$.
Then $L^{\boxtimes N}$ is an even ample line bundle on $A^{N}$.
Let $v$ be a place at which $A$ is degenerate.
By Lemma~\ref{lem:canoicalsubsethaspositivedimension},
we may take $v$ in such a way that $S_{X_v^{\an}}$ has positive dimension.
Let $\mu_{X_v^{\an}}$ be the canonical measure
on $X_v^{\an}$ of $L$.
Let $\mu_{Z_v^{\an}}$ 
be the canonical measure 
on $Z_v^{\an}$
of $L^{\boxtimes N}$
and let
$\mu_{Y_v^{\an}}$ be
the
canonical measures on $Y_v^{\an}$
of an even ample line bundle on $A^{N-1}$.
Then the argument using the equidistribution theorem
as in the proof of Zhang and the proof of Gubler
gives us
\[
\alpha^{\an}_{\ast} ( \mu_{Z_v^{\an}} ) = \mu_{Y_v^{\an}}
.
\]
Note in particular that $\alpha^{\an} (S_{Z_v^{\an}}) = 
S_{Y_v^{\an}}$.

By Lemma~\ref{lem:exist:MM},
we take Mumford models of $(A_v^N)^{\an}$ and $(A_v^{N-1})^{\an}$
respectively such that the morphism $\alpha^{\an}$ extends
to a morphism $\psi$ between the Mumford models.
We take a semistable alteration
$( \mathscr{Z}' , h : Z' \to Z_v )$
compatible with this Mumford model of $(A_v^N)^{\an}$
(cf. Lemma~\ref{lem:semistable_alteration}).
Since $\alpha : Z \to Y$ is a generically finite 
proper surjective
morphism, setting $g : = \alpha \circ h$,
we see that
$( \mathscr{Z}' , g : Z' \to Y_v )$
is also a semistable alteration compatible with a Mumford model.

Note that the restricted morphism 
$\alpha^{\an} : S_{Z_v^{\an}} \to S_{Y_v^{\an}}$ is a piecewise linear map.
Indeed,
let $\Delta_S$ be any non-degenerate stratum with respect to $h$.
It suffices to show that 
$\alpha^{\an}$ is piecewise linear on $h^{\an} ( \Delta_S )$.
By the definition of the piecewise linear structure on $S_{Z_v^{\an}}$,
this is equivalent to
$\alpha^{\an} \circ h^{\an} |_{\Delta_S} : \Delta_S \to
S_{Y_v^{\an}}$ being piecewise linear.
Since 
$( \mathscr{Z}' , g : Z' \to Y_v )$
is a semistable alteration compatible with a Mumford model
and $g^{\an} ( \Delta_{S}) \subset \alpha^{\an} (S_{Z_v^{\an}} ) =
S_{Y_v^{\an}}
$,
this follows from Lemma~\ref{lem:behaviorofdegenerate}~(1).

Let $|X_{v}^{\an}|^{N}$ denote the 
direct product of $N$-copies of $X_v^{\an}$ in the category of topological
spaces.
Since $Z_{v}^{\an}$ is the direct product of $N$-copies of $X_v^{\an}$
in the category of Berkovich spaces,
we have a natural continuous map $\beta : |Z_v^{\an}| \to |X_{v}^{\an}|^{N}$.
By \cite[Proposition~4.5]{yamaki5},
we have $\beta_{\ast} ( \mu_{Z_{v}^{\an}}) = \mu_{X_v^{\an}}^{N}$,
the product measure of $N$-copies of $\mu_{X_v^{\an}}$.
Further, we have
$\beta (S_{Z_v^{\an}}) = S_{X_v^{\an}}^{N}$.

We show that there exists a strict support 
$\tau$
of $\mu_{Z_v^{\an}}$
such that 
\addtocounter{Claim}{1}
\begin{align}
\label{align:main:proof:1}
\dim ( \alpha^{\an} (\tau) ) < \dim (\tau).
\end{align}
Here, this is not precise, but
we argue as if the above $\beta$
were an isomorphism,
that is,
$S_{Z_v^{\an}} = S_{X_v^{\an}}^{N}$
and
$\mu_{Z_{v}^{\an}} = \mu_{X_v^{\an}}^{N}$,
to explain the idea;
see \S~\ref{subsect:morepreciseargument}
below
for more precise argument.
Since $\dim ( S_{X_v^{\an}} ) > 0$,
there exists a strict support $\sigma$
of $\mu_{X_v^{\an}} $ with $\dim ( \sigma ) > 0$.
Then
$\sigma^{N}$ is a strict support of $\mu_{X_v^{\an}}^{N}$
with $\dim ( \sigma^{N} ) > 0$.
Since $\alpha$ contracts the diagonal of $Z = X^{N}$ to a point,
the diagonal of $\sigma^{N}$ contracts to a point.
Furthermore, since $\alpha^{\an} |_{\sigma^{N}}$ is piecewise linear,
we can therefore take a
polytope $\tau$
such that $\tau \subset \sigma^{N}$,
$\dim ( \tau ) = \dim ( \sigma^{N} )$,
$\alpha^{\an} |_{\tau}$ is linear,
and such that
$\alpha^{\an} ( \tau )$ is a polytope with
$\dim ( \alpha^{\an} ( \tau )) < \dim ( \tau )$.
By the choice of $\tau$,
this is a strict support of $\mu_{Z_v^{\an}} $.

Since $\tau$ is a strict support of $\mu_{Z_v^{\an}}$,
Lemma~\ref{lem:existenceofnondegstr}
give us 
an $S \in \str (\widetilde{\mathscr{Z}'})$
such that 
$
\tau \subset h^{\an} ( \Delta_S)
$,
$h^{\an} ( \Delta_S) \subset S_{Z_v^{\an}}$
and such that
\addtocounter{Claim}{1}
\begin{align}
\label{align:main:proof:2}
\dim ( \tau )
= \dim ( \Delta_S).
\end{align}
Since $\dim (\tau) \leq \dim ( h^{\an} ( \Delta_S) )
\leq \dim ( \Delta_S ) = \dim ( \tau)$,
we remark that
\addtocounter{Claim}{1} 
\begin{align}
\label{align:main:proof:3}
\dim ( \tau ) = \dim ( h^{\an} ( \Delta_S) ).
\end{align}

Recall that $( \mathscr{Z}' , g : Z' \to Y_v )$
is also a semistable alteration for $Y_v$
compatible with a Mumford model.
We would like to apply Proposition~\ref{prop:anyhasgooddimension}
in place of $\sigma$ and $f$ with $\alpha^{\an}(\tau)$ and $g$,
respectively,
so that we need to check the conditions in this proposition.
Since 
$\tau$ is a strict support of $\mu_{Z_v^{\an}} $
and
$\alpha^{\an}_{\ast} ( \mu_{Z_v^{\an}} ) = \mu_{Y_v^{\an}}$,
$\alpha^{\an} ( \tau )$ is a strict support of $\mu_{Y_v^{\an}}$.
Since $\tau \subset h^{\an} ( \Delta_S)$,
we have
\[
\alpha^{\an} ( \tau ) \subset \alpha^{\an} (h^{\an} ( \Delta_S)) =
g^{\an} ( \Delta_S),
\]
which is condition (a) in Proposition~\ref{prop:anyhasgooddimension}.
Since 
$g^{\an} (  \Delta_S ) =
\alpha^{\an} ( h^{\an} ( \Delta_S ))
\subset \alpha^{\an} ( S_{Z_v^{\an}} ) = S_{Y_v^{\an}}$,
condition (b) in Proposition~\ref{prop:anyhasgooddimension} is satisfied.
Further,
we take 
by Lemma~\ref{lem:behaviorofdegenerate}~(1)
a polytope $P \subset \Delta_S$ such that 
$\dim (P) = \dim ( \Delta_S )$,
$h^{\an} (P) \subset \tau$,
and such that $h^{\an} |_{P}$ is an affine map.
By Lemma~\ref{lem:behaviorofdegenerate}~(2),
we note $\dim ( h^{\an} (P)) = \dim (h^{\an} (\Delta_S) )$,
which equals $\dim ( \tau )$ by (\ref{align:main:proof:3}).
Since $\alpha^{\an}|_{\tau}$ is affine,
it follows that
\[
\dim ( g^{\an} ( P )) = \dim ( \alpha^{\an} ( h^{\an} ( P )))
= 
\dim (\alpha^{\an} ( \tau))
.
\]
By Lemma~\ref{lem:behaviorofdegenerate}~(2),
this means that $\dim (\alpha^{\an} ( \tau)) =  \dim ( g^{\an} ( \Delta_S ))$,
which is
condition (c) in Proposition~\ref{prop:anyhasgooddimension}.
Thus applying
Proposition~\ref{prop:anyhasgooddimension},
we obtain
$\dim ( \alpha^{\an} ( \tau) ) = \dim ( \Delta_S)$.

On the other hand, it follows from
(\ref{align:main:proof:1})
and
(\ref{align:main:proof:2})
that
$\dim ( \alpha^{\an} ( \tau )) < \dim ( \Delta_S)$,
which contradicts what we have show above.
Thus we complete the proof of Proposition~\ref{prop:main111}.

\subsection{Complement} \label{subsect:morepreciseargument}
In the above proof, 
we took a strict support $\tau$ of $\mu_{Z_v^{\an}}$
such that $\dim ( \alpha^{\an} ( \tau )) < \dim ( \tau )$.
In that argument,
we pretended that we had $\mu_{Z_v^{\an}} = 
\mu_{X_v^{\an}}^{N}$, which was not precise.
Here, we 
give a precise argument to obtain $\tau$ above.
We use the valuation map $\overline{\val} : A_v^{\an} \to \RR^{n} / \Lambda$,
where $n$ is the dimension of the torus part of $A_v^{\an}$
(cf. \S~\ref{subsection:ap2}).
Since the valuation map is compatible with direct product,
the map $\overline{\val}^{N} : (A_{v}^{N})^{\an} \to  (\RR^{n} / \Lambda)^{N}$
and $\overline{\val}^{N-1} : (A_{v}^{N-1})^{\an} \to  (\RR^{n} / \Lambda)^{N-1}$
are the valuation maps.
By \cite[Lemma~4.1 and Proposition~4.5]{yamaki5},
we have $(\overline{\val}^{N})_{\ast} ( \mu_{Z_v^{\an}}) = 
\overline{\val}_{\ast} (\mu_{X_v^{\an}})^{N}$.

Let $\sigma$ be a 
positive dimensional strict support of $\mu_{X_v^{\an}}$.
Then $\overline{\val} (\sigma)$ is a strict support of 
$\overline{\val}_{\ast} (\mu_{X_v^{\an}})$ of positive dimension.
It follows that
$(\overline{\val} (\sigma))^{N}$ is a strict support of 
$(\overline{\val}^{N})_{\ast} ( \mu_{Z_v^{\an}})$ positive dimension.
On the other hand, since the difference map induces the difference
map from $(\RR^{n} / \Lambda)^{N}$ to $(\RR^{n} / \Lambda)^{N-1}$,
we see that $\overline{\alpha}_{\aff}$
contracts the diagonal of $(\overline{\val} (\sigma))^{N}$ to a point.
Since $\overline{\alpha}_{\aff}$ is piecewise linear,
it follows that
there exists a polytope $\tau' \subset (\overline{\val} (\sigma))^{N}$
such that 
$\dim ( \tau' ) = \dim ((\overline{\val} (\sigma))^{N})$
and
$\dim ( \overline{\alpha}_{\aff} ( \tau' )) < \dim ( \tau')$.
Since 
$(\overline{\val} (\sigma))^{N}$ is a strict support of 
$(\overline{\val}^{N})_{\ast} ( \mu_{Z_v^{\an}})$,
so is $\tau'$.
Since the restriction $S_{Z_v^{\an}} \to 
\overline{\val}^{N} (S_{Z_v^{\an}})$
of $\overline{\val}^{N}$ is a finite surjective piecewise linear map,
there exists a strict support $\tau$ of $\mu_{Z_v^{\an}}$ such that
$\overline{\val}^{N} ( \tau ) \subset \tau'$
and
$\dim (\tau) = \dim (\tau')$.
Noting that $\overline{\val}^{N-1} \circ \alpha^{\an} = 
\overline{\alpha}_{\aff} \circ \overline{\val}^{N}$,
we see that
\[
\overline{\val}^{N-1} \circ \alpha^{\an} ( \tau ) = 
\overline{\alpha}_{\aff} ( \overline{\val}^N ( \tau ))
\subset \overline{\alpha}_{\aff} ( \tau' )
<
\dim (\tau') = \dim (\tau).
\]
Again since $\overline{\val}^{N-1}|_{S_{Y_v^{\an}}}$ is a finite 
piecewise linear map,
that shows $\dim ( \alpha^{\an} ( \tau )) < \dim (\tau)$.
Thus we have a required $\tau$.

\begin{Remark}
(This is a remark for the proof of \cite[Theorem~6.2]{yamaki6}.)
In the proof
of Proposition~\ref{prop:main111}, we take an even ample line bundle $L$ on $A$
and consider the even ample line bundle $L^{\boxtimes N}$
on $A^N$.
However, we can argue with any even ample line bundle on $A^{N}$.
(In the proof of \cite[Theorem~6.2]{yamaki6},
we actually did so.)
Indeed,
in the argument
we assume that $Z$ has dense small points,
so that by the equidistribution theorem
(Theorem~\ref{prop:ned}),
the canonical measure does not depend on the choice of
the even ample line bundle on $A^N$. 
(But in fact, we should also begin with $L^{\boxtimes N}$
in the proof of \cite[Theorem~6.2]{yamaki6} 
for the sake of
logical simplicity.)
\end{Remark}

\section{Canonical height of closed subvarieties} \label{sect:heightofsubvarietiels}

\subsection{Height of a closed subvariety}

Let $A$ be an abelian variety over $\overline{K}$ and let
$L$ be an even ample line bundle on $A$.
Let $Z$ be a cycle on $A$
of dimension $d$.
Then
the canonical height of $Z$
with respect to $L$
is defined.
We briefly explain what it is when 
$K$ is the function field
 of a variety
$\mathfrak{B}$ with $\dim ( \mathfrak{B}) = 1$.
The construction is similar to the construction
of the canonical metric;
we construct a global version of
a sequence of models as in (\ref{eq:sequenceofmodelsforcanonial}).
Fix an integer $m$ with $m > 1$
and fix a rigidification of $L$.
Let $K'$ be a finite extension of $K$
over which $Z$ can be defined.
Let $\mathfrak{B}'$ be the normalization of $\mathfrak{B}$
in $K'$.
We begin with any proper flat model $(\pi_1 : \mathcal{A}_1 
\to \mathfrak{B}' , \mathcal{L}_1)$
such that
$\mathcal{L}_1$ is nef.
There exist such $\pi_1$ and $\mathcal{L}_1$ because $L$ is ample.
Next we take a proper flat model
$\pi_2 : \mathcal{A}_2 
\to \mathfrak{B}'$
of $A$ such that the morphism $[m] : A \to A$
extends to $f_1 : \mathcal{A}_2 \to \mathcal{A}_1$.
Set $\mathcal{L}_2 : = f_1^{\ast} ( \mathcal{L}_1 )$.
Then the pair
$( \mathcal{A}_2 , \mathcal{L}_2 )$
is a model of $(A , L^{\otimes m^{2}})$ 
via the isomorphism $[m]^{\ast} ( L ) \cong L^{\otimes m^{2}}$
that respects the rigidification.
Again, we take a model $\pi_3 : \mathcal{A}_3 \to \mathfrak{B}'$
of $A$ such that the morphism $[m] : A \to A$
extends to $f_2 : \mathcal{A}_3 \to \mathcal{A}_2$,
and set $\mathcal{L}_3 : = f_2^{\ast} ( \mathcal{L}_2 )$.
Then the pair
$( \mathcal{A}_3 , \mathcal{L}_3 )$
is a model of $(A . L^{\otimes m^{4}})$.
Repeating this process, we obtain a sequence
$
\left(
\pi_n : \mathcal{A}_n \to \mathfrak{B}' , \mathcal{L}_n
\right)_{n \in \NN}
$
whose $n$th term is a model of $(A , L^{\otimes m^{2 (n-1)}})$
with $\mathcal{L}_n$ 
 nef on $\mathcal{A}_n$.

Let $\mathcal{Z}_n$ be the closure of $Z$ in $\mathcal{A}$.
Note that $\mathcal{Z}_n$ is a model of $Z$.
We consider the sequence 
\addtocounter{Claim}{1}
\begin{align} \label{eq:limitcanonicalheightofvarieties}
\left(
\frac{\deg ( \cherncl_1 ( \mathcal{L}_n)^{\cdot (d+1)} 
\cdot [ \mathcal{Z}_n])}{m^{2(n-1)(d+1)} [K' : K]}
\right)_{n \in \NN}
\end{align}
of rational numbers.
Then one shows that this sequence converges to a real number
and that this limit depends only on $(Z , L)$.
This number is the \emph{canonical height} of $Z$ with respect to $L$,
denoted by $\widehat{h}_L (Z)$.
When $X$ is a closed subvariety of $A$, 
we regard $X$ as a cycle and define the canonical height $\widehat{h}_L (X)$
to be the canonical height of this cycle.

\begin{Remark} \label{rem:nonnegativecanonicalheight}
In the above construction,
each
model $\mathcal{L}_n$ 
is nef on $\mathcal{A}_n$,
and hence
the sequence in (\ref{eq:limitcanonicalheightofvarieties})
is nonnegative.
It follows that
$\widehat{h}_L (X) \geq 0$.
This inequality holds for $\mathfrak{B}$ of
arbitrary dimension,
in fact.
\end{Remark}

One finds 
 in the
following proposition
an important property of the canonical height
in view of the density of small points.

\begin{Proposition} [Corollary~4.4 of \cite{gubler2}] \label{prop:dense-height0}
Let $X$ be a closed subvariety of $A$.
Then it has dense small points if and only if
$\widehat{h}_L (X) = 0$.
\end{Proposition}

\begin{Remark} \label{rem:heightofablianvariety}
Since $A$ is an abelian variety over $\overline{K}$,
$A$ has dense small points,
and hence
$\widehat{h}_L (A) = 0$
by the above proposition.
\end{Remark}

We explain the idea of the proof of Proposition~\ref{prop:dense-height0},
which is essentially due to Zhang.
We consider the essential minima
\[
e_{1} (X , L):=
\sup_{Y} \inf_{s \in X \left( \overline{K} \right) \setminus Y } \widehat{h}_L (x),
\]
where $Y$ runs through 
all closed subsets of codimension $1$ in $X$.
Then in general,
one shows that
\[
\deg_{L} (X) e_1 (X , L)
\leq 
\widehat{h}_L (X)
\leq (\dim (X) + 1) \deg_{L} (X) e_1 (X , L)
.
\]
By the definition of $e_1 (X , L)$,
$X$ has dense small points
if and only if
$e_1 (X , L) = 0$.
It follows from the above inequality
that $X$ has dense small points if and only if
$\widehat{h}_L (X) = 0$.

\subsection{Zhang's admissible pairing and the canonical height} \label{subsect:zcnon}

We mentioned in \S~\ref{subsubsect:adimsiiblepairing}
and \S~\ref{subsubsect:cinkir}
that
some partial answers to the geometric Bogomolov conjecture for curves
have been obtained before,
and that those results are obtained by showing the positivity
of the admissible pairing $( \omega_a , \omega_a)_a$
of the admissible dualizing sheaf $\omega_a$ of $C$.
This argument using the positivity of the admissible pairing
concerns Proposition~\ref{prop:dense-height0}.
Here, we explain why $( \omega_a , \omega_a)_a > 0$
implies the conjecture for the curve.

Let us recall the setting of the geometric Bogomolov conjecture
for curves.
Let $C$ be a smooth projective curve of genus $g \geq 2$
over $\overline{K}$
and
let $J_C$ be the Jacobian variety of $C$.
Fix a divisor $D$ on $C$ of degree $1$
and let
$
\jmath_{D} : C \hookrightarrow 
J_C
$
be the embedding
defined
by $\jmath_D (x) := x - D$.
Put $X := \jmath_D(C)$.
Assume that $C$ cannot be defined over $k$.
What we should show is that $X$ does not have dense small points.

Recall that $h_{NT}$ equals the canonical height 
$\widehat{h}_L$ associated
to $L = \OO_{J_C} ( \theta )$,
where $\theta$ is a symmetric theta divisor
(cf. Remark~\ref{rem:NTheight}).
In \cite{zhang1},
it is shown that
$\widehat{h}_L ( X ) \geq \alpha ( \omega_a , \omega_a)_a$
for some $\alpha > 0$,
and the equality holds if $(2 g - 2) D_1$ is a canonical divisor
on $C$.
It follows that
if one shows $( \omega_a , \omega_a)_a > 0$,
then $\widehat{h}_L ( X ) > 0$,
so that 
by Proposition~\ref{prop:dense-height0},
$X$ does not have dense small points.
That is the reason why 
$( \omega_a , \omega_a)_a > 0$ suffices for the Bogomolov
conjecture for curves.

\section{Proof of Theorem~\ref{thm:yamaki7}} \label{sect:proofofGBCforcurves}

In this section, we 
describe the idea of the proof of Theorem~\ref{thm:yamaki7}.

\subsection{Model for a nowhere degenerate abelian variety} \label{subsect:modelforNDAV}

For a nowhere degenerate abelian variety, we have
a good model as follows.

\begin{Proposition} [Proposition~2.5 of \cite{yamaki7}] \label{prop:model11}
Let $A$ be a nowhere degenerate abelian variety over $\overline{K}$
and let $L$ be a line bundle on $A$.
Then there exist a finite extension $K'$ of $K$,
a proper morphism $\pi : \mathcal{A} \to \mathfrak{B}'$,
where $\mathfrak{B}'$ is the 
normalization
of $\mathfrak{B}$ in $K'$,
and 
a line bundle $\mathcal{L}$ on $\mathcal{A}$
satisfying the following conditions.
\begin{enumerate}
\renewcommand{\labelenumi}{(\alph{enumi})}
\item
The pair $\left( \pi , \mathcal{L} \right)$
is a model of $(A ,L)$.
\item
There exists an open subset $\mathfrak{U} \subset \mathfrak{B}'$
with $\codim ( \mathfrak{B}' \setminus \mathfrak{U}
,  \mathfrak{B}' ) \geq 2$
such that
the restriction $\pi' : 
\pi^{-1} (\mathfrak{U}) \to \mathfrak{U}$ 
of $\pi$ is 
an abelian scheme.
\item
Let $0_{\pi'}$ be the zero-section of the abelian
scheme $\pi'$ in (b).
Then $0_{\pi'}^{\ast} ( \mathcal{L} ) \cong \OO_{\mathfrak{U}}$.
\end{enumerate}
\end{Proposition}

\begin{Remark} \label{rem:heightoverdimensin2}
Consider the case where $\dim ( \mathfrak{B} ) = 1$.
Then Proposition~\ref{prop:model11} claims the following.
Let $A$ and $L$ be as in this proposition.
Then, there exist a finite extension $K'$ of $K$
and a model
$( \pi : \mathcal{A} \to \mathfrak{B}',
\mathcal{L} )$
of $(A,L)$,
where $\mathfrak{B}'$ is the 
normalization
of $\mathfrak{B}$ in $K'$,
satisfying the following conditions.
\begin{enumerate}
\renewcommand{\labelenumi}{(\alph{enumi})}
\item
The morphism
$\pi$ is an abelian scheme.
\item
Let $0_{\pi}$ be the zero-section of the abelian
scheme $\pi$.
Then $0_{\pi}^{\ast} ( \mathcal{L} ) \cong \OO_{\mathfrak{B}'}$.
\end{enumerate}
\end{Remark}

Assume that $A$ is nowhere degenerate.
The next lemma shows 
that
the canonical height $\widehat{h}_L (X)$
is given by intersection on a suitable model.

\begin{Lemma} \label{lem:canonicalheightofsubvarnondeg}
Let $A$ be an abelian variety over $\overline{K}$ and let $L$ be an even
ample line bundle.
Let $X$ be a closed subvariety.
Assume that $A$ is nowhere degenerate.
Let $(\pi : \mathcal{A} \to \mathfrak{B}' , \mathcal{L} )$
be a model of $(A,L)$ as in Proposition~\ref{prop:model11}
and let $K'$ denote the function field of $\mathfrak{B}'$.
Assume that $X$ can be defined over $K'$
and let $\mathcal{X}$ be the closure of $X$ in $\mathcal{A}$.
Then we have
\[
\widehat{h}_L (X)
=
\frac{\deg_{\mathcal{H}'} \pi_{\ast} ( \cherncl_1 ( \mathcal{L})^{\cdot (d+1)} 
\cdot [ \mathcal{X}])}{ [K' : K]}
,
\]
where $\mathcal{H}'$ is the pullback of $\mathcal{H}$
to $\mathfrak{B}'$.
\end{Lemma}

Let us explain how Lemma~\ref{lem:canonicalheightofsubvarnondeg}
is verified under the assumption that $\dim ( \mathfrak{B} ) =1$.
Let $A$, $L$, and $X$ be as in Lemma~\ref{lem:canonicalheightofsubvarnondeg}.
Let $(\mathcal{A} ,\mathcal{L})$ be 
a model of $(A , L)$ as in Proposition~\ref{prop:model11};
see also Remark~\ref{rem:heightoverdimensin2}.
Recall that
in the limiting process to define the canonical height of subvarieties,
we consider the sequence (\ref{eq:limitcanonicalheightofvarieties}).
Here, we begin with $(\mathcal{A}_1 , \mathcal{L}_1) :=
(\mathcal{A} ,\mathcal{L})$.
Then we see from the condition of $(\mathcal{A} ,\mathcal{L})$
that $(\mathcal{A}_2 , \mathcal{L}_2)$
coincides with $(\mathcal{A}_1 , \mathcal{L}_1^{\otimes m^2})$.
By repeating the argument,
for any $n \in \NN$, we see that 
$(\mathcal{A}_n , \mathcal{L}_n)$
coincides with $(\mathcal{A}_1 , \mathcal{L}_1^{\otimes m^2(n-1)})$.
This shows that in this case
the sequence  (\ref{eq:limitcanonicalheightofvarieties})
is constant.
Thus Lemma~\ref{lem:canonicalheightofsubvarnondeg} holds.

\begin{Remark} \label{rem:heightoverdimensin1}
Under the setting of Remark~\ref{rem:heightoverdimensin2},
if $X$ can be defined over $K'$,
then we have
\[
\widehat{h}_L (X)
=
\frac{\deg ( \cherncl_1 ( \mathcal{L})^{\cdot (d+1)} 
\cdot [ \mathcal{X}])}{ [K' : K]}
\]
by
Lemma~\ref{lem:canonicalheightofsubvarnondeg},
where $d:=\dim (X)$. 
\end{Remark}

\subsection{Idea of the proof of Theorem~\ref{thm:yamaki7}}

We describe the idea of the proof of the assertion that 
(c) implies (b) in Theorem~\ref{thm:yamaki7}.
It suffices to show the following.
Let $A$ be a nowhere degenerate abelian variety.
Let $\mathfrak{t}$ be the image of the $\overline{K}/k$-trace homomorphism.
Suppose that the geometric Bogomolov conjecture holds for $A/\mathfrak{t}$.
Then it holds for $A$.

\subsubsection{Case of constant abelian varieties} \label{subsubsect:caseconstant}
First, we describe the idea of the proof of the above
assertion under the assumption that the
abelian variety is a constant
abelian variety.
Let $B$ be a constant abelian variety.
Let $Y$ be a closed subvariety of $B$.
Suppose that $Y$ has dense small points.
Then we want to show that $Y$ is a constant subvariety.
Remark that since $B$ is a constant abelian variety,
any torsion subvariety is a constant abelian variety
(cf. \cite[Proposition~3.7]{yamaki7}).

By definition,
we may write
$B = \widetilde{B} \otimes_k \overline{K}$
with some abelian variety $\widetilde{B}$ over $k$.
Let $\widetilde{M}$ be an even ample line bundle on $\widetilde{B}$.
Set $M := \widetilde{M} \otimes_k \overline{K}$, which is an even ample line bundle
on $B$.
We consider the standard model
$(\pi : \mathcal{B} \to \mathfrak{B}, \mathcal{M})$
of $(B,M)$ 
(cf. Example~\ref{ex:height0points}).
Then 
$( \mathcal{B} , \mathcal{M})$
satisfies the conditions
of Proposition~\ref{prop:model11}.
Since $Y$ has dense small points,
Proposition~\ref{prop:dense-height0}
tells us that $\widehat{h}_M ( Y ) = 0$.
By Lemma~\ref{lem:canonicalheightofsubvarnondeg},
it follows that
$\deg_{\mathcal{H}} \pi_{\ast}
( \cherncl_1 ( \mathcal{M})^{\cdot (d+1)} \cdot [\mathcal{Y}]) = 0$,
where $\mathcal{Y}$ is the closure of $Y$ in $\mathcal{B}$
and $d := \dim (Y)$.
Using this,
we can show that $Y$ is a constant subvariety
without difficulties;
see \cite[\S~3]{yamaki7} for the details.

\subsubsection{Reduction to the case of a product of a constant abelian variety
and a nowhere degenerate abelian variety with trivial trace} \label{subsect:rednnnndnd}

In this subsection, we remark that 
we may assume that 
$A$ is 
the direct product of a constant abelian variety and a
nowhere degenerate abelian variety with trivial $\overline{K}/k$-trace.

We begin by recalling the following basic fact.

\begin{Lemma} [cf. Corollary~6.7 of \cite{yamaki6}] \label{lem:GBCisogeny}
Let $\psi : A \to B$ be an isogeny of abelian varieties
over $\overline{K}$.
Then the geometric Bogomolov conjecture
holds for $A$ if and only if it holds for $B$.
\end{Lemma}

Let $A$ be a nowhere degenerate abelian variety over $\overline{K}$.
Then it is isogenous to the direct product of a constant abelian variety
and 
an abelian variety with trivial $\overline{K}/k$-trace.
Indeed,
by
the Poincar\'e complete reducibility theorem,
there exists a closed subvariety $C$ of $A$ such that
the natural homomorphism $\mathfrak{t} \times C \to A$
is an isogeny.
Note that $C$ is isogenous to $A / \mathfrak{t}$.
Furthermore, since the $\overline{K}/k$-trace homomorphism is finite,
we see that $A$ is isogenous to
$\left(
\widetilde{A}^{\overline{K} /k} 
\otimes_k
\overline{K}
\right)
\times A / \mathfrak{t}$.
Since $A / \mathfrak{t}$ has trivial $\overline{K}/k$-trace
(cf. \cite[Remark~5.4]{yamaki7}),
it follows that $A$ is isogenous to the product
of a constant abelian variety and a nowhere degenerate abelian variety
with trivial $\overline{K}/k$-trace.

Thus
by Lemma~\ref{lem:GBCisogeny}, 
our goal is the following assertion.
Let $A$ be a nowhere degenerate abelian variety over $\overline{K}$
with trivial $\overline{K}/k$-trace
and let $B = \widetilde{B} \otimes_k \overline{K}$
be a constant abelian variety.
Let $X$ be a closed subvariety of $B \times A$.
Assume that the geometric Bogomolov conjecture holds for $A$.
Suppose that $X$ has dense small points.
Then there exist a constant abelian subvariety $Y$ of $B$
and a torsion subvariety $T$ of $A$ such that $X = Y \times T$.
Let $\pr_{B} : B \times A \to B$ be the projection
and set $Y := \pr_B (X)$.
Since $X$ has dense small points, so does $Y$ by
Lemma~\ref{lem:density-quotient}.
By \S~\ref{subsubsect:caseconstant},
there exists a closed subvariety $\widetilde{Y}$
such that $Y = \widetilde{Y} \otimes_k \overline{K}$.
Thus we are reduced to showing the following proposition.

\begin{Proposition} \label{prop:onlytorsionpart}
Let $A$ be a nowhere degenerate abelian variety over $\overline{K}$
with trivial $\overline{K}/k$-trace.
Let $Y = \widetilde{Y} \otimes_k \overline{K}$
be a closed constant subvariety of a constant abelian variety
$B = \widetilde{B} \otimes_k \overline{K}$.
Let $p : Y \times A \to Y$ be the projection.
Let $X$ be a closed subvariety of $Y \times A$
with $p (X) = Y$.
Assume that the geometric Bogomolov conjecture holds for $A$.
Suppose that $X$ has dense small points.
Then there exists a torsion subvariety $T$ of $A$ such that $X = Y \times T$.
\end{Proposition}

\subsubsection{Relative height}

We keep the notation in Proposition~\ref{prop:onlytorsionpart}:
$A$ is a nowhere degenerate abelian variety, $B = \widetilde{B} \otimes_k
\overline{K}$ is a constant abelian variety,
$X$ is a closed subvariety of $B \times A$,
$Y = \widetilde{Y} \otimes_k \overline{K}$
is a constant closed subvariety of $B$,
and
$Y$ equals the image of $X$ by the projection $B \times A \to B$.

In the proof of Proposition~\ref{prop:onlytorsionpart},
we use the relative height.
Let $L$ be an even ample line bundle on $A$.
The relative height is a function $\mathbf{h}_{X/Y}^{L}$ defined over a
dense open subset of $Y(k)$
in the following way.
(In fact, it is defined over a dense open subset of $Y$, but we omit that; see \cite[\S~4.2]{yamaki7}
for the details.)
Take a point $\tilde{y} \in \widetilde{Y} (k)$.
Then $\tilde{y}$ 
can be naturally regarded  as a point of $Y ( \overline{K} )$,
and
we denote by $\overline{\tilde{y}_{K}}$
the corresponding point in $Y ( \overline{K} )$.
Let $X_{\overline{\tilde{y}_{K}}} := (p|_{X})^{-1} ( \overline{\tilde{y}_{K}} )$.
This is a closed subscheme of an abelian variety $p^{-1 } (
\overline{\tilde{y}_K}) = A$.
If $\tilde{y} \in \widetilde{Y} (k)$ is general,
$X_{\overline{\tilde{y}_{K}}} $ is of pure dimension
$d-e$, where $e:= \dim (Y)$.
Therefore, we consider the canonical height
$\widehat{h}_{L} ( X_{\overline{\tilde{y}_{K}}} )$,
and we set $\mathbf{h}_{X/Y}^{L}  ( \tilde{y})
:= \widehat{h}_{L} ( X_{\overline{\tilde{y}_{K}}} )$.

Since the relative height at $\tilde{y}$
is given by the height of the fiber $X_{\overline{\tilde{y}_{K}}}$,
we can describe it in terms of intersection products.
We take $K'$, $\mathfrak{B}'$
and
a model $(\mathcal{A} , \mathcal{L})$ over $\mathfrak{B}'$
of $(A, L)$
as in Proposition~\ref{prop:model11}.
Replacing $K'$ by a finite extension if necessary,
we may and do assume that $X$ can be defined over $K'$.
Let $\mathcal{B} \to \mathfrak{B}'$ be the standard
model of $B$,
that is,
$\mathcal{B} := \widetilde{B} \times_{\Spec (k)} \mathfrak{B}'$
and $\mathcal{B} \to \mathfrak{B}'$ is the canonical projection
(cf. Example~\ref{ex:height0points}).
Then $\mathcal{B}
\times_{\mathfrak{B}'}
\mathcal{A}$
is a model of $B \times A$
which satisfies the conditions
in Proposition~\ref{prop:model11}.
Let $\mathcal{X}$ be the closure of $X$ in $\mathcal{B}
\times_{\mathfrak{B}'}
\mathcal{A}$.
Set $\mathcal{Y} := \widetilde{Y} \times_{\Spec (k)} \mathfrak{B}'$.
Then $\mathcal{Y}$ is the model of $Y$ and
the canonical projection gives a surjective morphism
$\mathcal{X} \to \mathcal{Y}$.
Take a general $\tilde{y} \in \widetilde{Y} (k)$.
Remark that $\overline{\tilde{y}_{K}}$ is the 
geometric generic point of $\{ \tilde{y} \}
\times_{\Spec (k)} \mathfrak{B}' = \mathfrak{B}'$.
Then one sees that the fiber $\mathcal{X}_{\tilde{y}} 
\subset \{ \tilde{y} \} \times \mathcal{A} = \mathcal{A}$ 
of the composite 
$\mathcal{X} \to \mathcal{Y} \to \widetilde{Y}$
is a model of $X_{\overline{\tilde{y}_{K}}} \subset 
\{ \overline{\tilde{y}_K} \} \times A = A$.
By Lemma~\ref{lem:canonicalheightofsubvarnondeg},
we express the height of $X_{\overline{\tilde{y}_{K}}}$
by using the intersection with $\mathcal{X}_{\tilde{y}}$; indeed
\begin{align} \label{eq:relativeheightintersection}
\mathbf{h}^L_{X/Y} (\tilde{y}) =
\widehat{h}_L ( X_{\overline{\tilde{y}_{K}}} )
=
\frac{\deg_{\mathcal{H}'}
\pi_{\ast} \left(
\cherncl_1 ( \mathcal{L} )^{\cdot (d-e + 1)}
\cdot
[\mathcal{X}_{\tilde{y}}]
\right)}{[K':K]}
.
\end{align}

The relative height is used via the following lemma.

\begin{Lemma} [cf. Proposition~4.6 of
\cite{yamaki7}] \label{lem:proposition4.6yamaki7}
Under the setting above,
suppose that there exists a dense subset $S$ of $\widetilde{Y} (k)$
that satisfies the following two conditions.
\begin{enumerate}
\renewcommand{\labelenumi}{(\alph{enumi})}
\item
For any $\tilde{y}
\in S$,
we have
$\mathbf{h}_{X/Y}^{L}  ( \tilde{y}) = 0$.
\item
Any irreducible component
of $X_{\overline{\tilde{y}_{K}}}$
with its induced reduced subscheme structure is a torsion subvariety
of $\{ \tilde{y} \} \times A = A$.
\end{enumerate}
Then there exists a torsion subvariety $T$ of $A$
such that $X = Y \times T$.
\end{Lemma}

We omit the proof Lemma~\ref{lem:proposition4.6yamaki7},
and we refer to \cite{yamaki7} for the detail.
We just remark that
this lemma is obtained by showing a kind of rigidity
of torsion subvarieties in a nowhere degenerate abelian variety
with trivial $\overline{K}/k$-trace.

\subsubsection{Proof of Proposition~\ref{prop:onlytorsionpart}}

We give an outline of the proof of Proposition~\ref{prop:onlytorsionpart}.
First, assuming Lemma~\ref{lem:proposition5.1yamaki7} below,
we prove Proposition~\ref{prop:onlytorsionpart}.
Afterwards, we prove this lemma.

\begin{Lemma} [cf. Proposition~5.1 of \cite{yamaki7}] \label{lem:proposition5.1yamaki7}
Under the setting above,
suppose that $X$ has dense small points.
Then there exists a dense subset $S$ of $\widetilde{Y} (k)$
such that $\mathbf{h}_{X/Y}^{L}  ( \tilde{y}) = 0$
for any $\tilde{y} \in S$.
\end{Lemma}

We deduce 
Proposition~\ref{prop:onlytorsionpart}
from Lemmas~\ref{lem:proposition5.1yamaki7}.
Suppose that $X$ has dense small points.
Then by Lemma~\ref{lem:proposition5.1yamaki7}, 
there exists a dense subset $S$ of $\widetilde{Y} (k)$
such that $\mathbf{h}_{X/Y}^{L}  ( \tilde{y}) = 0$
for any $\tilde{y} \in S$.
This means that for any irreducible component $Z$ of 
$X_{\overline{\tilde{y}_K}}$, we have $\widehat{h}_L ( Z ) = 0$
(cf. Remark~\ref{rem:nonnegativecanonicalheight}).
It follows that $Z$ has dense small points.

Now assume that the geometric Bogomolov conjecture
holds for $A$.
Since $A$ has trivial $\overline{K}/k$-trace,
it follows that $Z$ is a torsion subvariety.
Thus we see that this $S$ satisfies conditions (a) and (b) 
in Lemma~\ref{lem:proposition4.6yamaki7}.
Therefore by this lemma,
we conclude the existence of a torsion subvariety $T$ as required.

\begin{Remark} \label{rem:keyoftheargument}
We keep the setting of Proposition~\ref{prop:onlytorsionpart}.
In the above argument,
we show that for general $\tilde{y} \in \widetilde{Y} (k)$,
any irreducible component $Z$ of the fiber $X_{\overline{\tilde{y}_K}}
\subset \{ \overline{\tilde{y}_K} \} \times A
= A$ has canonical height $0$.
\end{Remark}

Finally, 
we give an idea of the proof of Lemma~\ref{lem:proposition5.1yamaki7}.
To avoid technical difficulties, we assume that $\dim ( \mathfrak{B} ) = 1$.
The proof of Lemma~\ref{lem:proposition5.1yamaki7}
uses the expression of the canonical heights
in terms of intersections.
Recall that
$(\mathcal{A} , \mathcal{L})$ is a model
over $\mathfrak{B}'$
of $(A , L)$ as in Proposition~\ref{prop:model11}.
Let $\widetilde{M}$ be an even ample line bundle on $\widetilde{B}$
and set $M := \widetilde{M} \otimes_k \overline{K}$,
which is an even ample line bundle on $B = \widetilde{B} \otimes_k
\overline{K}$.
We assume that $\widetilde{M}$
is very ample, here.
Let
$( \mathcal{B} , \mathcal{M} )$ be the standard
model over $\mathfrak{B}'$ of $(B,M)$,
that is,
$\mathcal{B} := \widetilde{B} \times_{\Spec (k)}
\mathfrak{B}'$ 
and $\mathcal{M} := \widetilde{M} \otimes_{k} \OO_{\mathcal{B}'}$ 
(cf. Example~\ref{ex:height0points}).
Remark that $M \boxtimes L$ is an even ample line bundle on $B \times A$,
and
$(\mathcal{B} \times_{\mathfrak{B}'} \mathcal{A} \to 
\mathfrak{B}' , \mathcal{M} \boxtimes \mathcal{L})$
is a model 
over $\mathfrak{B}'$ of
$(B \times A , M \boxtimes L)$.
Remark also that the closure $\mathcal{X}$ of $X$ in 
$\mathcal{B} \times_{\mathfrak{B}'} \mathcal{A}$
is a model
of $X$.
Further,
$\widetilde{Y} \times_{\Spec (k)} \mathfrak{B}'$
is a model of $Y$
and 
equals the image of $\mathcal{X}$
by the projection $\mathcal{B} \times_{\mathfrak{B}'} \mathcal{A}
\to \mathcal{B}$.

Suppose that $X$ has dense small points.
Then by Proposition~\ref{prop:dense-height0},
we have $\widehat{h}_{M \boxtimes L} (X) = 0$.
By Lemma~\ref{lem:canonicalheightofsubvarnondeg}, 
it follows that
\addtocounter{Claim}{1}
\begin{align} \label{eq:vanishingwithmodels}
\deg_{\mathcal{H}'} \pi_{\ast} ( \cherncl_1 ( \mathcal{M} \boxtimes
\mathcal{L})^{\cdot (d+1)} 
\cdot [ \mathcal{X}])
=
0
.
\end{align}

We prove that
for any integer $e$ with $0 \leq e \leq d+1$,
\addtocounter{Claim}{1}
\begin{align} \label{eq:followingvanishing1}
\deg_{\mathcal{H}'} \pi_{\ast} ( 
\cherncl_1
(\pr_{\mathcal{A}}^{\ast} (\mathcal{L}))^{\cdot (d -e +1)} 
\cdot (
\cherncl_1 ( \pr_{\mathcal{B}}^{\ast} (\mathcal{M})
)^{e} 
\cdot [ \mathcal{X}]))
=
0,
\end{align}
where $\pr_{\mathcal{A}} : \mathcal{B} \times_{\mathfrak{B}'} \mathcal{A}
\to \mathcal{A}$
and $\pr_{\mathcal{B}} : \mathcal{B} \times_{\mathfrak{B}'} \mathcal{A}
\to \mathcal{B}$
are the canonical projections.
First, note that the left-hand side in (\ref{eq:followingvanishing1})
is non-negative.
Indeed,
it is, up to a positive multiple, the canonical height of $X$
with respect to $(d -e +1)$-copies of $\OO_B \boxtimes L$
and $e$-copies of $M \boxtimes \OO_A$
(cf. \cite[Theorem~3.5~(d)]{gubler3}),
and this canonical height is non-negative by 
(\cite[Theorem~11.18~(e)]{gubler0}).
Here, we have
\[
\deg_{\mathcal{H}'} \pi_{\ast} ( \cherncl_1 ( \mathcal{M} \boxtimes
\mathcal{L})^{\cdot (d+1)} 
\cdot [ \mathcal{X}])
=
\sum_{e = 0}^{d+1}
\begin{pmatrix}
d+1 \\ e
\end{pmatrix}
\deg_{\mathcal{H}'} \pi_{\ast} ( 
\cherncl_1
(\pr_{\mathcal{A}}^{\ast} (\mathcal{L}))^{\cdot (d -e + 1)} 
\cdot
(
\cherncl_1 ( \pr_{\mathcal{B}}^{\ast} (\mathcal{M})
)^{e} 
\cdot [ \mathcal{X}]))
,
\]
which equals $0$ by (\ref{eq:vanishingwithmodels}).
Then we obtain
(\ref{eq:followingvanishing1}).

We consider the composite morphism 
$\varphi : \mathcal{B} \times_{\mathfrak{B}'} \mathcal{A}
\to \mathcal{B} \to \widetilde{B}$.
This gives us
$\varphi|_{\mathcal{X} } : \mathcal{X} \to 
\widetilde{Y}$
by restriction.
We note
$\pr_{\mathcal{B}}^{\ast} (\mathcal{M}) = \varphi^{\ast} (\widetilde{M} )$.
Here,
recall that $\widetilde{M}$ is very ample.
Then
there exists a dense open subset $S \subset \widetilde{Y}(k)$
such that for any $\tilde{y} \in S$,
there exist
a finite number of points $\tilde{y}_1 ,
\ldots , \tilde{y}_m$ with $\tilde{y}_1 = \tilde{y}$
such that 
\[
\cherncl_1 ( \widetilde{M}
)^{e} \cdot [ \widetilde{Y} ]
= \sum_{i = 1}^{m}
[  \tilde{y}_i   ]
\]
as cycles classes
and such that
$\varphi|_{\mathcal{X}}$ is flat over $\tilde{y}_i$ for any $i=1,\ldots,m$.
Put $\mathcal{X}_{\tilde{y}_{i}} := (\varphi|_{\mathcal{X} })^{-1} (
\tilde{y}_i) = \varphi^{-1} ( \tilde{y}_i ) \cap \mathcal{X}$. 
Then by the choice of $\tilde{y}_1 
,\ldots , \tilde{y}_m$,
\[
\cherncl_1 ( \pr_{\mathcal{B}}^{\ast} (\mathcal{M})
)^{e} 
\cdot [ \mathcal{X}]
=
\sum_{i=1}^{m}
[\mathcal{X}_{\tilde{y}_{i}}]
,
\]
and hence by (\ref{eq:followingvanishing1}),
\[
\deg_{\mathcal{H}'} \pi_{\ast} 
\left( 
\cherncl_1
(\pr_{\mathcal{A}}^{\ast} (\mathcal{L}))^{\cdot (d -e + 1)} 
\cdot
\sum_{i=1}^{m}
[\mathcal{X}_{\tilde{y}_{i}}]
\right)
=0
.
\]
Since 
\[
\deg_{\mathcal{H}'} \pi_{\ast} ( 
\cherncl_1
(\pr_{\mathcal{A}}^{\ast} ( \mathcal{L} ))^{\cdot (d -e + 1)} 
\cdot
[ ( \mathcal{X}_{\tilde{y}_{i}} ])
\geq 0
\]
by the same reason as in the proof of (\ref{eq:followingvanishing1}),
it follows that
\[
\deg_{\mathcal{H}'} \pi_{\ast} ( 
\cherncl_1
(\pr_{\mathcal{A}}^{\ast} ( \mathcal{L} ))^{\cdot (d -e + 1)} 
\cdot
[ ( \mathcal{X}_{\tilde{y}_{i}} ])
=0
\]
for any $i = 1 , \ldots , m$.
Since $\tilde{y}_1 = \tilde{y}$,
\[
\deg_{\mathcal{H}'} \pi_{\ast} ( 
\cherncl_1
(\pr_{\mathcal{A}}^{\ast} ( \mathcal{L} ))^{\cdot (d -e + 1)} 
\cdot
[\mathcal{X}_{\tilde{y}} ])
=0
\]
holds, in particular.
Thus by (\ref{eq:relativeheightintersection}),
we obtain
$\mathbf{h}_{X/Y}^{L}  ( \tilde{y}) 
= 0$.
This proves Lemma~\ref{lem:proposition5.1yamaki7}.

\section{Proof of the conjecture for curves} \label{sect:proofforcurve}

In this section,
we give an idea of the proof of the 
geometric Bogomolov
conjecture for curves.
We also give a remark on the conjecture for abelian varieties.

\subsection{For curves and for divisors} \label{subsect:CD}

As we noted before, 
we actually prove 
in \cite{yamaki8}
the
following more general result.

\begin{Theorem} [Theorem~1.3 of \cite{yamaki8}] \label{thm:GBC_dimX=1}
Let $A$ be an abelian variety over $\overline{K}$
and let $X$ be a closed subvariety of $A$.
Assume that $\dim (X) = 1$.
Suppose that $X$ has dense small points.
Then $X$ is a special subvariety.
\end{Theorem}

It should be remarked that
the above theorem follows from the following.

\begin{Theorem} [Theorem~1.4 of \cite{yamaki8}] \label{thm:GBC_codimX=1}
Let $A$ be an abelian variety over $\overline{K}$.
Let $X$ be a closed subvariety of $A$.
Assume that $\codim (X,A) = 1$.
Suppose that $X$ has dense small points.
Then $X$ is a special subvariety.
\end{Theorem}

The principle to connect Theorem~\ref{thm:GBC_dimX=1}
to Theorem~\ref{thm:GBC_codimX=1}
is as follows:
if we take the sum of some copies of a curve in
an abelian variety,
then it will be a divisor of the abelian variety.
Let us describe the idea a little more precisely.
Some argument using
Theorem~\ref{thm:yamaki7}
shows that
we may assume that $A$ has trivial $\overline{K}/k$-trace
to show Theorem~\ref{thm:GBC_codimX=1}.
Now 
suppose that $X$ is a curve in $A$ with dense small points.
Since the $\overline{K}/k$-trace of $A$ is trivial,
our goal is to show that $X$ is the
translate of an abelian subvariety of dimension $1$
by a torsion point.
Consider $X - \tau$ for any $\tau \in A ( \overline{K} )_{tor}$.
We take a $\tau_1$ in such a way that
the dimension of the minimal
abelian subvariety 
containing
$X - \tau_1$ is minimal.
Set $X_1 := X - \tau_1$
and let
$A_{1}$ be the minimal abelian subvariety containing $X_1$.
Note that $X_1$ has dense small points,
since $\tau_1$ is a torsion point.
For a nonnegative integer $l$, we consider
$Z_l:= \sum_{j=0}^{l} (-1)^{j}X_1$, 
where the sum means the addition of the abelian variety.
Then
$Z$ is a closed subvariety of $A_1$,
and
since $X_1$
has dense small points, so does $Z$.
Furthermore, one sees that $\codim (Z_l , A_1) = 1$ for some
$l$.
Thus we reach the setting of Theorem~\ref{thm:GBC_codimX=1},
and a few more arguments leads to the conclusion.
See \cite[Proof of Theorem~5.11]{yamaki8} for the detail.

\subsection{Three steps of the proof}

By the argument in \S~\ref{subsect:CD},
the geometric Bogomolov conjecture for curves is
reduced to Theorem~\ref{thm:GBC_codimX=1}.
In this section, we give an outline of the proof of this theorem.
The proof
consists of three steps,
and each of them basically corresponds to \cite{yamaki6}, \cite{yamaki7}, and
\cite{yamaki8}, respectively.

\subsubsection*{Step~1}
We reduce 
Theorem~\ref{thm:GBC_codimX=1} to the case where 
the abelian variety is nowhere degenerate, i.e., to the following proposition.

\begin{Proposition} \label{prop:redtonowheredegenerate}
Assume that $A$ is nowhere degenerate.
Let $X$ be a closed subvariety of $A$ of codimension $1$.
Suppose that $X$ has dense small points.
Then $X$ is a special subvariety.
\end{Proposition}

We explain how
the theorem
follows from Theorem~\ref{thm:yamaki6_7.21} 
and Proposition~\ref{prop:redtonowheredegenerate}.
Let $A$ be any abelian variety over $\overline{K}$.
Let $\mathfrak{m}$ be the maximal nowhere degenerate abelian 
subvariety of $A$.
Recall that there exists a surjective homomorphism
$\phi : A \to \mathfrak{m}$;
see the argument just after Theorem~\ref{thm:yamaki6_7.21}.
Let $X$ be a closed subvariety of $A$ of codimension $1$
with dense small points.
Put $Y := \phi (X)$.
Then
$Y =\mathfrak{m}$ or 
$\codim ( Y , \mathfrak{m}) =1$.
If  
$Y = \mathfrak{m}$,
then
$Y$ is special.
If
$\codim ( Y , \mathfrak{m}) = 1$,
then applying Proposition~\ref{prop:redtonowheredegenerate}
to $\mathfrak{m}$ and $Y$ (in place of $A$ and $X$ respectively),
we find that $Y$ is special.
In any case,  Theorem~\ref{thm:yamaki6_7.21} shows that
$X$ is a special subvariety.

\subsubsection*{Step~2}

By the argument in Step~1,
we are reduced to showing Proposition~\ref{prop:redtonowheredegenerate}.
In this step, we reduce Proposition~\ref{prop:redtonowheredegenerate}
to the assertion for nowhere degenerate abelian varieties
with trivial $\overline{K}/k$-trace.
As we noted in \S~\ref{subsect:rednnnndnd},
any nowhere degenerate abelian variety
is isogenous to the product of a constant abelian variety
and a nowhere degenerate abelian variety with trivial $\overline{K}/k$-trace.
Therefore, 
by Lemma~\ref{lem:GBCisogeny},
it suffices to show the following:

\begin{Proposition} \label{prop:step2}
Let $A$ be a nowhere degenerate abelian variety over $\overline{K}$
with trivial $\overline{K}/k$-trace,
$B$ a constant abelian variety,
and let $X \subset B \times A$ be a closed subvariety of codimension
$1$.
Suppose that $X$ has dense small points.
Then $X$ is a special subvariety.
\end{Proposition}

Proposition~\ref{prop:step2} follows from
the proposition below.

\begin{Proposition} \label{prop:tothecasendavwtt}
Let $A$ be a nowhere degenerate abelian variety over $\overline{K}$
with trivial $\overline{K}/k$-trace
and let $X$ be a closed subvariety of $A$ of codimension $1$.
Then $X$ does not have dense small points.
\end{Proposition}

We show that Proposition~\ref{prop:step2}
is deduced form Proposition~\ref{prop:tothecasendavwtt}.
Let $Y$ be the image of $X$ by the projection $B  \times A \to B$.
Since $X$ has dense small points, do does $Y$,
and hence by the argument of \S~\ref{subsubsect:caseconstant},
$Y$ is a constant subvariety of $B$.

We prove that $\dim (Y) < \dim (B)$ by contradiction.
Suppose that $\dim (Y) = \dim (B)$,
i.e., $Y = B$.
Write $B = \widetilde{B} \otimes_k \overline{K}$.
As is noted in Remark~\ref{rem:keyoftheargument},
for a general $\tilde{y} \in \widetilde{B} (k)$,
any irreducible component $Z$ of $X_{\overline{\tilde{y}_K}}$
has canonical height $0$.
Thus by Proposition~\ref{prop:dense-height0},
$Z$ has dense small points.
On the other hand,
the dimension counting shows that $Z$ is a divisor 
on $\{ \overline{\tilde{y}_K} \} \times A = A$
for general $\tilde{y}$.
That contradicts Proposition~\ref{prop:tothecasendavwtt}.

Thus we have $\dim (Y) < \dim (B)$.
Then $Y$ has codimension $1$ in $B$.
Since $X$ has codimension $1$ in $B \times A$,
it follows that $X = Y \times A$.
Since $Y$ is a constant subvariety,
this proves that
$X$ is a special subvariety of $B \times A$.

\subsubsection*{Step~3}

Now, our goal is Proposition~\ref{prop:tothecasendavwtt}.
We see that 
this
follows from the proposition below.

\begin{Proposition} \label{prop:ampledivisorhaspositiveheight}
Let $A$ be a nowhere degenerate abelian variety 
over $\overline{K}$ with trivial $\overline{K}/k$-trace,
and let $X$ be an irreducible effective ample divisor
on $A$ of codimension $1$.
Set $D := X + [-1]^{\ast} (X)$
and $L := \OO_{A} ( D )$,
where $+$ is the addition of divisors.
Then $L$ is an even ample line bundle on $A$,
and we have $\widehat{h}_L(X) > 0$.
\end{Proposition}

Let us check that Proposition~\ref{prop:tothecasendavwtt}
is deduced from Proposition~\ref{prop:ampledivisorhaspositiveheight}.
Let $A$ be a nowhere degenerate abelian variety 
over $\overline{K}$ with trivial $\overline{K}/k$-trace
and let $X$ be a closed subvariety of $A$ of codimension $1$.
Note that any effective divisor is the pullback of an
ample divisor by some homomorphism,
which follows from \cite[p.88, Remarks on effective divisors by Nori]{mumford};
see \cite[Step~1 of the proof of Theorem~5.7]{yamaki8}
for the precise argument.
This means that there exists a surjective homomorphism
$\phi : A \to A'$
and an effective ample divisor $X'$ on $A'$
such that $X = \phi^{-1} (X')$.
By Proposition~\ref{prop:ampledivisorhaspositiveheight}
together
with Proposition~\ref{prop:dense-height0},
$X'$ does not have dense small points.
By Lemma~\ref{lem:density-quotient},
it follows that $X$ does not have dense small points.

We explain the
idea of the proof of 
Proposition~\ref{prop:ampledivisorhaspositiveheight}.
It is obvious that $L$ in the proposition is even and ample.
Now,
to avoid technical difficulties,
we make the following assumptions:
\begin{enumerate}
\renewcommand{\labelenumi}{(\alph{enumi})}
\item
$\dim (\mathfrak{B}) = 1$;
\item
there exists an abelian scheme $\pi : \mathcal{A} \to \mathfrak{B}$
with zero-section $0_\pi$
and with
geometric generic fiber $A$,
and
$X$ is defined over $K$;
\item
$0 \notin X$;
\item
$\# k > \aleph_0$,
i.e.,
$k$ has uncountably infinite cardinality.
\end{enumerate}
In fact, 
it is not very difficult to see that
we may assume
(b), (c), and (d) without loss of generality;
we make assumption (a) 
to avoid technical difficulties.

Let us give a sketch of the proof
of Proposition~\ref{prop:ampledivisorhaspositiveheight}.
(This is based on the argument in \cite[\S~1.3]{yamaki8}.)
Since $L$ is even,
we have
$\widehat{h}_{L} ( D) = 2 \widehat{h}_L (X)$.
Let $\mathcal{D}$ be the closure of $D$ in $\mathcal{A}$.
Since $0 \notin D$,
$0_\pi^{\ast} ( \mathcal{D} )$ is a well-defined 
effective divisor on $\mathfrak{B}$.
Set $\mathcal{N} := \OO_{\mathfrak{B}}
\left( 0_\pi^{\ast} ( \mathcal{D})
\right)$
and
$\mathcal{L} := \OO_{\mathcal{A}} \left( \mathcal{D}
\right) 
\otimes
\pi^{\ast} ( \mathcal{N}^{\otimes (-1)} )$.
Then $0_{\pi}^{\ast} ( \mathcal{L} ) = \OO_{\mathfrak{B}}$.
Put $n := \dim (A)$.
Then
one sees that
\[
\widehat{h}_{L} ( D ) = \deg
\left( 
\cherncl_1 ( \mathcal{L} )^{\cdot n}
\cdot  [ \mathcal{D} ] \right)
=
\deg
\left( 
\cherncl_1 ( \mathcal{L} )^{\cdot n}
\cdot  \cherncl_1 \left(
\OO_{\mathcal{A}}
\left(
\mathcal{D}
\right)
\right)
\cdot
[ \mathcal{A} ]
\right)
.
\]
Since $A$ has canonical height $0$
(cf. Remark~\ref{rem:heightofablianvariety}),
we have $\deg \left( \cherncl_1 ( \mathcal{L} )^{\cdot (n + 1)}
\cdot  [ \mathcal{A} ] \right) = 0$
(cf. Remark~\ref{rem:heightoverdimensin1}).
It follows that 
\begin{align*}
\deg
\left( 
\cherncl_1 ( \mathcal{L} )^{\cdot n}
\cdot  \cherncl_1 \left(
\OO_{\mathcal{A}}
\left(
\mathcal{D}
\right)
\right)
\cdot
[ \mathcal{A} ]
\right)
&=
\deg
\left( 
\cherncl_1 ( \mathcal{L} )^{\cdot (n+1)}
\cdot [
\mathcal{A}
]
\right)
+
\deg
\left( 
\cherncl_1 ( \mathcal{L} )^{\cdot n}
\cdot 
\pi^{\ast} \cherncl_1 ( \mathcal{N})
\cdot
[
\mathcal{A}
]
\right)
\\
&=
\deg_L (A) \cdot \deg (\mathcal{N})
.
\end{align*}
Since $\deg_L (A) > 0$ by the ampleness of $L$ on $A$,
it remains to show $\deg (\mathcal{N}) > 0$.
In fact, we will see below
that $\mathcal{N}$ is non-trivial.
Since $
0_\pi^{\ast} ( \mathcal{D})$ is an effective
divisor,
we conclude
$\deg (\mathcal{N}) > 0$.

The outline of the proof of the non-triviality
of $\mathcal{N}$
is as follows.
We prove the non-triviality by contradiction.
Suppose that it is trivial.
Then we can show that 
there exists a finite covering $\mathfrak{B}' \to \mathfrak{B}$
such that
the complete linear system $| 2 \mathcal{D}' |$
on $\mathcal{A}'$
is base-point free
(cf. \cite[Proposition~4.4]{yamaki8}),
where $\pi' : \mathcal{A}' \to \mathfrak{B}'$ and $\mathcal{D}'$
are the base-change of $\pi$ and $\mathcal{D}$ by this
$\mathfrak{B}' \to \mathfrak{B}$, respectively.
Let
$\varphi : \mathcal{A}' \to Z$
be the surjective morphism
associated to $| 2 \mathcal{D}' |$,
where $Z$ is a closed subvariety of the dual space of $|2 \mathcal{D}' |$.
Remark that for any curve $\gamma \subset \mathcal{A}'$,
$\deg \left(
\cherncl_1 ( \mathcal{L}' ) \cdot \gamma
\right) = 0$
if and only if $\varphi ( \gamma )$ is a point,
where $\mathcal{L}' = \OO_{\mathcal{A}'}( \mathcal{D}')$,
which is 
the pull-back of $\mathcal{L}$ to $\mathcal{A}'$.
Further, remark that for any $a \in A \left( 
\overline{K} \right)$,
$\widehat{h}_{L} ( a ) = 0$
if and only  if $\deg \left(
\cherncl_1 ( \mathcal{L}' ) \cdot \Delta_{a}
\right) = 0$,
where $\Delta_{a}$ is the closure of $a$
in $\mathcal{A}'$
(cf. (\ref{eq:heightnondegenerate})).

Let $\Gamma$ be the set of irreducible curves in $\mathcal{A}'$ such that
$\dim (\varphi ( \gamma)) = 0$.
We set 
\[
A(0 ; L) := \left\{ a \in A 
\left( \overline{K} \right)
\left| 
\ 
\widehat{h}_L (a) = 0 \right.
\right\}.
\]
The above argument
shows that
$\Gamma \supset \{ \Delta_a \mid  a \in A (0;L) \}$.
Since $L^{\otimes 2}$ is ample
and $\pi'$ is an abelian scheme,
$\left( \mathcal{L}'\right)^{\otimes 2}$ 
is relatively ample with respect to $\pi'$.
It follows that for any curve $\gamma$ in a fiber of $\pi'$,
we have $\deg \left(
\cherncl_1 ( \mathcal{L}' ) \cdot \gamma
\right) > 0$, and hence
$\varphi$ is finite on any fiber of $\pi'$.
This means that if $\gamma \in \Gamma$,
then $\gamma$ is flat over 
$\mathfrak{B}'$.
Noting that an irreducible curve in $\mathcal{A}'$ that is flat over 
$\mathfrak{B}'$ is the closure of some point in $A 
\left( \overline{K} \right)$,
we have
$\Gamma \subset
 \left\{ \Delta_a \mid  a \in A \left(
\overline{K}\right) \right\} 
$.
Further,
since $\dim (\varphi ( \gamma)) = 0$
for any $\gamma \in \Gamma$,
it follows from what we note in the previous paragraph that 
$\Gamma \subset
\{ \Delta_a \mid  a \in A (0;L) \} 
$.
Thus 
$\Gamma =
\{ \Delta_a \mid  a \in A (0;L) \} 
$.

On the other hand,
since
$A (0;L)$ 
is dense in $A$,
the set
$\Gamma
=
\{ \Delta_a \mid a \in A (0;L) \}
$
is dense in $\mathcal{A}'$.
Since $\dim ( \varphi (\gamma)) = 0$ for any $\gamma \in \Gamma$,
it follows that
$\varphi$ is not generically finite 
on $\mathcal{A}'$.
Here, recall that $\# k > \aleph_0$.
Then we have 
$\# \Gamma > \aleph_0$.
Note that the map $A \left( \overline{K} \right) \to 
\{ \Delta_a \mid a \in A \left( \overline{K} \right) \}$
given by $a \mapsto \Delta_a$ induces a bijection
$A( 0 ; L) \to \{ \Delta_a \mid a \in A (0;L) \} = \Gamma$.
Then
the above shows that
$A \left(
\overline{K}
\right)$ has uncountably many
points of height $0$.
However,
since $A$ has trivial $\overline{K}/k$-trace,
a point of $A \left( 
\overline{K}
\right)$ has height $0$ 
if and only if it is torsion
(cf. Proposition~\ref{prop:height0points_ff}),
and there are only countably many such points.
This is a contradiction.
Thus Proposition~\ref{prop:ampledivisorhaspositiveheight} 
is proved.

\subsection{Application to the conjecture for abelian varieties}
Here are remarks on some contributions to the geometric Bogomolov
conjecture for abelian varieties.
Using Theorems~\ref{thm:GBC_codimX=1}
and \ref{thm:GBC_dimX=1},
we obtain the following results on the conjecture
for abelian varieties.

\begin{Theorem} [Corollary~6.4 of \cite{yamaki8}] \label{thm:GBCyamaki8}
Let $A$ be an abelian variety over $\overline{K}$
and let $\mathfrak{m}$ be the maximal nowhere degenerate
abelian subvariety of $A$. Let $\mathfrak{t}$ be the image of
the $\overline{K}/k$-homomorphism.
Assume that $\dim ( \mathfrak{m} / \mathfrak{t}) \leq 3$.
Then the geometric Bogomolov conjecture holds for $A$.
\end{Theorem}

Indeed, since $\dim ( \mathfrak{m} / \mathfrak{t}) \leq 3$,
it follows by Remark~\ref{rem:trivialcaseforGBC},
Theorems~\ref{thm:GBC_codimX=1}
and \ref{thm:GBC_dimX=1}
that the geometric Bogomolov conjecture holds
for $\mathfrak{m} / \mathfrak{t}$.
By Theorem~\ref{thm:yamaki7}, the geometric Bogomolov conjecture holds
for $A$.

In conclusion,
the geometric Bogomolov conjecture is reduced to the
the following conjecture.

\begin{Conjecture}
Let $A$ be a nowhere degenerate abelian variety 
over $\overline{K}$ with trivial $\overline{K}/k$-trace.
Let $X$ be a closed subvariety of $A$.
Assume that $2 \leq \dim (X) \leq \dim (A) - 2$.
Suppose that $X$ has dense small points.
Then $X$ is a torsion subvariety.
\end{Conjecture}

This is not what we discuss in this paper,
but
we have shown in \cite{yamaki9} that
the geometric Bogomolov conjecture holds for nowhere degenerate
abelian varieties of dimension $5$ with trivial $\overline{K}/k$-trace
(cf. \cite[Theorem~1.3]{yamaki9}),
and thus we have the following.

\begin{Theorem} [Theorem~1.4 of \cite{yamaki9}] \label{thm:yamaki9}
Let $A$, $\mathfrak{m}$, and
$\mathfrak{t}$ be as in Theorem~\ref{thm:GBCyamaki8}.
Assume that $\dim ( \mathfrak{m} / \mathfrak{t}) = 5$.
Then the geometric Bogomolov conjecture holds for $A$.
\end{Theorem}

It would be interesting to obtain the result for 
an abelian variety $A$ as above
with
$\dim (A) = 4$
and $X$ with $\dim (X) = 2$.

\section{Manin--Mumford conjecture in positive characteristic} \label{sect:MM}
In this section,
we give a remark on the relationship between the geometric Bogomolov
conjecture and the Manin--Mumford conjecture in positive characteristic.

Let $\mathfrak{K}$ be an algebraically closed field.
Let $A$ be an abelian variety over $\mathfrak{K}$.
Manin--Mumford conjecture,
which is Raynaud's theorem, asserts
that if $\ch (\mathfrak{K}) = 0$,
then any closed subvariety of $A$
that has dense torsion points
is a torsion subvariety
(cf. \S~\ref{subsubsect:raynaud}).

In the case of $\ch (\mathfrak{K}) > 0$,
on the other hand,
the same assertion does not hold.
If $\mathfrak{K}$ is an algebraic closure of a finite field
and $X$ is any closed subvariety of $A$,
then any $\mathfrak{K}$-point of $X$ is a torsion point,
and thus $X$ always has dense torsion points.

However,
up to such influence 
that stems from finite fields,
one may expect 
that a similar statement should also hold in positive characteristics.
In fact, the following precise result
is known.

\begin{citeTheorem} \label{thm:MM}
Assume that $\ch (\mathfrak{K}) > 0$.
Let $k$ be the algebraic closure in $\mathfrak{K}$ 
of the prime field of $\mathfrak{K}$.
Let $A$ be an abelian variety over $\mathfrak{K}$.
Let $X$ be a closed subvariety of $A$.
If $X$ has dense torsion points,
then there exist an abelian subvariety $G$ of $A$,
a closed subvariety $\widetilde{Y}$ of $\widetilde{A}^{\mathfrak{K}/k}$
and a torsion point $\tau \in A \left( \mathfrak{K} \right)$ 
such that $X = G + \Tr^{\mathfrak{K}/k}_{A} \left( \widetilde{Y}
\otimes_k \mathfrak{K} \right) + \tau$,
where $\left( 
\widetilde{A}^{\mathfrak{K}/k} , 
\Tr^{\mathfrak{K}/k}_{A}
\right)$
is the $\mathfrak{K}/k$-trace of $A$.
\end{citeTheorem}

Here, remark that the $\mathfrak{K}/k$-trace
of $A$ is defined by the same way as the $\overline{K}/k$-trace
in \S~\ref{subsection:height0}
(in place of $\overline{K}$ with $\mathfrak{K}$);
it is known that 
there exists a unique $\mathfrak{K}/k$-trace of $A$
(cf. \cite[Ch.VIII, \S~3]{lang1}).

This theorem is due to the following authors:
In 2001,
Scanlon gave a sketch of the model-theoretic proof of this theorem
(\cite{scanlon0}).
In 2004,
Pink and Roessler gave an algebro-geometric proof
(\cite{PR}).
In 2005, Scanlon gave a detailed model-theoretic proof
in \cite{scanlon1}
based on the argument in \cite{scanlon0}.
Note that
in those papers, they prove a generalized version 
for semiabelian varieties $A$, in fact.

Here, we explain that Theorem~\ref{thm:MM} can be deduced from the geometric
Bogomolov conjecture for $A$, as Moriwaki mentioned
in \cite{moriwaki5}
in the case of characteristic $0$.
Let $A$ be an abelian variety over $\mathfrak{K}$ 
and let $X$ be a closed subvariety
of $A$.
Then there exist $t_1 , \ldots , t_n \in \mathfrak{K}$
such that $A$ and $X$ can be defined over $K := k ( t_1 , \ldots , t_n)$,
that is,
there exist an abelian variety $A_0$ over $K$
and a closed subvariety $X_0$ of $A_0$ such that
$A = A_0 \otimes_K \mathfrak{K}$ and $X = X_0 \otimes_K \mathfrak{K}$.
Further,
there exists a normal projective variety $\mathfrak{B}$ over $k$
with function field $K$.
Let $\overline{K}$ be the algebraic closure of $K$ in $\mathfrak{K}$
and set $A_{\overline{K}} := A_0 \otimes_K \overline{K}$
and $X_{\overline{K}} := X_0 \otimes_K \overline{K}$.
Let $\mathcal{H}$ be an ample line bundle on $\mathfrak{B}$.
Then we have a notion of height
over $K$, and we can consider the canonical height
on 
$A_{\overline{K}}$
associated to an even ample line bundle.
Suppose that $X$ has dense torsion points.
Then $X_{\overline{K}}$ has dense torsion points
and hence has dense small points.
If we assume the geometric Bogomolov conjecture,
it follows that $X_{\overline{K}}$ is a special subvariety,
and this implies the conclusion of Theorem~\ref{thm:MM}.

Since 
the geometric Bogomolov conjecture 
is still open, the above argument 
does not give a new proof of Theorem~\ref{thm:MM}.
However,
we can actually deduce
the theorem in the following special cases:
\begin{enumerate}
\item
$\dim (X) = 1$ or $\codim (X , A) = 1$
(from 
Theorems~\ref{thm:GBC_codimX=1}
and \ref{thm:GBC_dimX=1});
\item
$\dim (A) \leq 3$
(as a consequence above).
\end{enumerate}
In particular, we have recovered the positive characteristic version
of Theorem~\ref{thm:raynaud-curve}.


\renewcommand{\thesection}{Appendix} 
\renewcommand{\theTheorem}{A.\arabic{Theorem}}
\renewcommand{\theClaim}{A.\arabic{Theorem}.\arabic{Claim}}
\renewcommand{\theequation}{A.\arabic{Theorem}.\arabic{Claim}}
\renewcommand{\theProposition}
{A.\arabic{Theorem}.\arabic{Proposition}}
\renewcommand{\theLemma}{A.\arabic{Theorem}.\arabic{Lemma}}
\setcounter{section}{0}
\renewcommand{\thesubsection}{A.\arabic{subsection}}

\section{Admissible formal schemes and the Raynaud generic fibers}
In this appendix,
we give a brief summary on admissible formal schemes and
the Raynaud generic fibers.
Further, we explain the Raynaud extension and the valuation map
for an abelian variety.
\subsection{Admissible formal schemes and Raynaud generic fibers}

First of all, we recall the notion of affinoid algebras
and associated affinoid spaces in the sense of Berkovich.
Let $\KK \langle x_{1} , \ldots , x_n \rangle$
be the Tate algebra 
over $\KK$, that is, 
the completion of the polynomial ring 
$\KK [ x_{1} , \ldots , x_n ]$
with respect to the Gauss norm.
By definition, $\KK \langle x_{1} , \ldots , x_n \rangle$
equals 
\[
\left\{
\left.
\sum_{\mathbf{m} = (m_1 , \ldots , m_n) \in \ZZ_{\geq 0}^n} a_{\mathbf{m}}
x_1^{m_1} \cdots x_n^{m_n} 
\in
\KK [[  x_{1} , \ldots , x_n ]]
\ 
\right|
\lim_{m_1 + \cdots + m_n \to \infty} |a_{\mathbf{m}}| = 0
\right\}
.
\]
A \emph{$\KK$-affinoid algebra}
is a
$\KK$-algebra isomorphic to 
$\KK \langle x_{1} , \ldots , x_n \rangle / I$
for some ideal $I$ of $\KK \langle x_{1} , \ldots , x_n \rangle$.
Let $\mathrm{Max} (R)$ be the maximal spectrum
of a $\KK$-affinoid algebra $R$, that is, the set of maximal ideals of $R$.
For each $p \in \mathrm{Max} (R)$, the residue field at $p$
is canonically isomorphic to $\KK$, and thus 
it
is endowed with
a norm.
Therefore,
we can consider the supremum semi-norm
$|\cdot|_{\sup} : R \to \RR$
over $\mathrm{Max} (R)$.
The \emph{Berkovich spectrum 
$\mathcal{M}(R)$
of $R$}
is
the set of multiplicative seminorms 
$R$
bounded with $|\cdot|_{\sup}$
endowed with the weakest topology such that
for any $f \in R$,
the function $\mathcal{M} (R) \to \RR$
given by $p \mapsto p(f) =: |f(p)|$ is continuous.
The Berkovich spectrum of $R$
is also called the Berkovich affinoid space
associated to $R$.
A Berkovich affinoid  $\mathcal{M} (R)$ 
have a non-archimedean analytic structure
whose ring of functions $R$,
but we do not explain it here.

A $\KK^{\circ}$-algebra 
is called an
 \emph{admissible $\KK^{\circ}$-algebra} if it does not have any
$\KK^{\circ}$-torsions and it is isomorphic to
$\KK^{\circ} \langle x_{1} , \ldots , x_{n} \rangle / I$
for some $n \in \NN$ and for some ideal
$I$ of $\KK^{\circ} \langle x_{1} , \ldots , x_{n} \rangle$.
Note that an admissible $\KK^{\circ}$-algebra is flat over $\KK^{\circ}$.
The formal spectrum 
of an admissible $\KK^{\circ}$-algebra
is called an \emph{affine admissible formal scheme}.

For an admissible $\KK^{\circ}$-algebra 
$\mathcal{R}$,
we can associate an affinoid algebra
$\mathcal{R} \otimes_{\KK^{\circ}} \KK$.
Thus
to an affine admissible formal scheme $\mathscr{U} = \Spf
( \mathcal{R} )$,
one associates an Berkovich affinoid space 
$\mathscr{U}^{\an} = \mathcal{M}(\mathcal{R} \otimes_{\KK^{\circ}} \KK)$.

A formal scheme over $\KK^{\circ}$ is called an 
\emph{admissible formal scheme}
if it has a locally finite open atlas of affine admissible
formal schemes.
Let $\mathscr{X}$ be an admissible formal scheme.
We take an affine covering $\{ \mathscr{U}_\lambda \}$ of $\mathscr{X}$.
Then we have a family $\{ \mathscr{U}_\lambda^{\an} \}$
of Berkovich affinoid spaces.
In fact, 
the patching data of the covering $\{ \mathscr{U}_\lambda \}$
give rise to patching data of the family 
$\{ \mathscr{U}_\lambda^{\an} \}$
of Berkovich affinoid spaces,
and hence those Berkovich affinoid spaces patch together to be
a topological space.
Further, this topological space does not depend on the choice
of the affine covering $\{ \mathscr{U}_\lambda \}$
and depends only on $\mathscr{X}$.
We denote this topological space by $\mathscr{X}^{\an}$
and call it the \emph{Raynaud generic fiber of $\mathscr{X}$}.

The terminology of Berkovich spaces here
is compatible with that in \S~\ref{subsect:berkovich}
when we consider proper schemes.
Let $\mathscr{X}$ be a proper flat scheme over
$\KK^{\circ}$ with irreducible and reduced generic fiber $X$.
Since $X$ is a variety over $\KK$, one associates a Berkovich space
$X^{\an}$ in the sense of \S~\ref{subsect:berkovich}.
On the other hand,
the formal completion 
$\widehat{\mathscr{X}}$
with respect to an element $\varpi \in \KK$ with $0 < | \varpi | < 1$
is an admissible formal scheme,
and hence
 one can associate the Raynaud generic fiber
$\widehat{\mathscr{X}}^{\an}$,
which is a Berkovich analytic space in the above sense.
Then one checks that $X^{\an} = \widehat{\mathscr{X}}^{\an}$ holds.
We remark, however, that
if $\mathscr{X}$ is not proper,
then $X^{\an}$ and $\widehat{\mathscr{X}}^{\an}$ are different.

\begin{Example}
Let $x_1 , \ldots , x_n$ be indeterminates,
and
consider 
$\Spec ( \KK^{\circ} [x_1^{\pm 1} , \ldots , x_n^{\pm 1}] )$.
Then the generic fiber equals the algebraic torus
$\mathbb{G}_m^n = \Spec ( \KK [x_1^{\pm 1} , \ldots , x_n^{\pm 1}] )$,
and we have
its Berkovich analytification
$(\mathbb{G}_m^n)^{\an}$.
On the other hand,
the formal completion of  
$\Spec ( \KK^{\circ} [x_1^{\pm 1} , \ldots , x_n^{\pm 1}] )$
is
$\Spf ( \KK^{\circ} \langle x_1^{\pm 1} , \ldots , x_n^{\pm 1}
\rangle )$,
called the \emph{formal torus},
and one sees that
the Raynaud generic fiber
$
\left( \Spf ( \KK^{\circ} \langle x_1^{\pm 1} , \ldots , x_n^{\pm 1}
\rangle )
\right)^{\an}
=: 
(\mathbb{G}_m^n)_1^{\an}$
equals
$\{ 
p \in (\mathbb{G}_m^n)^{\an} 
\mid
|x_1 (p)| = \cdots = |x_n (p) | = 1 \}$.
\end{Example}

\subsection{Raynaud extensions and valuation maps}
\label{subsection:ap2}

Let $A$ be an abelian variety over $\KK$.
Then there exists a
unique admissible formal group scheme
$\mathscr{A}^{\circ}$ over $\KK^{\circ}$
with a homomorphism
$i : (\mathscr{A}^{\circ})^{\an} \to A^{\an}$
having the following properties.
\begin{itemize}
\item
The image $A^{\circ}$ of $i$ is an analytic subdomain of $A^{\an}$,
and
$i$ is an isomorphism to its image.
\item
There exist a formal abelian scheme $\mathscr{B}$
and an exact sequence
\[
1 \to
\mathscr{T}
\to \mathscr{A}^{\circ}
\to
\mathscr{B}
\to 0
\]
of admissible formal group schemes,
where $\mathscr{T}$ is a formal torus.
\end{itemize}
We remark that the exact sequence in (\ref{align:raynaud:special})
arises by
restricting
the above exact sequence of formal group schemes
to the special fibers.

Taking the Raynaud generic fiber of the above exact sequence
and identifying $\mathscr{T}^{\an}$ with $(\mathbb{G}_m^n)_1^{\an}$,
we obtain an exact sequence
\[
1 \to
(\mathbb{G}_m^n)_1^{\an}
\to (\mathscr{A}^{\circ})^{\an}
\to
\mathscr{B}^{\an}
\to 0
\]
of Berkovich analytic group spaces.
Pushing out this exact sequence by the natural inclusion
$(\mathbb{G}_m^n)_1^{\an} \hookrightarrow (\mathbb{G}_m^n)_1^{\an}$,
we obtain an exact sequence
\[
1 \to
(\mathbb{G}_m^n)^{\an}
\to E
\to
\mathscr{B}^{\an}
\to 0
.
\]
These two exact sequences of analytic groups are called
\emph{Raynaud extension of $A^{\an}$}.
One shows  that the  injective homomorphism 
$(\mathscr{A}^{\circ})^{\an} \to A^{\an}$ extends to a 
unique homomorphism $p : E \to A^{\an}$.
Further, $M:= \Ker (p)$ is a discrete subgroup of $E$,
and $p$ induces an isomorphism $E / M \cong A^{\an}$.
The nonnegative integer $n$ above is called the dimension of the
torus part of $A^{\an}$.
The number $r$ in (\ref{align:raynaud:special})
equals this number $n$
for $A = A_v$.

Recall that we have a homomorphism 
$(\mathbb{G}_m^n)^{\an} \to \RR^n$
(cf. (\ref{align:def:val})).
In fact, one shows that this homomorphism extends to 
a unique continuous homomorphism $\val : E \to \RR^n$,
called the valuation map.
Set $\Lambda : = \val (M)$.
Then one also shows that $\Lambda$ is a complete lattice of $\RR^n$.
(If fact, we denote by $V$ the value group of $\KK$, then $\Lambda \subset
V^n$.)
Further, $\val$ induces a continuous homomorphism
$\overline{\val} : A^{\an} \to \RR^n / \Lambda$.




\small{

}

\end{document}